\documentclass[preprint,3p]{elsarticle}

\usepackage{amssymb}
\usepackage{amsmath}
\usepackage{amsthm}
\usepackage{enumitem}
\usepackage{array}
\usepackage{tikz-cd}
\usepackage[colorlinks=true]{hyperref}

\setlist[enumerate]{font=\upshape}

\theoremstyle{plain}
\newtheorem{thm}{Theorem}[section]
\newtheorem{lem}[thm]{Lemma}
\newtheorem{prop}[thm]{Proposition}
\newtheorem{cor}[thm]{Corollary}
\theoremstyle{definition}
\newtheorem{defn}[thm]{Definition}
\newtheorem{exmp}[thm]{Example}
\newtheorem{rem}[thm]{Remark}

\newcommand{\longhookrightarrow}{\lhook\joinrel\longrightarrow}
\newcommand{\TopMV}{\mathbf{TopMV}}
\newcommand{\Topulg}{\mathbf{Topu\ell G}}
\newcommand{\uLG}{\mathbf{u\ell G}}
\newcommand{\MV}{\mathbf{MV}}
\newcommand{\LCMV}{\mathbf{LCMV}}
\newcommand{\ULS}{\mathbf{ULS\ell G}_{u}}
\newcommand{\TopMVkw}{\mathbf{TopMV}_{k_\omega}}
\newcommand{\Topulgkw}{\mathbf{Topu\ell G}_{k_\omega,\mathrm{cob}}}
\newcommand{\TopulgM}{\mathbf{Topu\ell G}_{\mathrm{M}}}
\newcommand{\QTopMV}{\mathbf{QTopMV}}
\newcommand{\QTopulg}{\mathbf{QTopu\ell G}}
\newcommand{\QTopulgM}{\mathbf{QTopu\ell G}_{\mathrm{qM}}}

\numberwithin{equation}{section}

\begin{document}

\begin{frontmatter}

\title{A topological Mundici representation theorem}

\author{Yifan Meng}
\ead{mengyifanmoe@qq.com}

\author{Lili Shen\corref{cor}}
\ead{shenlili@scu.edu.cn}

\author{Dexue Zhang}
\ead{dxzhang@scu.edu.cn}

\address{School of Mathematics, Sichuan University, Chengdu 610064, China}
\cortext[cor]{Corresponding author.}

\begin{abstract}
We lift Mundici's equivalence between MV-algebras and abelian
$\ell$-groups with a strong unit to the topological setting.
To this end, we prove that, for every topological MV-algebra,
the associated Mundici group can be endowed with a topology that
extends the topology of the MV-algebra and makes it into an abelian
topological $\ell$-group.
Among such topologies there is a finest one;
we identify this topology as the topological-group reflection of
the bounded-stage final topology, and construct an explicit
zero-neighbourhood base in terms of the original MV-algebra topology.
Using this topology, we establish a topological Mundici adjunction
that restricts to an equivalence between topological MV-algebras and
Mundici-fine abelian topological $\ell$-groups with a strong unit.
Here, Mundici-fine means that the group topology is the finest
$\ell$-group extension of its unit-interval topology.
We also identify extremal $\ell$-group extensions under additional
hypotheses, and establish a quasitopological correspondence clarifying
the distinction between separate and joint continuity.
\end{abstract}

\begin{keyword}
%% keywords here, in the form: keyword\sep keyword
Topological MV-algebra \sep Topological lattice-ordered group
\sep Topological Mundici adjunction \sep Topological Mundici correspondence
\sep Quasitopological MV-algebra \sep Quasitopological lattice-ordered group
%% PACS codes here, in the form:\PACS code\sep code

%% MSC codes here, in the form:\MSC code\sep code
%% or\MSC[2008] code\sep code(2000 is the default)
\MSC[2020] 06D35 \sep 06F20 \sep 06F30 \sep 18A40 \sep 22A30
\end{keyword}

\end{frontmatter}

\tableofcontents

\section{Introduction}
\label{sec-introduction}
\emph{MV-algebras} were introduced by Chang to provide an algebraic semantics for
\L ukasiewicz's infinite-valued propositional logic \cite{Chang1958}.
Their prototypical example is the real unit interval $[0,1]$, with
\[
x\oplus y=\min\{x+y,1\}\quad\text{and}\quad \neg x=1-x;
\]
Boolean algebras are precisely the MV-algebras satisfying $x\oplus x=x$.  
Since then, MV-algebras have 
become central objects in many-valued logic, with a rich theory drawing on 
ordered algebra and piecewise-linear geometry. Standard references are the 
monographs \cite{Cignoli2000a,Mundici2011}.

One of the central structural results in this theory is 
Mundici's theorem \cite{Mundici1986},
which identifies MV-algebras with unit intervals in certain
\emph{lattice-ordered groups} (\emph{$\ell$-groups} for short), namely
abelian ones equipped with a strong unit.  For background on
lattice-ordered groups, see
\cite{Birkhoff1942,Birkhoff1967,Darnel1994,Kopytov1994}.  If $(G,u)$ is such a group,
then
\begin{equation} \label{eq-Gamma-def}
 \Gamma(G,u)=[0,u],\qquad
 x\oplus y=(x+y)\wedge u,\qquad
 \neg x=u-x
\end{equation}
is an MV-algebra, and every MV-algebra arises in this way.  The resulting
functor
\[
 \Gamma\colon\uLG\longrightarrow\MV
\]
is an equivalence of categories \cite{Mundici1986,Cignoli1998}.  We write
$\Xi\colon\MV\longrightarrow\uLG$ for its classical quasi-inverse, so that
\[
\Xi A=(G_A,u_A)
\]
is the \emph{Mundici group} of $A$.  This correspondence originated in the 
study of \L ukasiewicz logic and
$K$-theoretic invariants of AF $C^*$-algebras.
Subsequent developments include connections with polyhedral geometry
\cite{Mundici1988}, alternative proofs and constructive treatments
\cite{Dubuc2015,MisselbeckWessel2026}, algebraic variants and generalizations
\cite{Abbadini2021,DiNola1994,Dvurecenskij2002,Kroupa2021},
and extensions to topos-theoretic settings \cite{Caramello2015}.

The question studied here is whether this classical representation
survives the addition of compatible topologies.  A \emph{topological
MV-algebra} \cite{Hoo1997} carries a topology making its MV-operations 
continuous, while a \emph{topological unital $\ell$-group} 
\cite{Smarda1967,Ball1979} carries a topology compatible with both
its group and lattice operations.  The unit interval of the latter is
automatically a topological MV-algebra.  In the reverse direction,
one must \emph{reconstruct} a compatible topology on the whole Mundici
group from the given topology on the unit interval.  The algebraic
Mundici construction involves finite sums and differences of interval
elements, with no fixed bound on the number of terms.  Continuity of
the MV-operations does not by itself give continuity of the group
operations.

Research on topological lattice-ordered groups has addressed
the construction of compatible topologies from the order structure
\cite{Kalton1974,Redfield1974,Yang2009a} and from group-valued
metrics \cite{Kopperman2011}, embeddings \cite{Madell1969},
lattices of compatible topologies \cite{Smarda1976},
completeness and completion \cite{Jordan2018},
and connectedness \cite{Jordan2023}.
For topological MV-algebras, studies have explored connections with
measure theory \cite{Weber2011}, locally convex and uniform structures
\cite{Weber2012,Li2026,Li2027}, relations with locally solid $\ell$-groups
\cite{Weber2023}, initial topologies and topological isomorphism
theorems \cite{Yin2025}, and regularity and separation properties
\cite{Gan2024}.
Compact Hausdorff and Stone MV-algebras are treated in
\cite{Nganou2017,Nganou2019}.
At a broader categorical level, topological theories and enriched
Stone-type dualities provide a common framework for interactions of
this kind between topology, order, and algebra
\cite{Hofmann2007,Hofmann2018}.
MV-algebras and unital $\ell$-groups also enter algebraic descriptions
of compact Hausdorff spaces and extensions of classical dualities
\cite{Cignoli2004,Marra2017,Abbadini2025};
these duality results concern special classes of locally finite,
complete, or compact objects.
The question we address is different: does every topological MV-algebra
admit an $\ell$-group topology on its Mundici group which induces the
prescribed topology on the unit interval?

We prove that the answer is affirmative:
\begin{itemize}[beginpenalty=10000]
\item Every topological MV-algebra $(A,\tau)$ admits an $\ell$-group
topology on its Mundici group $(G_A,u_A)$ inducing $\tau$ on the unit
interval $[0,u_A]$.
\item Among such topologies there is a finest one, denoted by $\sigma_A$,
for which we give an explicit zero-neighbourhood base in terms of $\tau$
(Corollary~\ref{cor-reconstruction-topology-chain}).
\item The assignment $(A,\tau)\mapsto(G_A,u_A,\sigma_A)$ extends to a
functor that yields a topological Mundici adjunction and an equivalence
with the full subcategory of Mundici-fine groups, whose topology is the
finest $\ell$-group extension of their unit-interval topology.
\end{itemize}

The paper is organized as follows.
In Section~\ref{sec-adjunction}, the topological unit-interval functor
\[
 \Gamma_t\colon\Topulg\longrightarrow\TopMV
\]
is shown to have a left adjoint
\[
 \Xi_t\colon\TopMV\longrightarrow\Topulg,
 \qquad \Xi_t\dashv\Gamma_t.
\]
This \emph{topological Mundici adjunction}
(Theorem~\ref{thm-topological-adjunction}) is constructed independently
of the interval-preserving reconstruction.
Its unit and counit have the same underlying algebraic maps as in the
classical Mundici equivalence.
At this stage, however, each component of the unit is known only to
be a continuous algebraic isomorphism; continuity of its inverse is
not yet available.

To phrase the question differently, we identify $A$ with $[0,u_A]$.
As in Definition~\ref{def-admissible-topology}, an $\ell$-group topology
$\rho$ on $G_A$ is called
\begin{itemize}
\item \emph{admissible} for $(A,\tau)$ if
$\rho|_{[0,u_A]}\subseteq\tau$;
\item \emph{exact} for $(A,\tau)$ if
$\rho|_{[0,u_A]}=\tau$.
In this case, $\rho$ is called an \emph{$\ell$-group extension of $\tau$}.
\end{itemize}
Every exact topology is admissible, but the converse need not hold.
The functor $\Xi_t$ sends $(A,\tau)$ to $(G_A,u_A,\sigma_A)$, where
$\sigma_A$ is the largest admissible $\ell$-group topology; we call it
the \emph{Mundici topology} associated with $(A,\tau)$.
Thus the adjunction isolates the reconstruction problem as the question
whether $\sigma_A$ is exact:
\[
 \sigma_A|_{[0,u_A]}=\tau.
\]
By maximality, it suffices to find any $\ell$-group extension $\rho$
of $\tau$: admissibility gives $\rho\subseteq\sigma_A$, and hence
forces the same interval equality for $\sigma_A$.
Conversely, exactness of $\sigma_A$ makes it such an extension.

Section~\ref{sec-coordinate-bounded-reconstruction} first uses the
unit-interval topology $\tau$ to define a topology $\tau_n$ on each
bounded stage
\[
 I_n=[-nu_A,nu_A]\qquad(n\geq1).
\]
Proposition~\ref{prop-lgroup-extension-bounded-stages} establishes
\emph{bounded-stage rigidity}: every $\ell$-group extension $\rho$
of $\tau$ satisfies $\rho|_{I_n}=\tau_n$ for every $n\geq1$.
These intervals cover $G_A$, but two $\ell$-group extensions can
agree on every $I_n$ and still differ on the whole group
(Example~\ref{exmp-product-cubes}).
We must therefore find a topology on $G_A$ that induces $\tau_n$
on every $I_n$ and makes addition and the lattice operations jointly
continuous.  This section constructs a lower bound $\kappa_A$ and an
upper bound $\lambda_A$: every $\ell$-group extension $\rho$ of $\tau$
satisfies
\[
 \kappa_A\subseteq\rho\subseteq\lambda_A.
\]
The two bounds are constructed as follows.
\begin{itemize}
\item The \emph{coordinate topology} $\kappa_A$ is the initial topology
on $G_A$ induced by all bounded coordinate maps into $(A,\tau)$.
On each $I_n$, $\kappa_A$ induces $\tau_n$
(Proposition~\ref{prop-kappa-lambda}).
\item The \emph{bounded-stage final topology} $\lambda_A$ declares a
subset of $G_A$ open precisely when its intersection with each $I_n$
is open in $\tau_n$.
It is the finest topology on $G_A$ agreeing with $\kappa_A$ on every $I_n$
(Proposition~\ref{prop-lambda-exact-space}).
\end{itemize}
Both recover the unit-interval topology, but they can differ globally
(Example~\ref{exmp-kappa-lambda-strict}); neither need be a group
topology (Examples~\ref{exmp-kappa-not-group}
and~\ref{exmp-lambda-not-group}).
The section also characterizes $\sigma_A$ as the finest $\ell$-group
topology contained in $\lambda_A$
(Theorem~\ref{thm-lambda-envelope}).

Nevertheless, $\lambda_A$ makes the Mundici group a quasitopological
$\ell$-group: inversion is continuous, while addition and the lattice
operations are separately continuous
(Propositions~\ref{prop-lambda-quasitopological}
and~\ref{prop-lambda-quasi-reconstruction}).
Even though addition is jointly continuous on products of bounded
stages, joint continuity on the whole group $(G_A,\lambda_A)$ does
not follow, because products need not preserve the final topology.
In fact, joint continuity fails for the rational product cube in
Example~\ref{exmp-lambda-not-group}.

Section~\ref{sec-group-envelope} passes to the topological-group
reflection $\gamma_A$ of $\lambda_A$, that is, the finest group
topology contained in $\lambda_A$.
Addition is now jointly continuous, but we must still prove that
$\gamma_A$ induces the prescribed topology on every $I_n$ and makes
the lattice operations continuous.
We give an explicit zero-neighbourhood base for $\gamma_A$ consisting
of unions of finite sums of symmetric coordinate boxes
(Proposition~\ref{prop-coordinate-box-topology}).

To use the topology on $A$, we write elements of these neighbourhoods
as finite signed sums of elements of $A$.
There is no fixed bound on the number of summands, so continuity on
a fixed bounded stage does not suffice.
The key technical step is the \emph{bounded-output lemma}
(Lemma~\ref{lem-strong-FBOT}):
\begin{itemize}
\item Fix a finite $\ell$-group term $t(x_1,\ldots,x_r)$, allowing
fixed parameters from $G_A$, with $t(0,\ldots,0)=0$.
For each stage $I_m$ and each zero-neighbourhood $V$ in $I_m$,
there is a sequence $(W_q)_{q\geq1}$ of zero-neighbourhoods in $A$
such that, for any pairwise disjoint finite sets $E_1,\ldots,E_r$
of positive integers, the sums
\[
 g_i=\sum_{q\in E_i}\epsilon_q a_q\quad(1\leq i\leq r),
 \qquad a_q\in W_q,\quad \epsilon_q\in\{-1,1\}
\]
satisfy
\[
 t(g_1,\ldots,g_r)\in I_m
 \implies
 t(g_1,\ldots,g_r)\in V.
\]
The same sequence works for all these choices; the input sums
need not lie in a common bounded stage.
\end{itemize}

Together with the coordinate-box neighbourhoods, this estimate
proves continuity of the bounded coordinates: their changes always
lie in $I_1$.  It follows that $\gamma_A$ retains the topology on
each $I_n$.
We also use this estimate to prove joint continuity of the
lattice operations (Theorem~\ref{thm-rho-lgroup}).
Thus $\gamma_A$ is an $\ell$-group extension of $\tau$.
By maximality of $\sigma_A$, we have $\gamma_A\subseteq\sigma_A$,
and hence
\begin{equation}
\label{eq-intro-topological-exactness}
 \tau=\gamma_A|_A\subseteq\sigma_A|_A\subseteq\tau,
\end{equation}
where the last inclusion follows from admissibility of $\sigma_A$.
In fact, $\gamma_A=\sigma_A$, since $\sigma_A$ is a group
topology contained in $\lambda_A$.
The resulting \emph{topological Mundici representation theorem}
(Theorem~\ref{thm-topological-representation}) states:
\begin{itemize}
\item For every topological MV-algebra $(A,\tau)$, the canonical
Mundici isomorphism
\[
 \eta_A\colon(A,\tau)\longrightarrow\Gamma_t\Xi_t(A,\tau)
\]
is a homeomorphism.  In particular, $\sigma_A|_{[0,u_A]}=\tau$,
so $(A,\tau)$ arises as the unit interval of the topological
unital abelian $\ell$-group $(G_A,u_A,\sigma_A)$.
\end{itemize}
The accompanying Corollary~\ref{cor-reconstruction-topology-chain}
records the maximality and group-reflection description of $\sigma_A$,
together with its explicit zero-neighbourhood base.

Since the unit is a natural isomorphism, $\Xi_t$ is fully faithful.
In the other direction, reconstructing a topological unital
$\ell$-group from its unit interval gives the same algebraic group
with a possibly finer topology.
Section~\ref{sec-correspondence} shows that the counit is a
homeomorphism precisely when the original topology is already the
finest extension of its unit-interval topology
(Proposition~\ref{prop-transported-maximality}).
Together, these results give the \emph{topological Mundici correspondence}
(Theorem~\ref{thm-topological-correspondence}):
\begin{itemize}
\item The essential image of $\Xi_t$ is the replete coreflective full
subcategory $\TopulgM$ of \emph{Mundici-fine} topological unital
abelian $\ell$-groups, whose topology is the finest $\ell$-group
extension of their given unit-interval topology.
Hence
\[
 \TopMV\simeq\TopulgM.
\]
\end{itemize}
The choice of this subcategory is essential.
The same unit-interval and bounded-stage topologies can occur in
nonisomorphic topological unital abelian $\ell$-groups
(Example~\ref{exmp-product-cubes}); the correspondence selects the
largest $\ell$-group extension functorially.
Together, Theorems~\ref{thm-topological-representation} and
\ref{thm-topological-correspondence} lift the classical Mundici
equivalence to the topological setting, as shown in the
commutative diagram~\eqref{eq-topological-equivalence-square}.

As a corollary of the representation theorem, we show at the end of
Section~\ref{sec-correspondence} that every
topological MV-algebra $(A,\tau)$ is a topological retract of its
Mundici group $(G_A,\sigma_A)$
(Corollary~\ref{thm-retract-Tychonoff}), and derive separation-axiom
consequences (Corollary~\ref{cor-separation-axioms}).
Related separation and regularity results for topological MV-algebras
appear in \cite{Hoo1997,Gan2024,Li2026}.

Section~\ref{sec-restricted-correspondence} studies two settings in
which the upper and lower bounds from
Section~\ref{sec-coordinate-bounded-reconstruction} become extremal
$\ell$-group extensions.  For a
$k_\omega$ topological MV-algebra,
Theorem~\ref{thm-komega-reconstruction} proves that the bounded-stage
final topology $\lambda_A$ is already an $\ell$-group topology,
hence $\lambda_A=\sigma_A$. 
Theorem~\ref{thm-komega-equivalence} then gives an equivalence
\[
 \TopMVkw\simeq\Topulgkw.
\]
The group-side category consists of $k_\omega$ topological unital
abelian $\ell$-groups whose compact subsets are order bounded.
For compact Hausdorff $A$, the bounded stages form a
$k_\omega$-decomposition of the reconstructed group
(Corollary~\ref{cor-compact-reconstruction}).

Local convexity gives the complementary result at the lower bound.
Theorem~\ref{thm-locally-convex-reconstruction} proves that $\kappa_A$
is an unbounded locally solid $\ell$-group topology extending $\tau$.
Hence, it is the least $\ell$-group extension of $\tau$
(Corollary~\ref{cor-extremal-lgroup-extensions}).
This leads to an alternative realization
(Theorem~\ref{thm-locally-convex-equivalence})
\[
 \LCMV\simeq\ULS
\]
through $\Xi_u(A,\tau)=(G_A,u_A,\kappa_A)$.
Unlike the topological Mundici correspondence, this equivalence
uses the least $\ell$-group extension $\kappa_A$ of $\tau$.

Section~\ref{sec-quasi-reflection} revisits the reconstruction
in a quasitopological setting.  Binary operations are required to be
separately continuous, while negation and inversion are continuous.
The initial topology $\kappa_A$ already gives a quasitopological
extension of the interval topology, making this reconstruction
more direct:

\begin{itemize}
\item For a quasitopological MV-algebra $(A,\tau)$, the initial topology
$\kappa_A$ is the least quasitopological $\ell$-group topology
inducing $\tau$ on $A$
(Proposition~\ref{thm-quasi-minimal-reconstruction}).
The translated unit intervals define a final topology
$\sigma_A^{\mathrm q}$, which is the largest quasi-admissible topology
(Proposition~\ref{prop-quasi-translated-final}).
Since $\kappa_A$ is quasi-admissible,
\begin{equation}
\label{eq-intro-quasi-exactness}
 \tau=\kappa_A|_A
 \subseteq(\sigma_A^{\mathrm q})|_A\subseteq\tau.
\end{equation}
The last inclusion is the admissibility condition.
Thus $\sigma_A^{\mathrm q}$, the topology used by the left adjoint
$\Xi_{\mathrm q}$, is exact, yielding the quasi correspondence
(Proposition~\ref{prop-quasi-adjunction} and
Corollary~\ref{thm-quasi-correspondence}).

\item For a topological MV-algebra $(A,\tau)$, however,
$\kappa_A$ need not be an $\ell$-group topology
(Example~\ref{exmp-kappa-not-group}),
so the same admissibility argument does not apply automatically.
Exactness of $\sigma_A$ instead follows from
\eqref{eq-intro-topological-exactness}, the analogue of
\eqref{eq-intro-quasi-exactness} with the group reflection
$\gamma_A$ in place of $\kappa_A$.
The difficult step, carried out in Section~\ref{sec-group-envelope},
is to prove that $\gamma_A$ is an exact $\ell$-group topology.
\end{itemize}

For topological inputs, the topological $\ell$-group reflection
of the maximal quasi reconstruction $(G_A,u_A,\sigma_A^{\mathrm q})$
is $(G_A,u_A,\sigma_A)$
(Proposition~\ref{thm-quasi-reflection-factorization}).
Thus this reflection recovers the topological Mundici representation
and preserves the given topology $\tau$ on the unit interval
(Corollary~\ref{cor-quasi-reflection-exactness}).

Section~\ref{sec-quasi-reflection} also identifies the box-subspace topology
$\nu_A$ induced by the bounded coordinates or their deviations
(Proposition~\ref{prop-box-quasi-reconstruction}).
It is an upper bound for all $\ell$-group extensions of $\tau$ and
is contained in $\lambda_A$.
Both $\nu_A$ and $\lambda_A$ have group reflection $\sigma_A$,
but $\nu_A$ need not itself be a group topology
(Example~\ref{exmp-box-not-group}).
Example~\ref{exmp-all-reconstruction-inclusions-strict} shows that
all four inclusions in the comparison
\eqref{eq-quasi-full-chain} can be strict simultaneously.
Finally, Section~\ref{sec-summary} summarizes the resulting topology
comparisons, the roles of the topologies in the proof, and
related examples.

\section{The Mundici topology \texorpdfstring{$\sigma_A$}{sigma\_A}
and the topological Mundici adjunction}
\label{sec-adjunction}

We begin by fixing the algebraic and topological categories used throughout
the paper.  We then construct the topological Mundici adjunction directly
from a maximality principle.  The point is that the adjunction itself does
not require the bounded-stage reconstruction developed later: one can first
form the finest $\ell$-group topology making the canonical map from
$(A,\tau)$ continuous and only afterwards ask whether this topology recovers
$\tau$ on the unit interval.

First recall the two categories underlying the algebraic Mundici
equivalence.  An \emph{MV-algebra} \cite{Chang1958,Cignoli2000a,Mundici2011} 
is an algebra $(A,\oplus,\neg,0)$ such that $(A,\oplus,0)$ is 
a commutative monoid and, for all $x,y\in A$,
\[
 \neg\neg x=x,\qquad x\oplus\neg0=\neg0,
 \qquad
 \neg(\neg x\oplus y)\oplus y
 =\neg(\neg y\oplus x)\oplus x.
\]
Put $1=\neg0$.  The canonical order on $A$ is given by
\[
 x\leq y\iff \neg x\oplus y=1;
\]
it is a lattice order.  An \emph{MV-homomorphism} is a map preserving
$\oplus$, $\neg$, and $0$.  We denote the category of MV-algebras and
MV-homomorphisms by 
\[
\MV.
\] 

An \emph{$\ell$-group} (short for \emph{lattice-ordered group})
\cite{Birkhoff1942,Birkhoff1967,Darnel1994,Kopytov1994} is a group equipped with a
lattice order such that all left and right translations are order
automorphisms.  All $\ell$-groups considered below are abelian.
We write their group operations additively and take empty sums to be $0$.
Thus, in such an $\ell$-group $G$, 
every translation
$x\mapsto a+x$ is an order automorphism; equivalently, addition
distributes over finite joins and meets.  We shall also use repeatedly
the standard fact that the lattice underlying an $\ell$-group is
distributive (see \cite[Section~2.2, Proposition~2]{Kopytov1994}).
Write 
\[
 |x|=x\vee(-x).
\]
A \emph{strong unit} of $G$ is an element $u\geq0$ such
that, for every $x\in G$, there is an integer $n\geq1$ with
$|x|\leq nu$.  A \emph{unital abelian $\ell$-group} is a pair $(G,u)$
consisting of an abelian $\ell$-group and a specified strong unit.  An
$\ell$-group homomorphism preserves the group and lattice operations; a
homomorphism $(G,u)\longrightarrow(H,v)$ is \emph{unit-preserving} if it
sends $u$ to $v$.  We denote the category of unital abelian $\ell$-groups
and unit-preserving $\ell$-group homomorphisms by 
\[
\uLG.
\]

For $(G,u)\in\uLG$, let
\[
 \Gamma(G,u)=[0,u]=\{x\in G:0\leq x\leq u\},
\]
where $[0,u]$ is equipped with the MV-operations displayed in \eqref{eq-Gamma-def}.
For a morphism $h\colon(G,u)\longrightarrow(H,v)$ in $\uLG$, the restriction
$h|_{[0,u]}$ is regarded as a map $[0,u]\longrightarrow[0,v]$, 
since $h$ is order-preserving and $h(u)=v$. 
Thus, the interval construction defines a functor
\begin{equation}
\label{eq-Gamma-definition}
 \Gamma\colon\uLG\longrightarrow\MV,
 \qquad
 \Gamma(G,u)=[0,u],
 \qquad
 \Gamma h=h|_{[0,u]}.
\end{equation}
Mundici's theorem states that the functor $\Gamma$ is an equivalence
of categories \cite{Mundici1986,Cignoli1998}.
Equivalently, we may choose a functor
\[
 \Xi\colon\MV\longrightarrow\uLG
\]
such that $\Xi\dashv\Gamma$ is an adjoint equivalence.  We denote its unit
and counit, both natural isomorphisms, by
\begin{equation}
\label{eq-Mundici-unit-counit}
 \eta\colon1_{\MV}\Longrightarrow\Gamma\Xi,
 \qquad
 \varepsilon\colon\Xi\Gamma\Longrightarrow1_{\uLG}.
\end{equation}
In particular, if $\Xi A=(G_A,u_A)$, then
\begin{equation}
\label{eq-Mundici-unit}
 \eta_A\colon A\longrightarrow[0,u_A]
\end{equation}
is an MV-isomorphism.  We identify $A$ algebraically with this unit
interval whenever convenient.

An \emph{MV-algebra topology} on an MV-algebra is a topology for which
$\oplus$ and $\neg$ are continuous; a \emph{topological MV-algebra} $(A,\tau)$ is an
MV-algebra $A$ equipped with such a topology $\tau$ \cite{Hoo1997}.  We write 
\[
\TopMV
\] 
for their
category with continuous MV-homomorphisms.  An \emph{$\ell$-group topology}
on an $\ell$-group is a topology for which $+$, unary $-$, $\vee$, and
$\wedge$ are continuous \cite{Smarda1967,Ball1979}.  A \emph{topological unital 
abelian $\ell$-group} $(G,u,\rho)$ is a unital abelian $\ell$-group $(G,u)$ equipped 
with an $\ell$-group topology $\rho$.
Their category, with continuous unit-preserving $\ell$-homomorphisms, is
denoted by 
\[
\Topulg.
\]
Unless explicitly stated otherwise, no separation
axiom is imposed.

\begin{exmp}
\label{exmp-function-algebras}
Let $X$ be a compact Hausdorff space.
\begin{enumerate}[label=(\alph*),ref=(\alph*)]
\item Let $C(X,\mathbb R)$ denote the abelian group of continuous
real-valued functions on $X$, and let $\rho_\infty$ be the topology
induced by the supremum norm
\[
 \|f\|_\infty=\sup_{x\in X}|f(x)|.
\]
With the pointwise lattice order, addition, negation, and the lattice
operations are continuous with respect to $\rho_\infty$.  Moreover, since
every continuous real-valued function on $X$ is bounded, the constant
function $\mathbf 1$ with value $1$ is a strong unit.  Hence
\[
 \bigl(C(X,\mathbb R),\mathbf 1,\rho_\infty\bigr)
\]
is a topological unital abelian $\ell$-group.

\item Let $C(X,[0,1])$ denote the continuous $[0,1]$-valued functions on
$X$, and let
\[
 \tau_\infty=\rho_\infty|_{C(X,[0,1])}
\]
be the topology induced by the supremum norm, equivalently the subspace
topology inherited from $C(X,\mathbb R)$.  Equip $C(X,[0,1])$ with the
pointwise MV-operations
\[
  (f\oplus g)(x)=\min\{f(x)+g(x),1\},\qquad
  (\neg f)(x)=1-f(x).
\]
These operations are continuous, so
\[
 \bigl(C(X,[0,1]),\tau_\infty\bigr)
\]
is a topological MV-algebra.
\end{enumerate}

If $X$ is an $n$-point discrete space, the two examples reduce to
\[
 (\mathbb R^n,\mathbf 1,\rho_\infty)
 \qquad\text{and}\qquad
 ([0,1]^n,\tau_\infty),
\]
respectively, where the displayed topologies are the usual Euclidean
topologies.  As we shall see immediately below, the topological
unit-interval functor relates the two examples by
\[
 \Gamma_t\bigl(C(X,\mathbb R),\mathbf 1,\rho_\infty\bigr)
 =
 \bigl(C(X,[0,1]),\tau_\infty\bigr).
\]
\end{exmp}

For $(G,u,\rho)\in\Topulg$, the subspace $[0,u]$ is a topological
MV-algebra, since the operations displayed in \eqref{eq-Gamma-def} are
restrictions of continuous maps on $G$.  Moreover, every continuous
unit-preserving $\ell$-homomorphism restricts to a continuous
MV-homomorphism between the corresponding unit intervals.  Thus the
algebraic interval functor lifts to
\begin{equation}
\label{eq-Gammat-functor}
 \Gamma_t\colon\Topulg\longrightarrow\TopMV,
\end{equation}
with object and morphism assignments
\begin{equation}
\label{eq-Gammat-definition}
 \Gamma_t(G,u,\rho)=([0,u],\rho|_{[0,u]}),
 \qquad
 \Gamma_t h=h|_{[0,u]}.
\end{equation}

The interval construction therefore presents no topological difficulty in
this direction.  Unlike its algebraic counterpart, however, $\Gamma_t$
is not full and can send nonisomorphic objects to the same topological
MV-algebra (Example~\ref{exmp-product-cubes}).  Hence it is not an
equivalence between $\Topulg$ and $\TopMV$.  Nevertheless, it
is still a right adjoint.  To construct its left adjoint, we equip the
Mundici group with the largest $\ell$-group topology for which the canonical
inclusion of the original MV-algebra is continuous.  Maximality will force
continuous MV-homomorphisms to extend continuously.

Let $(A,\tau)\in\TopMV$ and write $\Xi A=(G_A,u_A)$.
Regard the classical unit from \eqref{eq-Mundici-unit} also as a map into
the ambient group:
\begin{equation}
\label{eq-canonical-Mundici-inclusion}
 \iota_A\colon
 A\xrightarrow{\eta_A}[0,u_A]\longhookrightarrow G_A.
\end{equation}

\begin{defn}
\label{def-admissible-topology}
Let $(A,\tau)\in\TopMV$ and write $\Xi A=(G_A,u_A)$.
An $\ell$-group topology $\rho$ on $G_A$ is called
\begin{itemize}
\item \emph{admissible for $(A,\tau)$} if the canonical map $\iota_A$
from \eqref{eq-canonical-Mundici-inclusion} is continuous;
equivalently, under the algebraic identification $A=[0,u_A]$,
\[
 \rho|_A\subseteq\tau;
\]
\item \emph{exact for $(A,\tau)$} if $\iota_A$ is a topological embedding;
equivalently,
\[
 \rho|_A=\tau.
\]
\end{itemize}
In the exact case, $\rho$ is called an
\emph{$\ell$-group extension of $\tau$}.
\end{defn}

Let $\mathcal T_A$ denote the set of all admissible
$\ell$-group topologies on $G_A$.  Since no separation
axiom is imposed, the indiscrete topology is admissible;
hence $\mathcal T_A\neq\varnothing$.

\begin{lem}
\label{lem-finest-admissible}
The join of any family of $\ell$-group topologies on an
abelian $\ell$-group is an $\ell$-group topology.
Moreover, the join of any family of admissible $\ell$-group
topologies for $(A,\tau)$ is admissible.
In particular, $\mathcal T_A$ has a largest element, namely
\begin{equation}
\label{eq-largest-admissible-topology}
 \sigma_A=\bigvee\mathcal T_A=\bigcup\mathcal T_A .
\end{equation}
\end{lem}

\begin{proof}
Let $(\rho_i)_{i\in I}$ be a family of $\ell$-group topologies
on $G$, and let $\rho=\bigvee_{i\in I}\rho_i$ be their join in the
lattice of topologies on the underlying set. Equivalently, $\rho$
is the topology generated by $\bigcup_{i\in I}\rho_i$.

Assume first that $I\neq\varnothing$. For a basic $m$-ary operation
$\omega$ and any subbasic $\rho$-open set $U$, there exists
$i\in I$ such that $U\in\rho_i$. Continuity in $\rho_i$ gives
\[
 \omega^{-1}(U)\in\rho_i^m\subseteq\rho^m .
\]
By the subbasis criterion, $\rho$ is therefore an
$\ell$-group topology. The case $I=\varnothing$ gives the
indiscrete topology, which is also an $\ell$-group topology.

Now take $G=G_A$.  If each $\rho_i$ is admissible for $(A,\tau)$, then
$\iota_A^{-1}(U)$ is $\tau$-open for every such subbasic
set $U$, so $\iota_A$ is continuous into $\rho$.
Thus the generated topology $\rho$ is admissible.

Applying this to the family $\mathcal T_A$, the topology
$\bigvee\mathcal T_A$ is itself an admissible $\ell$-group topology
for $(A,\tau)$ and contains every member of $\mathcal T_A$. Hence
$\bigvee\mathcal T_A\in\mathcal T_A$.
Since $\bigvee\mathcal T_A$ contains every member of $\mathcal T_A$
and belongs to $\mathcal T_A$, we obtain
$\bigcup\mathcal T_A
\subseteq
\bigvee\mathcal T_A
\subseteq
\bigcup\mathcal T_A$,
which proves \eqref{eq-largest-admissible-topology}.
\end{proof}

We call $\sigma_A$ the \emph{Mundici topology} associated with
$(A,\tau)$.  By Lemma~\ref{lem-finest-admissible}, it is the largest
admissible $\ell$-group topology on $G_A$; at this point only
$\sigma_A|_A\subseteq\tau$ is known, and equality remains to be proved.

The following elementary fact is the mechanism behind both functoriality
and the adjunction.  Pulling a target topology back---that is, taking the
initial topology along an algebraic homomorphism---produces a compatible
source topology; maximality can then be used to prove continuity from
$\sigma_A$.

\begin{lem}
\label{lem-pullback-topology}
Let $h\colon(H,v)\longrightarrow(G,u)$ be a unit-preserving
$\ell$-group homomorphism, and let $\rho$ be an $\ell$-group topology on
$G$.  Then the initial (inverse-image) topology $h^*\rho$ on $H$ induced by
$h$ is an $\ell$-group topology.
\end{lem}

\begin{proof}
As above, write $\rho^m$ and $(h^*\rho)^m$ for the product topologies on
$G^m$ and $H^m$, respectively.  By definition,
$h\colon(H,h^*\rho)\longrightarrow(G,\rho)$ is continuous, so the product
map
\[
 h^m\colon(H^m,(h^*\rho)^m)\longrightarrow(G^m,\rho^m)
\]
is continuous.
For every basic $m$-ary $\ell$-group operation $\omega$,
\[
 h\omega_H=\omega_Gh^m.
\]
For each $U\in\rho$, continuity of $\omega_G$ and $h^m$ gives
\[
 \omega_H^{-1}h^{-1}(U)
 =(h^m)^{-1}\omega_G^{-1}(U)
 \in(h^*\rho)^m.
\]
Since every $h^*\rho$-open set is of the form $h^{-1}(U)$ for some
$U\in\rho$, this proves continuity of $\omega_H$.
\end{proof}

Define $\Xi_t$ on objects and morphisms by
\begin{equation}
\label{eq-Xit-definition}
 \Xi_t(A,\tau)=(G_A,u_A,\sigma_A),\qquad \Xi_t f=\Xi f.
\end{equation}

\begin{prop}
\label{prop-Xit-functor}
The assignments in \eqref{eq-Xit-definition} lift the algebraic functor $\Xi$ to a functor
\[
 \Xi_t\colon\TopMV\longrightarrow\Topulg.
\]
\end{prop}

\begin{proof}
The object assignment is well defined by
Lemma~\ref{lem-finest-admissible}.  Let
$f\colon(A,\tau_A)\longrightarrow(B,\tau_B)$ be a continuous
MV-homomorphism, and put
\[
 h=\Xi f\colon(G_A,u_A)\longrightarrow(G_B,u_B).
\]
Lemma~\ref{lem-pullback-topology} makes $h^*\sigma_B$ an
$\ell$-group topology on $G_A$.  Naturality of the algebraic unit gives
\[
 (\Gamma h)\eta_A=\eta_Bf,
 \qquad\text{hence}\qquad h\iota_A=\iota_Bf.
\]
The right-hand side is continuous into $(G_B,\sigma_B)$.  By the
definition of the initial topology, this means that $\iota_A$ is
continuous into $(G_A,h^*\sigma_B)$.  Thus $h^*\sigma_B$ is admissible
for $(A,\tau_A)$, and maximality gives
\[
 h^*\sigma_B\subseteq\sigma_A.
\]
In particular, $h^{-1}(U)$ is $\sigma_A$-open for every
$\sigma_B$-open $U$, which proves continuity of $h$.
Preservation of identities and composition follows from the same
properties of $\Xi$.
\end{proof}

The same pullback argument applies to a homomorphism into any
topological unital abelian $\ell$-group, whose topology need not itself
be maximal.  This gives the following universal property.

\begin{thm}
\label{thm-topological-adjunction}
The functors $\Xi_t\colon\TopMV\longrightarrow\Topulg$ and
$\Gamma_t\colon\Topulg\longrightarrow\TopMV$ form an adjunction
\[
 \Xi_t\dashv\Gamma_t.
\]
Explicitly, for $(A,\tau)\in\TopMV$ and $(G,u,\rho)\in\Topulg$, the classical
Mundici correspondence restricts to a natural bijection
\begin{equation}
\label{eq-topological-Mundici-bijection}
 \Topulg(\Xi_t(A,\tau),(G,u,\rho))
 \cong
 \TopMV((A,\tau),\Gamma_t(G,u,\rho)),
 \qquad h\mapsto(\Gamma_t h)\eta_A.
\end{equation}
The unit and counit have the same underlying maps as the classical
Mundici unit and counit in \eqref{eq-Mundici-unit-counit}.
\end{thm}

\begin{proof}
The classical Mundici adjunction gives a natural bijection
\begin{equation}
\label{eq-algebraic-Mundici-bijection}
 \uLG(\Xi A,(G,u))
 \cong
 \MV(A,\Gamma(G,u)),
 \qquad
 h\mapsto f=(\Gamma h)\eta_A.
\end{equation}
We show that it restricts to continuous morphisms.  Let
$\rho$ be an $\ell$-group topology on $G$, and put
$h^*\rho$ for the initial topology on $G_A$ induced by
the algebraic $\ell$-homomorphism $h$.
By Lemma~\ref{lem-pullback-topology}, $h^*\rho$ is
an $\ell$-group topology.  Since $\sigma_A$ is the largest
admissible topology and admissibility is preserved upon
passing to coarser $\ell$-group topologies, we have
\[
\begin{aligned}
 h\colon(G_A,\sigma_A)\longrightarrow(G,\rho)
   \text{ is continuous}
 &\iff
   h^*\rho\subseteq\sigma_A\\
 &\iff
   h^*\rho\text{ is admissible}\\
 &\iff
   h\iota_A\colon(A,\tau)\longrightarrow(G,\rho)
   \text{ is continuous}\\
 &\iff
   f\colon(A,\tau)\longrightarrow([0,u],\rho|_{[0,u]})
   \text{ is continuous}.
\end{aligned}
\]
For the last equivalence, use
$h\iota_A=j_Gf$, where
$j_G\colon[0,u]\longhookrightarrow G$ is the subspace inclusion.
This proves the claimed bijection of topological hom-sets.
Naturality and the triangle identities follow from the classical
adjunction, since the unit and counit retain their underlying
algebraic maps.
\end{proof}

Theorem~\ref{thm-topological-adjunction} is therefore a topological lifting
of the classical Mundici adjoint equivalence.  Let
\[
 U_{\mathrm{MV}}\colon\TopMV\longrightarrow\MV,
 \qquad
 U_{\ell}\colon\Topulg\longrightarrow\uLG
\]
denote the forgetful functors.  Then
\begin{equation}
\label{eq-topological-lifting-square}
\begin{tikzcd}[column sep=4.5em,row sep=3em]
\TopMV
  \arrow[r,shift left=1.4ex,"\Xi_t"]
  \arrow[r,phantom,"\bot" description]
  \arrow[d,"U_{\mathrm{MV}}"']
&
\Topulg
  \arrow[l,shift left=1.4ex,"\Gamma_t"]
  \arrow[d,"U_{\ell}"]
\\
\MV
  \arrow[r,shift left=1.4ex,"\Xi"]
  \arrow[r,phantom,"\simeq" description]
&
\uLG
  \arrow[l,shift left=1.4ex,"\Gamma"]
\end{tikzcd}
\end{equation}
commutes.  The upper pair is the adjunction
$\Xi_t\dashv\Gamma_t$, while the lower pair is the classical Mundici
equivalence, with quasi-inverse functors $\Xi$ and $\Gamma$, fixed in
\eqref{eq-Mundici-unit-counit}.

\begin{rem}
\label{rem-taut-lift}
The vertical forgetful functors in
\eqref{eq-topological-lifting-square} are both \emph{topological functors} 
\cite{Adamek1990,Hofmann2014}.
Indeed, for either of the two vertical functors, write $U$ for that
functor and consider a source of algebraic homomorphisms
\[
 \bigl(f_i\colon A\longrightarrow U B_i\bigr)_{i\in I}
\]
into topological objects of the corresponding kind.
Equip $A$ with the initial topology induced by the family
$(f_i)_{i\in I}$.  For every basic $m$-ary
operation $\omega$, one has
\[
 f_i\omega_A=\omega_{B_i}f_i^m.
\]
Since $f_i^m$ is continuous for the corresponding product topologies and
$\omega_{B_i}$ is continuous, the initial property shows that
$\omega_A$ is continuous.  Thus the initial lift again belongs to
$\TopMV$ or $\Topulg$, respectively.

For $U_{\ell}$, Lemma~\ref{lem-pullback-topology} is precisely the
special case of this argument for a source consisting of a single
unit-preserving $\ell$-homomorphism.  The same one-map initial-lift
mechanism is used in Proposition~\ref{prop-Xit-functor} and in the proof
of Theorem~\ref{thm-topological-adjunction} when the target topology is
pulled back along an algebraic extension.

In particular, the two vertical functors have the familiar discrete and
indiscrete adjoints:
\[
 D_{\mathrm{MV}}\dashv U_{\mathrm{MV}}\dashv I_{\mathrm{MV}},
 \qquad
 D_{\ell}\dashv U_{\ell}\dashv I_{\ell}.
\]
These adjoints, however, do not by themselves produce
$\Xi_t\dashv\Gamma_t$ by composition: for example,
$D_{\ell}\Xi U_{\mathrm{MV}}$ and
$I_{\ell}\Xi U_{\mathrm{MV}}$ equip the algebraic Mundici group with
the discrete and indiscrete topologies, respectively, rather than with
the Mundici topology $\sigma_A$.

The relevant general lifting principle is instead the
\emph{Taut Lift Theorem}; see, e.g.,
\cite[Theorem~21.28]{Adamek1990} and
\cite[Theorem~II.5.11.1]{Hofmann2014}.
To apply it here, observe that $\Gamma_t$ sends
$U_{\ell}$-initial sources to $U_{\mathrm{MV}}$-initial sources.
Indeed, suppose that the topology $\rho$ on $(G,u)$ is initial
for the family of unit-preserving $\ell$-homomorphisms
\[
 \bigl(h_i\colon(G,u)\longrightarrow(G_i,u_i,\rho_i)\bigr)_{i\in I}.
\]
The subspace topology $\rho|_{[0,u]}$ is then generated by sets of the
form
\[
 [0,u]\cap h_i^{-1}(O)
 =
 (\Gamma h_i)^{-1}\bigl(O\cap[0,u_i]\bigr),
 \qquad O\in\rho_i.
\]
Hence $\rho|_{[0,u]}$ is exactly the initial topology induced by
the family of restricted MV-homomorphisms
\[
 \bigl(\Gamma h_i\colon
 [0,u]\longrightarrow([0,u_i],\rho_i|_{[0,u_i]})\bigr)_{i\in I}.
\]
Since $\Xi\dashv\Gamma$ algebraically, the Taut Lift Theorem therefore
gives a left adjoint to $\Gamma_t$ whose underlying algebraic functor is
$\Xi$.

In fact, the construction in the proof of the Taut Lift Theorem
recovers precisely the topology $\sigma_A$.  
For $(A,\tau)\in\TopMV$, consider all morphisms $f_i$
and their algebraic adjoints $\widehat f_i$:
\[
 f_i\colon(A,\tau)\longrightarrow\Gamma_t(G_i,u_i,\rho_i),
 \qquad
 \widehat f_i\colon\Xi A\longrightarrow(G_i,u_i).
\]
The $U_{\ell}$-initial lift of the source $(\widehat f_i)_i$
equips $\Xi A$ with the topology
\[
 \theta_A=\bigvee_i\widehat f_i^{\,*}\rho_i.
\]
For each $i$, let $j_i\colon[0,u_i]\longhookrightarrow G_i$
denote the canonical inclusion.  By
Lemma~\ref{lem-pullback-topology}, each
$\widehat f_i^{\,*}\rho_i$ is an $\ell$-group topology, and it is
admissible because $\widehat f_i\,\iota_A=j_i f_i$.
Hence $\theta_A\subseteq\sigma_A$.
Conversely, if $\rho$ is any admissible $\ell$-group topology on
$\Xi A$, then the canonical unit
\[
 \eta_A\colon(A,\tau)\longrightarrow
 \Gamma_t(G_A,u_A,\rho)
\]
is continuous and therefore occurs among the above source maps.  Its
algebraic adjoint is $1_{\Xi A}$, so
$\rho\subseteq\theta_A$.  Since this holds for every admissible
$\rho$, we obtain $\sigma_A\subseteq\theta_A$ and hence
\[
 \theta_A=\sigma_A.
\]
Thus, in the present setting, the topology produced by the Taut Lift
construction coincides exactly with the Mundici topology
$\sigma_A$ constructed above.
\end{rem}

The diagram \eqref{eq-topological-lifting-square} records what has
been lifted: the classical Mundici equivalence yields a topological
adjunction, but it does not yet show that its unit is a homeomorphism.
The main problem is whether the Mundici topology $\sigma_A$
associated with every $(A,\tau)\in\TopMV$ is an $\ell$-group extension
of $\tau$.  Equivalently, we must decide whether
\[
 \sigma_A|_{[0,u_A]}=\tau.
\]
Indeed, if an $\ell$-group extension $\rho$ exists, its
admissibility gives $\rho\subseteq\sigma_A$, and hence
\[
 \tau=\rho|_A\subseteq\sigma_A|_A\subseteq\tau.
\]
Conversely, equality for $\sigma_A$ makes it an $\ell$-group
extension.  Thus asking whether an $\ell$-group extension exists is
equivalent to asking whether the Mundici topology is such an extension;
this is the topological reconstruction problem for the adjunction.

\section{Coordinate and bounded-stage reconstruction}
\label{sec-coordinate-bounded-reconstruction}

We now seek an $\ell$-group topology on the Mundici group that induces
the given topology on its unit interval.  By Section~\ref{sec-adjunction},
it suffices to show that the Mundici topology $\sigma_A$ has this property.
We first use the unit-interval topology to define a topology
$\tau_n$ on each bounded stage and prove that every $\ell$-group
extension must induce $\tau_n$ there.
For an arbitrary topological MV-algebra, we then construct the
coordinate topology $\kappa_A$ and the bounded-stage final topology
$\lambda_A$ on its Mundici group, both inducing $\tau_n$ on each
stage.  We conclude by relating $\lambda_A$ to the Mundici
topology $\sigma_A$.

\subsection{Uniqueness of bounded-stage topologies}
\label{sec-preliminaries}

Let $(G,u)$ be an abelian $\ell$-group with strong unit $u$,
and let $A=[0,u]$ be its unit-interval MV-algebra.
For $n\geq1$ set
\[
 I_n=[-nu,nu]=\{g\in G:-nu\leq g\leq nu\}.
\]
We call these order intervals $I_n$ the \emph{bounded stages} of $G$.
Since $u$ is a strong unit,
\begin{equation}
\label{eq-bounded-stage-exhaustion}
 G=\bigcup_{n\geq1}I_n.
\end{equation}
Following \cite[Definition~2.10]{Ball2002}, for each $k\in\mathbb Z$, define the \emph{bounded coordinate}
\begin{equation}
\label{eq-pi-k}
 \pi_k\colon G\longrightarrow A,\qquad
 \pi_k(x)=((x-ku)\wedge u)\vee0.
\end{equation}
We collect the coordinate identities that will be used below from
\cite[Proposition~2.11]{Ball2002}.
Here 
\[
 a\odot b=\neg(\neg a\oplus\neg b)
\]
for $a,b\in A$, and $\bigoplus$ denotes the iterated MV-sum.

\begin{lem}[Bounded-coordinate identities]
\label{lem-coordinate-identities}
Let $(G,u)$ be an abelian $\ell$-group with strong unit,
$A=[0,u]$, and $\pi_k$ be as in \eqref{eq-pi-k}.
For $f,g\in G$ and $k\in\mathbb Z$, the following identities hold:
\begin{enumerate}[label=(\arabic*),ref=(\arabic*)]
\item\label{coordinate-identities-interval}
For $a\in A$,
\[
 \pi_k(a)=
 \begin{cases}
  u,&k<0,\\
  a,&k=0,\\
  0,&k>0.
 \end{cases}
\]
\item\label{eq-bounded-coordinate-tails}
For $n\geq1$ and $g\in I_n$,
\[
 \pi_k(g)=
 \begin{cases}
  0,&k\geq n,\\
  u,&k\leq-n-1.
 \end{cases}
\]
\item\label{coordinate-identities-lattice}
The coordinates preserve both lattice operations:
\[
 \pi_k(f\vee g)=\pi_k(f)\vee\pi_k(g),
 \qquad
 \pi_k(f\wedge g)=\pi_k(f)\wedge\pi_k(g).
\]
\item\label{coordinate-identities-absorption}
Adjacent coordinates satisfy
\[
 \pi_k(g)\oplus\pi_{k+1}(g)=\pi_k(g).
\]
\item\label{coordinate-identities-inversion}
Negation in $G$ satisfies
\[
 \pi_k(-g)=\neg\pi_{-k-1}(g).
\]
\item\label{coordinate-identities-addition}
Addition in $G$ satisfies
\[
 \pi_k(f+g)
 =\bigoplus_{i+j=k-1}\bigl(\pi_i(f)\odot\pi_j(g)\bigr).
\]
Only finitely many summands are nonzero.
\item\label{eq-coordinate-reconstruction}
For $n\geq1$ and $g\in I_n$,
\[
 g=\sum_{k=-n}^{n-1}\pi_k(g)-nu,
\]
where the sum is taken in $G$.
\end{enumerate}
\end{lem}

\begin{proof}
\ref{coordinate-identities-interval} and~\ref{eq-bounded-coordinate-tails}
follow from \cite[Proposition~2.11(1)--(3)]{Ball2002}, and
\ref{coordinate-identities-lattice}--\ref{coordinate-identities-inversion}
are \cite[Proposition~2.11(8)--(10)]{Ball2002}, respectively.

The addition identity in~\ref{coordinate-identities-addition} is
\cite[Proposition~2.11(11)]{Ball2002}.
To see explicitly that only finitely many summands can be nonzero,
let $f\in I_n$ and $g\in I_m$.
By~\ref{eq-bounded-coordinate-tails}, a nonzero summand requires
$i\leq n-1$ and $j\leq m-1$.
Since $i+j=k-1$, it follows that
\[
 k-m\leq i\leq n-1,\qquad j=k-1-i.
\]
Thus, on fixed bounded stages, the formula is a finite MV-term
in the input coordinates.

Finally, \ref{eq-coordinate-reconstruction} is
\cite[Proposition~2.11(16)]{Ball2002}.
It also follows directly by summing the $\ell$-group identity
\[
 \pi_k(g)=(g\wedge(k+1)u)-(g\wedge ku).
\]
For $g\in I_n$, this gives
\[
 \sum_{k=-n}^{n-1}\pi_k(g)
 =(g\wedge nu)-(g\wedge(-nu))
 =g+nu,
\]
as required.
\end{proof}

For each $n\geq1$ define the maps
\begin{equation}
\label{eq-stage-coordinate-map}
 s_n\colon I_n\longrightarrow A^{2n},
 \qquad
 s_n(g)=(\pi_{-n}(g),\ldots,\pi_{n-1}(g)),
\end{equation}
and
\begin{equation}
\label{eq-stage-sum-map}
 q_n\colon A^{2n}\longrightarrow I_n,
 \qquad
 q_n((a_k)_{k=-n}^{n-1})
 =\sum_{k=-n}^{n-1}a_k-nu.
\end{equation}
The latter map is well defined: each $a_k$ belongs to $[0,u]$,
so
\[
 0\leq\sum_{k=-n}^{n-1}a_k\leq2nu,
 \qquad
 -nu\leq\sum_{k=-n}^{n-1}a_k-nu\leq nu.
\]
By Lemma~\ref{lem-coordinate-identities}\ref{eq-coordinate-reconstruction},
\begin{equation}
\label{eq-stage-section}
 q_ns_n=1_{I_n}.
\end{equation}
Now equip $A$ with a topology $\tau$.
For each $n\geq1$, define the \emph{bounded-stage topology} $\tau_n$
on $I_n$ to be the initial topology induced by the map
\begin{equation}
\label{eq-stage-initial-map}
 s_n\colon I_n\longrightarrow(A^{2n},\tau^{2n}),
\end{equation}
where $\tau^{2n}$ is the product topology on $A^{2n}$.
Equivalently, $\tau_n$ is the coarsest topology making the
coordinates occurring in $s_n$ continuous.
The following proposition shows that every $\ell$-group extension
of $\tau$ must induce these stage topologies.

\begin{prop}[Uniqueness of bounded-stage topologies]
\label{prop-lgroup-extension-bounded-stages}
Let $(G,u)$ be a unital abelian $\ell$-group, and let $\tau$ be
a topology on $A=[0,u]$.
Every $\ell$-group topology $\rho$ on $G$ with $\rho|_A=\tau$
satisfies
\[
 \rho|_{I_n}=\tau_n
 \qquad(n\geq1).
\]
In particular, any two $\ell$-group topologies on $G$ agreeing
on $A$ agree on every bounded stage.
\end{prop}

\begin{proof}
Fix $n\geq1$.
We first verify continuity of the maps
\[
 s_n\colon(I_n,\rho|_{I_n})
 \longrightarrow(A^{2n},\tau^{2n}),
 \qquad
 q_n\colon(A^{2n},\tau^{2n})
 \longrightarrow(I_n,\rho|_{I_n}).
\]
Every coordinate $\pi_k$ is continuous from $(G,\rho)$
to $(A,\rho|_A)=(A,\tau)$, since translations and lattice
operations are continuous for an $\ell$-group topology.
Restricting these coordinates to $I_n$ proves continuity of $s_n$.
Continuity of $q_n$ follows from joint continuity of addition and
translation by $-nu$, followed by corestriction to the subspace $I_n$.

By \eqref{eq-stage-section}, $q_n$ is a continuous left inverse
of $s_n$ on its image, so $s_n$ is a topological embedding.
Thus $\rho|_{I_n}$ is the initial topology induced by $s_n$,
which is $\tau_n$ by definition.
\end{proof}

Thus the unit-interval topology uniquely determines the topology
on every bounded stage of any $\ell$-group extension.  The remaining
task is to construct suitable topologies on the whole group while
retaining these prescribed restrictions.

\subsection{The coordinate topology
\texorpdfstring{$\kappa_A$}{kappa\_A}: a lower bound}
\label{sec-coordinate-topology}

Now let $(A,\tau)\in\TopMV$ and write $\Xi A=(G_A,u_A)$.
For simplicity, write 
\[
 (G,u)=(G_A,u_A)
\]
and identify
$A=[0,u]$ via \eqref{eq-Mundici-unit}.
The bounded stages $I_n$, their topologies $\tau_n$, and the
algebraic maps $\pi_k$, $s_n$, and $q_n$ above are now understood
on this Mundici group.
Whenever a bounded stage $I_n$ is regarded as a topological space,
it carries $\tau_n$ unless otherwise stated.

We first construct a topology on $G$ that recovers these stage
topologies.  For $g\in G$, define the
\emph{translated unit-interval truncation}
\begin{equation}
\label{eq-translated-truncation}
 t_g\colon G\longrightarrow A,
 \qquad t_g(x)=((x+g)\wedge u)\vee0.
\end{equation}
Let $\kappa_A$ be the initial topology induced by
\begin{equation}
\label{eq-translated-initial-family}
 \bigl(t_g\colon G\longrightarrow(A,\tau)\bigr)_{g\in G}.
\end{equation}
Thus $\kappa_A$ is the coarsest topology making all $t_g$
continuous.  Equivalently, for every topological space $X$, a map
$f\colon X\longrightarrow(G,\kappa_A)$ is continuous if and
only if all $t_gf$ are continuous into $(A,\tau)$.
The bounded coordinates \eqref{eq-pi-k} are special cases:
\[
 \pi_k=t_{-ku}\qquad(k\in\mathbb Z).
\]
We shall show that these coordinates alone determine $\kappa_A$,
which explains the name \emph{coordinate topology}.

By Lemma~\ref{lem-coordinate-identities}\ref{coordinate-identities-lattice},
\ref{coordinate-identities-inversion}, and~\ref{coordinate-identities-addition}, every bounded-stage
coordinate of an $\ell$-group operation is a finite MV-term
in the input coordinates.  Such terms are continuous in
$(A,\tau)$, giving the following stage properties.

\begin{lem}
\label{lem-stage-topology}
For all $n,m\geq1$:
\begin{enumerate}[label=(\alph*),ref=(\alph*)]
\item\label{lem-stage-topology-quotient}
$q_n\colon(A^{2n},\tau^{2n})\longrightarrow(I_n,\tau_n)$
is a quotient map, and
$s_n\colon(I_n,\tau_n)\longrightarrow(A^{2n},\tau^{2n})$
is a topological embedding and a section of $q_n$;
\item\label{lem-stage-topology-embedding}
if $n\leq m$, the inclusion
$(I_n,\tau_n)\longhookrightarrow(I_m,\tau_m)$ is a topological embedding;
\item\label{lem-stage-topology-operations}
the restrictions of unary $-$, $\vee$, $\wedge$, and $+$ to bounded
stages are continuous:
\[
 -\colon I_n\longrightarrow I_n,\qquad
 \vee,\wedge\colon I_n^2\longrightarrow I_n,\qquad
 +\colon I_n\times I_m\longrightarrow I_{n+m}.
\]
\end{enumerate}
\end{lem}

\begin{proof}
\ref{lem-stage-topology-quotient}
The map $s_n$ is injective by \eqref{eq-stage-section},
and is a topological embedding by the definition of $\tau_n$
as its initial topology.  Each coordinate of
\[
 s_nq_n\colon(A^{2n},\tau^{2n})
 \longrightarrow(A^{2n},\tau^{2n})
\]
is a finite MV-term by
Lemma~\ref{lem-coordinate-identities}\ref{coordinate-identities-interval},
\ref{coordinate-identities-inversion}, and~\ref{coordinate-identities-addition}; hence $s_nq_n$
is continuous.  As $s_n$ embeds $(I_n,\tau_n)$
into $(A^{2n},\tau^{2n})$, this makes $q_n$ continuous.
Now $q_ns_n=1_{I_n}$ by \eqref{eq-stage-section}, so the continuous map $q_n$
has a continuous section.  Consequently, $q_n$ is a
quotient map and $s_n$ is an embedding.

\ref{lem-stage-topology-embedding}
For $n\leq m$, let $e_{n,m}\colon A^{2n}\longrightarrow A^{2m}$
be the map that prepends $m-n$ copies of $u$ and
appends $m-n$ copies of $0$ to each coordinate tuple.
It is a topological embedding, with inverse on its image
given by projection onto the middle $2n$ coordinates.
Lemma~\ref{lem-coordinate-identities}\ref{eq-bounded-coordinate-tails} shows that,
for the inclusion $i_{n,m}\colon I_n\longhookrightarrow I_m$,
the diagram
\[
\begin{tikzcd}[column sep=large,row sep=large]
(I_n,\tau_n)
  \arrow[r,"i_{n,m}"]
  \arrow[d,"s_n"']
&
(I_m,\tau_m)
  \arrow[d,"s_m"]
\\
(A^{2n},\tau^{2n})
  \arrow[r,"e_{n,m}"']
&
(A^{2m},\tau^{2m})
\end{tikzcd}
\]
commutes.
Since $s_n$ and $e_{n,m}$ are topological embeddings,
so is their composite $e_{n,m}s_n=s_mi_{n,m}$.
As $s_m$ is also a topological embedding,
it follows that $i_{n,m}$ is a topological embedding.

\ref{lem-stage-topology-operations} Consider one of the displayed
bounded operations, say $\omega$, with target $I_r$ (where $r=n$ for unary
negation or either lattice operation, while $r=n+m$ for addition).  Since $\tau_r$ is
the initial topology induced by
$s_r\colon(I_r,\tau_r)\longrightarrow(A^{2r},\tau^{2r})$,
it is enough to prove that $s_r\omega$ is continuous.  Since $\tau^{2r}$
is the product topology, this in turn is equivalent to
continuity of each coordinate of $s_r\omega$.  By
Lemma~\ref{lem-coordinate-identities}\ref{coordinate-identities-lattice},
\ref{coordinate-identities-inversion}, and~\ref{coordinate-identities-addition}, every such coordinate is a finite
MV-term in the input coordinates.  The input coordinate maps are
continuous by the definition of the corresponding stage topologies, and
every finite MV-term operation is continuous because $A$ is a topological
MV-algebra.  Hence every coordinate of $s_r\omega$ is continuous, and therefore
so is $\omega$.
\end{proof}

\begin{rem}
Lemma~\ref{lem-stage-topology} shows that each $I_n$ is a quotient
of $A^{2n}$, that its topology agrees with the subspace topology
inherited from every larger stage, and that the group and lattice
operations are continuous on bounded stages. We still need a topology
on the whole group that induces these stage topologies and makes
the operations jointly continuous.
\end{rem}

For $n\geq1$, let
\begin{equation}
\label{eq-truncation-map}
 r_n\colon G\longrightarrow I_n,
 \qquad r_n(x)=(-nu)\vee(x\wedge nu).
\end{equation}

\begin{lem}[Truncation identity]
\label{lem-truncation-identity}
For all integers $m,C\geq1$ and $w\in I_C$,
\[
 r_m(z+w)=r_m\bigl(r_{m+C}(z)+w\bigr)
 \qquad(z\in G).
\]
\end{lem}

\begin{proof}
Since $w\in I_C$, we have
\[
 w-(m+C)u\leq -mu\leq mu\leq w+(m+C)u.
\]
Translations preserve joins and meets, so distributivity gives
\[
 r_{m+C}(z)+w
 =\bigl((z+w)\vee(w-(m+C)u)\bigr)\wedge(w+(m+C)u).
\]
Truncating this expression to $I_m$, distributivity and the displayed
bounds yield
\[
 r_m\bigl(r_{m+C}(z)+w\bigr)
 =((z+w)\vee(-mu))\wedge mu
 =r_m(z+w).\qedhere
\]
\end{proof}

We now show that the bounded coordinates alone determine $\kappa_A$
and identify its stage restrictions.

\begin{prop}
\label{prop-kappa-lambda}
For $(A,\tau)\in\TopMV$, the topology $\kappa_A$ is also
the initial topology induced by each of the following two families:
\begin{equation}
\label{eq-coordinate-initial-family}
 \bigl(\pi_k\colon G\longrightarrow(A,\tau)\bigr)_{k\in\mathbb Z},
 \qquad
 \bigl(r_n\colon G\longrightarrow(I_n,\tau_n)\bigr)_{n\geq1}.
\end{equation}
Moreover,
\begin{equation}
\label{eq-kappa-stage-restriction}
 \kappa_A|_{I_n}=\tau_n\quad(n\geq1),
\end{equation}
and in particular
\begin{equation}
\label{eq-coordinate-interval-exactness}
 \kappa_A|_A=\tau.
\end{equation}
Finally, every $\ell$-group extension $\rho$ of $\tau$ on $G$
satisfies $\kappa_A\subseteq\rho$.
\end{prop}

\begin{proof}
Let $\theta$ and $\zeta$ be the initial topologies induced by
$(\pi_k)_{k\in\mathbb Z}$ and $(r_n)_{n\geq1}$, respectively,
with the target topologies in \eqref{eq-coordinate-initial-family}.
We prove that all three topologies coincide by establishing
\[
 \zeta\subseteq\theta\subseteq\kappa_A\subseteq\zeta.
\]
We first determine the stage restrictions of $\theta$.

The restriction of an initial topology to a subset is the initial
topology induced by the restricted maps: its subbasic sets are the
inverse images under those restrictions of open sets in the targets.
On $I_n$, the coordinates with $-n\leq k<n$ induce $\tau_n$ by
\eqref{eq-stage-initial-map}, and all remaining coordinates are constant
by Lemma~\ref{lem-coordinate-identities}\ref{eq-bounded-coordinate-tails}.  Hence
\begin{equation}
\label{eq-theta-stage-restriction}
 \theta|_{I_n}=\tau_n.
\end{equation}
In particular, $s_1(a)=(u,a)$ for $a\in A$.
The embedding $s_1$ identifies $A\subseteq I_1$ with
$\{u\}\times A\subseteq A^2$, whose topology is $\tau$; thus
\begin{equation}
\label{eq-stage-one-restriction}
 \tau_1|_A=\tau.
\end{equation}

To prove $\zeta\subseteq\theta$, we show that every
$r_n\colon(G,\theta)\longrightarrow(I_n,\tau_n)$ is continuous.
For $-n\leq k<n$, nested lattice truncation gives
\begin{equation}
\label{eq-truncation-coordinate-factorization}
 \pi_k=(\pi_k|_{I_n})r_n,
\end{equation}
because $[ku,(k+1)u]\subseteq[-nu,nu]$.
Applying Lemma~\ref{lem-coordinate-identities}\ref{eq-coordinate-reconstruction}
to $r_n(x)$ yields
\begin{equation}
\label{eq-truncation-coordinates}
 r_n(x)=q_n\bigl(\pi_{-n}(x),\ldots,\pi_{n-1}(x)\bigr).
\end{equation}
The coordinate map into $(A^{2n},\tau^{2n})$ is $\theta$-continuous,
and $q_n$ is continuous by
Lemma~\ref{lem-stage-topology}\ref{lem-stage-topology-quotient}.
Thus every $r_n$ is $\theta$-continuous, and the defining initial
property of $\zeta$ gives $\zeta\subseteq\theta$.

Next, $\theta\subseteq\kappa_A$ because
$\pi_k=t_{-ku}$ makes every bounded coordinate $\kappa_A$-continuous.
To prove $\kappa_A\subseteq\zeta$, it remains to show that every
$t_g\colon(G,\zeta)\longrightarrow(A,\tau)$ is continuous.
Fix $g\in G$, choose $C\geq1$ with $g\in I_C$, and put $R=C+1$.
By Lemma~\ref{lem-truncation-identity} with $m=1$ and the identity
$\pi_0r_1=\pi_0$, we obtain
\[
 \pi_0(r_R(x)+g)=\pi_0(x+g).
\]
Consequently,
\begin{equation}
\label{eq-translated-truncation-factorization}
 t_g=(t_g|_{I_R})r_R.
\end{equation}
The map $r_R\colon(G,\zeta)\longrightarrow(I_R,\tau_R)$
is continuous by the definition of $\zeta$.
For the other factor, bounded-stage addition sends $x\in I_R$
continuously to $x+g\in I_{R+C}$, and meet with $u$ followed by join
with $0$ is continuous on that stage by
Lemma~\ref{lem-stage-topology}\ref{lem-stage-topology-operations}.
The range lies in $A$, whose subspace topology in $I_{R+C}$ is $\tau$
by Lemma~\ref{lem-stage-topology}\ref{lem-stage-topology-embedding}
and \eqref{eq-stage-one-restriction}.
Corestriction to $(A,\tau)$ is therefore continuous.
Thus every $t_g$ is $\zeta$-continuous, and the defining initial
property of $\kappa_A$ gives $\kappa_A\subseteq\zeta$.
Together with $\zeta\subseteq\theta\subseteq\kappa_A$, this yields
$\kappa_A=\theta=\zeta$.
The stage and interval equalities
\eqref{eq-kappa-stage-restriction} and \eqref{eq-coordinate-interval-exactness}
now follow from \eqref{eq-theta-stage-restriction} and
\eqref{eq-stage-one-restriction}.

Finally, let $\rho$ be an $\ell$-group extension
of $\tau$ on $G$.  Since translations and lattice operations
are continuous for $\rho$ and $\rho|_A=\tau$, each translated
truncation $t_g\colon(G,\rho)\longrightarrow(A,\tau)$
is continuous.  By the defining initial property of $\kappa_A$,
we obtain $\kappa_A\subseteq\rho$.
\end{proof}

Although $\kappa_A$ provides a necessary lower bound for all
$\ell$-group extensions of $\tau$, it need not itself be an
$\ell$-group topology, or even a group topology,
as the following example shows.

\begin{exmp}[The coordinate topology need not be a group topology]
\label{exmp-kappa-not-group}
Let
\[
 S=\mathbb N_{\geq1}\times\{0,1\},
 \qquad P_n=\{(n,0),(n,1)\}\quad(n\geq1),
\]
and let $A$ be the Boolean algebra of finite and cofinite subsets of
$S$.  Regarded as a Boolean MV-algebra, it has operations
\[
 E\oplus F=E\cup F,\qquad
 \neg E=S\setminus E,\qquad 0=\varnothing.
\]
We first equip $A$ with a Hausdorff topology $\tau$ for which these
MV-operations are continuous.

Write $E\triangle F$ for symmetric difference.  For $m\geq1$, put
\[
 \mathcal U_m=
 \left\{\bigcup_{n\in D}P_n:
 D\subseteq\{m,m+1,\ldots\}\text{ is finite}\right\}.
\]
Each $\mathcal U_m$ is a subgroup of the Boolean group
$(A,\triangle)$, the subgroups decrease with $m$, and
$\bigcap_m\mathcal U_m=\{\varnothing\}$.  We take the sets
\[
 E\triangle\mathcal U_m
 :=\{E\triangle H:H\in\mathcal U_m\}
 \qquad(E\in A,\ m\geq1)
\]
as a neighbourhood base at each $E$.  The subgroup properties show
that this defines a translation-invariant group topology $\tau$ on
$(A,\triangle)$, and the trivial intersection makes it Hausdorff.
Here symmetric difference is used only to define the topology, not
as the MV-operation $\oplus$.

To check the MV-operations, note first that complementation is the
translation $E\mapsto E\triangle S$ and is therefore continuous.
It remains to verify joint continuity of intersection, since union
then follows by De Morgan's law.  Fix $E,F\in A$ and $m\geq1$.
Because $E$ and $F$ are finite or cofinite, there is an $N\geq m$
such that each of $E$ and $F$ either contains or is disjoint from
every $P_n$ with $n\geq N$.  If $H,K\in\mathcal U_N$, the sets
$E\cap K$, $H\cap F$, and $H\cap K$ are all unions of finitely many
blocks $P_n$ with $n\geq N$, hence belong to $\mathcal U_N$.
Using the distributive law for intersection over symmetric
difference, we obtain
\[
 \bigl((E\triangle H)\cap(F\triangle K)\bigr)
       \triangle(E\cap F)
 =
 (E\cap K)\triangle(H\cap F)\triangle(H\cap K)
 \in\mathcal U_N\subseteq\mathcal U_m.
\]
Thus intersection is continuous at $(E,F)$.  Together with
continuous complementation, this proves
$(A,\tau)\in\TopMV$.

We now identify its algebraic Mundici group explicitly.
Let
\[
 G=\mathbb Z\mathbf1+\mathbb Z^{(S)}
   =\{g\in\mathbb Z^S:
       g\text{ is constant outside a finite subset of }S\},
 \qquad u=\mathbf1,
\]
with pointwise addition and lattice operations.
The element $u$ is a strong unit, and $[0,u]$ consists precisely
of characteristic functions of finite or cofinite subsets of $S$.
Hence we may identify $A=[0,u]$ and $(G,u)=\Xi A$ algebraically.
For $E\subseteq S$, denote its characteristic function by
$\mathbf1_E$, and set
\[
 b_n=\mathbf1_{P_n},
 \qquad c_n=\mathbf1_{\{(n,0)\}}\qquad(n\geq1).
\]
By construction, $b_n\to0$ in $(A,\tau)$: for every $m$ we have
$P_n\in\mathcal U_m$ whenever $n\geq m$.
In contrast, $c_n\not\to0$, since the singleton $\{(n,0)\}$
does not belong to $\mathcal U_1$.

Consider the following elements of $G$:
\[
 x_n=n b_n,\qquad y_n=-n b_n+c_n.
\]
Their sum is $x_n+y_n=c_n$.  We claim that both $x_n$ and $y_n$
converge to $0$ in $\kappa_A$.
To see this, recall that the coordinate maps are evaluated pointwise
by
\[
 \pi_k(g)(s)=\max\{0,\min\{1,g(s)-k\}\}
 \qquad(s\in S,\ k\in\mathbb Z).
\]
Outside $P_n$, both $x_n$ and $y_n$ vanish.  On $P_n$, the values of
$x_n$ equal $n$, while those of $y_n$ are $1-n$ and $-n$.
Consequently, for each fixed $k\in\mathbb Z$ and all sufficiently
large $n$ (depending on $k$),
\[
 \pi_k(x_n)=
 \begin{cases}
 b_n,&k\geq0,\\
 u,&k<0,
 \end{cases}
 \qquad
 \pi_k(y_n)=
 \begin{cases}
 0,&k\geq0,\\
 u-b_n,&k<0.
 \end{cases}
\]
But
\begin{equation}
\label{eq-zero-coordinate-values}
 \pi_k(0)=
 \begin{cases}
 0,&k\geq0,\\
 u,&k<0.
 \end{cases}
\end{equation}
Since $b_n\to0$ in $\tau$ and complementation is continuous,
$u-b_n\to u$ in $\tau$.  Therefore
$\pi_k(x_n)\to\pi_k(0)$ and
$\pi_k(y_n)\to\pi_k(0)$ for every fixed $k$.
By the initial definition of $\kappa_A$, this gives
\[
 x_n\to0,\qquad y_n\to0
 \quad\text{in }(G,\kappa_A).
\]
On the other hand,
\[
 \pi_0(x_n+y_n)=\pi_0(c_n)=c_n\not\to0
 \quad\text{in }(A,\tau).
\]
Thus $x_n+y_n$ does not converge to $0$ in $\kappa_A$.
Addition is not jointly continuous at $(0,0)$, so $\kappa_A$
is not a group topology, and in particular not an
$\ell$-group topology.  Notice that $\kappa_A$ is nevertheless
Hausdorff: the maps $\pi_k$ separate points by
Lemma~\ref{lem-coordinate-identities}\ref{eq-coordinate-reconstruction},
and $(A,\tau)$ is Hausdorff.
\end{exmp}

\subsection{The bounded-stage final topology
\texorpdfstring{$\lambda_A$}{lambda\_A}: an upper bound}

By Proposition~\ref{prop-kappa-lambda}, the topology $\kappa_A$
induces $\tau_n$ on every bounded stage.  These are precisely
the stage topologies forced on all $\ell$-group extensions by
Proposition~\ref{prop-lgroup-extension-bounded-stages}.
We therefore retain them while allowing the finest possible
global topology.

\begin{defn}
\label{def-Mundici-final-topology}
The \emph{bounded-stage final topology} $\lambda_A$ is the
\emph{bounded-stage finalization} of $\kappa_A$, defined by
\begin{equation}
\label{eq-lambda-stage-test}
 \lambda_A
 =\{U\subseteq G:U\cap I_n\in\kappa_A|_{I_n}
                   \text{ for every }n\geq1\}.
\end{equation}
Equivalently, it is the final topology
of the family of inclusions
\[
 \bigl(j_n\colon(I_n,\tau_n)\longhookrightarrow G\bigr)_{n\geq1}.
\]
\end{defn}

The displayed family of sets is a topology: the stagewise openness
condition is preserved by arbitrary unions and finite intersections.
Every $\kappa_A$-open set satisfies that condition, so
\begin{equation}
\label{eq-kappa-lambda-comparison}
 \kappa_A\subseteq\lambda_A.
\end{equation}
Since, by
Lemma~\ref{lem-stage-topology}\ref{lem-stage-topology-quotient}, each map
$q_n\colon(A^{2n},\tau^{2n})\longrightarrow(I_n,\tau_n)$ is a quotient
map, the preceding
condition is equivalent to
\begin{equation}
\label{eq-lambda-quotient-test}
\begin{aligned}
 U\in\lambda_A
 &\iff
 q_n^{-1}(U\cap I_n)\text{ is open in }(A^{2n},\tau^{2n})
 \quad\text{for every }n \\
 &\iff
 (j_nq_n)^{-1}(U)\text{ is open in }(A^{2n},\tau^{2n})
 \quad\text{for every }n.
\end{aligned}
\end{equation}
Thus openness for $\lambda_A$ can be tested entirely on finite powers of
the original space $A$.

We next verify that $\lambda_A$ retains the prescribed topology on
every bounded stage.  Recall that the truncations $r_n$ from
\eqref{eq-truncation-map} satisfy
\begin{equation}
\label{eq-truncation-retraction}
 r_nj_n=1_{I_n}.
\end{equation}

\begin{prop}
\label{prop-lambda-exact-space}
For every $n\geq1$, the truncation
$r_n\colon(G,\lambda_A)\longrightarrow(I_n,\tau_n)$
is continuous.  Moreover, the inclusion
$(I_n,\tau_n)\longhookrightarrow(G,\lambda_A)$
is a topological embedding; equivalently,
\begin{equation}
\label{eq-lambda-stage-restriction}
 \lambda_A|_{I_n}=\kappa_A|_{I_n}=\tau_n
 \qquad(n\geq1).
\end{equation}
Consequently, $\lambda_A$ is the finest topology on $G$ inducing
the same topology as $\kappa_A$ on every bounded stage: for any topology
$\rho$ on $G$,
\[
 \rho|_{I_n}=\kappa_A|_{I_n}\quad(n\geq1)
 \implies \rho\subseteq\lambda_A.
\]
In particular,
\begin{equation}
\label{eq-stage-one-exactness}
 \lambda_A|_{[0,u]}
 =
 \tau_1|_{[0,u]}
 =
 \tau.
\end{equation}
\end{prop}

\begin{proof}
By Proposition~\ref{prop-kappa-lambda}, the map
$r_n\colon(G,\kappa_A)\longrightarrow(I_n,\tau_n)$ is continuous.
Since $\kappa_A\subseteq\lambda_A$, the map
$r_n\colon(G,\lambda_A)\longrightarrow(I_n,\tau_n)$
is also continuous.
On the other hand, the inclusion
$j_n\colon(I_n,\tau_n)\longrightarrow(G,\lambda_A)$ is continuous
by the final-topology definition.  As
$r_nj_n=1_{I_n}$ by \eqref{eq-truncation-retraction},
$j_n$ is a topological embedding.  Hence
$\lambda_A|_{I_n}=\tau_n=\kappa_A|_{I_n}$,
where the second equality follows from \eqref{eq-kappa-stage-restriction}.

Now suppose a topology $\rho$ on $G$ induces
$\kappa_A|_{I_n}$ on every stage.  For each
$\rho$-open set $U$, all the sets
$U\cap I_n$ are $\tau_n$-open, so the stagewise
test \eqref{eq-lambda-stage-test} gives $U\in\lambda_A$.
Thus $\rho\subseteq\lambda_A$, and the stage equality
just proved makes $\lambda_A$ the finest such topology.

Finally, restricting the stage equality to
$A\subseteq I_1$ and using
\eqref{eq-stage-one-restriction} yields
$\lambda_A|_A=\tau$.
\end{proof}

Proposition~\ref{prop-lambda-exact-space} shows that $\lambda_A$ gives an
exact reconstruction at the level of spaces: all bounded stages retain
their prescribed topologies, and in particular the original topology on
$A$ is recovered without change.  What is not yet known is its
compatibility with the $\ell$-group operations.  Although
Lemma~\ref{lem-stage-topology}\ref{lem-stage-topology-operations} gives
continuous maps
\[
 I_n\times I_m\longrightarrow I_{n+m}
\]
for addition, and analogous bounded-stage maps for the lattice operations,
this does not by itself imply joint continuity on $G\times G$, since
products need not preserve final topologies.

\begin{rem}
\label{rem-ball-bounded-convergence}
The bounded-stage final topology $\lambda_A$ is closely related
to the bounded $\ell u$-convergence in \cite{Ball2002}.
For background on convergence on lattice-ordered
groups, see Ball \cite{Ball1980}. 
Ball, Georgescu, and Leu{\c s}tean \cite[Theorem~3.3]{Ball2002}
establish a one-to-one correspondence between MV-convergences on $A$ and
bounded $\ell u$-convergences on its Mundici group $G$.
For the convergence induced by $\tau$, write $\xi_A^b$ for the
corresponding bounded convergence.  A filter $\mathcal F$ on $G$
satisfies $\mathcal F\xrightarrow{\xi_A^b}g$ if and only if
for some $n\geq1$,
\[
 g\in I_n\in\mathcal F,
 \qquad
 \pi_k(\mathcal F)\longrightarrow\pi_k(g)
 \text{ in }(A,\tau)\quad(-n\leq k<n),
\]
where $\pi_k(\mathcal F)$ denotes the image filter.

The \emph{topological modification} of a convergence
structure $\xi$ is the topology $\operatorname{Top}(\xi)$ in which
$U$ is open if and only if $U\in\mathcal F$ whenever
$\mathcal F\xrightarrow{\xi}g$ and $g\in U$.
For the bounded convergence above,
\[
 \operatorname{Top}(\xi_A^b)=\lambda_A.
\]
Indeed, by the initial definition of $\tau_n$, the bounded
convergence rule says that, for some $n$ with $g\in I_n\in\mathcal F$,
the trace filter $\mathcal F|_{I_n}$ converges to $g$ in
$(I_n,\tau_n)$.
Thus $\xi_A^b$ is the final convergence structure induced by the inclusions
$j_n\colon(I_n,\tau_n)\longhookrightarrow G$, with each stage viewed as a
convergence space.  The same bounded stages therefore yield $\xi_A^b$
in the category of convergence spaces and $\lambda_A$ in the category
of topological spaces.

If $U\in\lambda_A$, stagewise openness therefore
gives $U\in\mathcal F$ whenever $\mathcal F\xrightarrow{\xi_A^b}g\in U$.
Conversely, if $U$ is open in $\operatorname{Top}(\xi_A^b)$, every
filter $\mathcal H$ converging to $g\in U\cap I_n$ in $(I_n,\tau_n)$
has image filter $j_n(\mathcal H)\xrightarrow{\xi_A^b}g$.
Hence $U\cap I_n\in\mathcal H$, so the filter characterization of
open sets gives $U\cap I_n\in\tau_n$.
The stagewise test \eqref{eq-lambda-stage-test} yields $U\in\lambda_A$.

The bounded convergence itself need not be topological, even for
compact Hausdorff metrizable inputs.  To see this, take
$A=[0,1]^{\mathbb N}$ with coordinatewise MV-operations and the
product topology.  Its Mundici group is
\[
 G=\ell^\infty(\mathbb N),\qquad u=\mathbf1.
\]
Let $b_j$ be the $j$th standard basis vector; then $b_j\to0$ in $A$.
For each integer $r\geq1$, let $\mathcal F_r$ be the tail filter of
$(r b_j)_{j\geq1}$, with multiples taken in $G$.
It contains $I_r$ and converges to $0$ under $\xi_A^b$, since
\[
 \pi_k(r b_j)=
 \begin{cases}
  u,&-r\leq k<0,\\
  b_j,&0\leq k<r.
 \end{cases}
\]
If $\xi_A^b$ were induced by a topology, its zero-neighbourhood
filter $\mathcal N(0)$ would itself converge to $0$.
Boundedness would give $I_m\in\mathcal N(0)$ for some $m$, and hence
$I_m\in\mathcal F_r$ for every $r$, since every filter converging
to $0$ contains $\mathcal N(0)$.
For $r>m$, however, no term $r b_j$ lies in $I_m$, a contradiction.

Thus the assertion in \cite[Remark~3.4(2)]{Ball2002} that
topological MV-convergences have topological bounded extensions
fails in general.  This concerns the additional topologicality
assertion, rather than the convergence correspondence of
\cite[Theorem~3.3]{Ball2002}.
The topological modification yields $\lambda_A$;
obtaining joint continuity of the $\ell$-group operations requires
the further reconstruction of Section~\ref{sec-group-envelope}.
\end{rem}

Although $\kappa_A$ and $\lambda_A$ induce the same topology on every
bounded stage, the inclusion
$\kappa_A\subseteq\lambda_A$ can be strict, as the
following example shows.

\begin{exmp}
\label{exmp-kappa-lambda-strict}
Let $A=[0,1]^{\mathbb N}$ have coordinatewise MV-operations and the
product topology $\tau$.  Its Mundici group is
\[
 G=\ell^\infty(\mathbb N),\qquad u=\mathbf1.
\]
Writing $\rho_{\mathrm{pt}}$ for the pointwise topology on $G$, we
show directly that
\[
 \kappa_A=\rho_{\mathrm{pt}}\subsetneq\lambda_A.
\]

First, let $e_j\colon G\longrightarrow\mathbb R$ be the $j$th
coordinate evaluation, and use the same notation for its restriction
to $A$.  For $k\in\mathbb Z$, set
\[
 \varphi_k\colon\mathbb R\longrightarrow[0,1],\qquad
 \varphi_k(t)=((t-k)\wedge1)\vee0.
\]
Then $e_j\pi_k=\varphi_ke_j$.  Each scalar component of $\pi_k$ is
therefore pointwise-continuous, so the product topology on $A$ makes
every $\pi_k\colon(G,\rho_{\mathrm{pt}})\longrightarrow(A,\tau)$
continuous.  Initiality gives $\kappa_A\subseteq\rho_{\mathrm{pt}}$.

For the reverse inclusion, the usual topology on $\mathbb R$ is
initial for $(\varphi_k)_{k\in\mathbb Z}$.  Indeed, all these maps
are usual-continuous.  If $x\in(k,k+1)$ and
$0<\varepsilon<\min\{x-k,k+1-x\}$, then
\[
 \varphi_k^{-1}\bigl((x-k-\varepsilon,x-k+\varepsilon)\bigr)
 =(x-\varepsilon,x+\varepsilon).
\]
At an integer $n$, for $0<\varepsilon<1$,
\[
 \varphi_{n-1}^{-1}\bigl((1-\varepsilon,1]\bigr)
 \cap\varphi_n^{-1}\bigl([0,\varepsilon)\bigr)
 =(n-\varepsilon,n+\varepsilon).
\]
These sets supply a usual neighbourhood base at every point.
Since every $\varphi_ke_j=e_j\pi_k$ is $\kappa_A$-continuous,
this scalar initial property makes each
$e_j\colon(G,\kappa_A)\longrightarrow\mathbb R$ continuous.
The pointwise topology is initial for the evaluations, so
$\rho_{\mathrm{pt}}\subseteq\kappa_A$ and hence
$\kappa_A=\rho_{\mathrm{pt}}$.

Next, consider the set
\[
 U=\{x=(x_j)\in\ell^\infty(\mathbb N):
                   |x_j|<j\text{ for every }j\geq1\}.
\]
On $I_n=[-n\mathbf1,n\mathbf1]$, the conditions with $j>n$ are
automatic, because $|x_j|\leq n<j$.  Thus
\[
 U\cap I_n=\bigcap_{j=1}^n\{x\in I_n:|x_j|<j\}.
\]
By Proposition~\ref{prop-kappa-lambda},
$\tau_n=\kappa_A|_{I_n}=\rho_{\mathrm{pt}}|_{I_n}$.
Thus $U\cap I_n$ is $\tau_n$-open for every $n$, and finality
\eqref{eq-lambda-stage-test} gives $U\in\lambda_A$.

Finally, $0\in U$, but every pointwise basic neighbourhood of $0$
restricts only finitely many scalar coordinates.  Choose an
unrestricted coordinate $j$ and let $\mathbf e_j$ be the $j$th
standard basis vector.  The element $j\mathbf e_j$ belongs
to that neighbourhood and lies outside $U$, since its $j$th
coordinate equals $j$.  Thus $U$ is not pointwise-open at $0$.
Since $\kappa_A=\rho_{\mathrm{pt}}$, this proves
$\kappa_A\subsetneq\lambda_A$.
\end{exmp}

The next example shows that bounded-stage finalization need not make
addition jointly continuous, even for a metrizable input.

\begin{exmp}[The bounded-stage final topology need not be a group topology]
\label{exmp-lambda-not-group}
Put $D=\mathbb Q\cap[0,1]$ with its usual subspace topology, and let
$A=D^{\mathbb N}$ have the product topology and coordinatewise
MV-operations.  Its algebraic Mundici group is
\[
 G=\{x\in\mathbb Q^{\mathbb N}:\sup_j|x_j|<\infty\},
 \qquad u=\mathbf1.
\]
Indeed, this is a unital abelian $\ell$-group with unit interval $A$;
the bounded-coordinate reconstruction expresses every element as a
finite sum of interval elements followed by an integer translation.
The scalar-truncation argument of
Example~\ref{exmp-kappa-lambda-strict}, restricted to rational
coordinates, shows that $\kappa_A$ is the pointwise topology.
Consequently, $\tau_N$ is the product topology on
\[
 I_N=(\mathbb Q\cap[-N,N])^{\mathbb N}.
\]

Let $\mathbf e_j$ denote the $j$th standard basis vector.
Choose irrational $\alpha_n\in(0,1)$ with $\alpha_n\to0$ and,
for each $n\geq1$, rationals $q_{n,m}\in(0,1)$ with
$q_{n,m}\to\alpha_n$ as $m\to\infty$ $(m\geq2)$.  Set
\[
 F_n=\{q_{n,m}\mathbf e_1+n\mathbf e_m:m\geq2\},
 \qquad F=\bigcup_{n\geq1}F_n,
 \qquad W=G\setminus F.
\]
For fixed $n$, the displayed sequence converges in
$\mathbb R^{\mathbb N}$ to $\alpha_n\mathbf e_1\notin G$.
The sequence together with its limit is compact, hence closed
in $\mathbb R^{\mathbb N}$; its intersection with $G$ is $F_n$.
Thus $F_n$ is pointwise-closed in $G$.
Every point of $F_n$ has sup norm exactly $n$, so
\[
 F\cap I_N=\bigcup_{n=1}^N F_n
\]
is closed in $(I_N,\tau_N)$ for every $N$.
The final-topology criterion therefore makes $F$ closed in
$\lambda_A$.  Since $0\notin F$, the set $W$ is a
$\lambda_A$-open zero-neighbourhood.

We show that no two $\lambda_A$-open zero-neighbourhoods $U,V$
satisfy $U+V\subseteq W$.
Openness of $U\cap I_1$ gives $\varepsilon>0$ such that
$q\mathbf e_1\in U$ for every rational $0\leq q<\varepsilon$.
Choose $n$ with $\alpha_n<\varepsilon$.  Then
$q_{n,m}\mathbf e_1\in U$ for all sufficiently large $m$.
For this fixed $n$, the sequence $(n\mathbf e_m)_{m\geq2}$ converges to $0$
in $(I_n,\tau_n)$, so $n\mathbf e_m\in V$ for all sufficiently
large $m$ as well.  Their sum belongs to $F_n$, hence lies
outside $W$.  Addition is therefore not continuous at $(0,0)$
for $\lambda_A$.
\end{exmp}

By Proposition~\ref{prop-lambda-exact-space}, the bounded reconstruction
$\lambda_A$ already recovers $\tau$ exactly on the unit interval at the
level of spaces.  Example~\ref{exmp-lambda-not-group} shows that
$\lambda_A$ need not even be a group topology, so a compatible
$\ell$-group extension of $\tau$ still has to be obtained.  The topological Mundici adjunction established in
Theorem~\ref{thm-topological-adjunction} uses the Mundici topology
$\sigma_A$, the largest admissible $\ell$-group topology.
We now show that $\sigma_A$ is the finest $\ell$-group topology
contained in $\lambda_A$.

\begin{thm}
\label{thm-lambda-envelope}
For $(A,\tau)\in\TopMV$, write
\[
 \Xi A=(G,u)
 \quad\text{and}\quad
 \Xi_t(A,\tau)=(G,u,\sigma_A),
\]
and identify $A=[0,u]$ via the classical Mundici unit.
Then:
\begin{enumerate}[label=(\alph*),ref=(\alph*)]
\item\label{thm-lambda-envelope-characterization}
an $\ell$-group topology $\rho$ on $G$ is admissible if and only if $\rho\subseteq\lambda_A$;
\item\label{thm-lambda-envelope-maximality}
$\sigma_A$ is the finest $\ell$-group topology on $G$ contained in
$\lambda_A$.
\end{enumerate}
\end{thm}

\begin{proof}
\ref{thm-lambda-envelope-characterization}
Suppose first that $\rho$ is admissible.  By
Definition~\ref{def-admissible-topology}, the inclusion
$\iota_A\colon(A,\tau)\longrightarrow(G,\rho)$ is continuous.
Since $\rho$ is a group topology, for every $n\geq1$ the map
\[
 j_nq_n\colon(A^{2n},\tau^{2n})\longrightarrow(G,\rho),
 \qquad
 (a_k)_{k=-n}^{n-1}\mapsto\sum_{k=-n}^{n-1}a_k-nu,
\]
is continuous by finite addition and translation.
Here $j_n$ is the inclusion from
Definition~\ref{def-Mundici-final-topology}.
For every $\rho$-open set $U$, each inverse image
$(j_nq_n)^{-1}(U)$ is therefore open in $(A^{2n},\tau^{2n})$.
The quotient test \eqref{eq-lambda-quotient-test} gives
$U\in\lambda_A$, so $\rho\subseteq\lambda_A$.

Conversely, let $\rho$ be an $\ell$-group topology on $G$ satisfying
$\rho\subseteq\lambda_A$.  Restricting the two topologies to
$A=[0,u]$ and using Proposition~\ref{prop-lambda-exact-space}, we obtain 
$\rho|_A\subseteq\lambda_A|_A=\tau$. 
Hence, for every $\rho$-open subset $U$ of $G$,
\[
 \iota_A^{-1}(U)=U\cap A
\]
is $\tau$-open.  Therefore $\iota_A\colon(A,\tau)\longrightarrow(G,\rho)$ 
is continuous, so $\rho$ is admissible.  This proves
\ref{thm-lambda-envelope-characterization}.

\ref{thm-lambda-envelope-maximality}
By \ref{thm-lambda-envelope-characterization},
\[
 \mathcal T_A
 =
 \{\rho:
   \rho\text{ is an $\ell$-group topology on }G
   \text{ and }\rho\subseteq\lambda_A\}.
\]
Lemma~\ref{lem-finest-admissible} says that $\sigma_A$ is the largest
element of $\mathcal T_A$.  Hence $\sigma_A$ is precisely the finest
$\ell$-group topology on $G$ contained in $\lambda_A$.
\end{proof}

The final topology $\lambda_A$ has more structure than
stage preservation alone suggests.  Its translations are
homeomorphisms and inversion is continuous.

Recall that a \emph{quasitopological group} \cite{Arhangelskii2008} 
is a group equipped with a
topology for which multiplication is separately continuous and inversion
is continuous.

\begin{prop}
\label{prop-lambda-quasitopological}
For every $a\in G$, the translation
\[
 T_a\colon(G,\lambda_A)\longrightarrow(G,\lambda_A),
 \qquad T_a(x)=a+x,
\]
is a homeomorphism, and inversion is continuous.  In particular,
$(G,\lambda_A)$ is a quasitopological abelian group.
\end{prop}

\begin{proof}
Fix $a\in G$, and choose $r\geq1$ with $a\in I_r$.
For every $n\geq1$, translation restricts to
$T_a|_{I_n}\colon(I_n,\tau_n)\longrightarrow(I_{n+r},\tau_{n+r})$.
This map is continuous: it is the composite of the continuous map
$c_a\colon I_n\longrightarrow I_r\times I_n$, $c_a(x)=(a,x)$,
and bounded-stage addition from
Lemma~\ref{lem-stage-topology}\ref{lem-stage-topology-operations}.
With the stage inclusions $j_n$, the diagram
\[
\begin{tikzcd}[column sep=large,row sep=large]
&
(I_n,\tau_n)
  \arrow[dl,"c_a"']
  \arrow[d,"T_a|_{I_n}"]
  \arrow[r,"j_n"]
&
(G,\lambda_A)
  \arrow[d,"T_a"]
\\
(I_r,\tau_r)\times(I_n,\tau_n)
  \arrow[r,"+"']
&
(I_{n+r},\tau_{n+r})
  \arrow[r,"j_{n+r}"']
&
(G,\lambda_A)
\end{tikzcd}
\]
commutes.  Thus $T_aj_n=j_{n+r}\circ+\circ c_a$ is continuous
for every $n$.  Finality of $\lambda_A$ gives continuity of $T_a$.
The same argument applies to its inverse $T_{-a}$, so $T_a$ is a
homeomorphism.

Inversion restricts to a continuous map
$-\colon(I_n,\tau_n)\longrightarrow(I_n,\tau_n)$ by
Lemma~\ref{lem-stage-topology}\ref{lem-stage-topology-operations}.
The diagram
\[
\begin{tikzcd}
(I_n,\tau_n)
  \arrow[r,"-"]
  \arrow[d,"j_n"']
&
(I_n,\tau_n)
  \arrow[d,"j_n"]
\\
(G,\lambda_A)
  \arrow[r,"-"']
&
(G,\lambda_A)
\end{tikzcd}
\]
commutes for every $n$, so $-\circ j_n=j_n\circ-$ is continuous.
Finality again gives continuity of inversion on $(G,\lambda_A)$.
\end{proof}

We have therefore constructed a quasitopological abelian
group $(G,\lambda_A)$ that induces $\tau_n$ on every
bounded stage.  At this point, the established
comparisons are
\[
 \kappa_A\subseteq\lambda_A,\qquad
 \sigma_A\subseteq\lambda_A,
\]
while
\[
 \kappa_A|_A=\lambda_A|_A=\tau,\qquad
 \sigma_A|_A\subseteq\tau.
\]
Every $\ell$-group extension $\rho$ of $\tau$ must satisfy
\[
 \kappa_A\subseteq\rho\subseteq\lambda_A,
\]
by Proposition~\ref{prop-kappa-lambda} and
Theorem~\ref{thm-lambda-envelope}\ref{thm-lambda-envelope-characterization},
respectively.
Section~\ref{sec-group-envelope} constructs such an extension
from the topological-group reflection of $\lambda_A$.

\section{From group reflection to topological
\texorpdfstring{$\ell$}{l}-group reconstruction}
\label{sec-group-envelope}

We continue with $(A,\tau)\in\TopMV$, $\Xi A=(G,u)$, and the
algebraic identification $A=[0,u]$ from
Section~\ref{sec-coordinate-bounded-reconstruction}.
Our goal is to show that the topological-group reflection of
$\lambda_A$ is an $\ell$-group extension of $\tau$ and coincides
with the Mundici topology $\sigma_A$.

Throughout, a \emph{group topology} requires joint continuity of
addition and continuity of inversion; joint continuity of the
lattice operations is a further property to be proved.

\subsection{The group reflection
\texorpdfstring{$\gamma_A$}{gamma\_A} and its coordinate-box neighbourhoods}
\label{subsec-bounded-control-machinery}

We first construct the group reflection $\gamma_A$ of $\lambda_A$,
then give a zero-neighbourhood base consisting of unions of finite
sums of symmetric coordinate boxes.  Let
\begin{equation}
\label{eq-gamma-A-def}
\gamma_A=\bigvee\{\theta\subseteq\lambda_A:
   \theta\text{ is a group topology on }G\}
\end{equation}
be the group reflection of
$(G,\lambda_A)$, that is, the finest group topology contained
in $\lambda_A$.  We also write
$\gamma_A=R_{\mathrm{grp}}(\lambda_A)$.  This join is a group topology by the same
subbasis argument as in Lemma~\ref{lem-finest-admissible};
the defining family contains the indiscrete topology.  Equivalently, let $H$ be
a topological group and let $f\colon(G,\lambda_A)\longrightarrow H$ be a 
continuous group homomorphism.  Then there is a unique continuous
group homomorphism $\overline f\colon(G,\gamma_A)\longrightarrow H$
for which the triangle
\[
\begin{tikzcd}[column sep=large,row sep=large]
(G,\lambda_A)
  \arrow[r,"1_G"]
  \arrow[dr,"f"']
&
(G,\gamma_A)
  \arrow[d,dashed,"\overline f"]
\\
&
H
\end{tikzcd}
\]
commutes.

To work concretely with $\gamma_A$, we need an explicit
neighbourhood base at $0$ for the group reflection.  The result below
is also a special case of \cite[Lemma~11]{Tkachenko2013}; see
Remark~\ref{rem-group-reflection-Tkachenko}.  For the reader's
convenience, however, we give a direct proof in the present abelian
setting.
For subsets $E,F\subseteq G$, we use the notation
\[
E+F=\{e+f:e\in E,\ f\in F\},\qquad
-F=\{-f:f\in F\},\qquad
E-F=E+(-F).
\]
A subset $E\subseteq G$ is called \emph{symmetric} if $E=-E$.

\begin{lem}
\label{lem-group-reflection-base}
Let $(U_i)_{i\geq1}$ be a sequence of symmetric open
zero-neighbourhoods in $(G,\lambda_A)$, and set
\begin{equation}
\label{eq-group-reflection-base}
 O(U_1,U_2,\ldots)
 =
 \bigcup_{n\geq1}(U_1+\cdots+U_{2n}).
\end{equation}
As $(U_i)_{i\geq1}$ ranges over all such sequences, the sets
$O(U_1,U_2,\ldots)$ form a zero-neighbourhood base for $\gamma_A$.
\end{lem}

\begin{proof}
First, we show that the sets in \eqref{eq-group-reflection-base}
form a zero-neighbourhood base for a group topology $\rho$.

Every such set contains $0$ and is symmetric, and termwise
intersections give common refinements:
\[
 O(U_1\cap V_1,U_2\cap V_2,\ldots)
 \subseteq O(U_1,U_2,\ldots)\cap O(V_1,V_2,\ldots).
\]
For the addition condition, put $C_i=U_{2i-1}\cap U_{2i}$.
Then
\[
 O(C_1,C_2,\ldots)+O(C_1,C_2,\ldots)
 \subseteq O(U_1,U_2,\ldots).
\]
Indeed, pad two representing sums with zeros to the same even
length $2N$.  Their $i$th summands both belong to
$C_i\subseteq U_{2i-1}\cap U_{2i}$, so interleaving them gives
a sum in $U_1+\cdots+U_{4N}$.

For the local-translation condition, if
$x\in U_1+\cdots+U_{2n}$, set $D_i=U_{2n+i}$.  Appending any
finite even sum from this tail gives
\[
 x+O(D_1,D_2,\ldots)\subseteq O(U_1,U_2,\ldots).
\]
These properties give the usual neighbourhood-base criterion for
a group topology; conjugation imposes no further condition because
$G$ is abelian.

Second, we show that $\rho$ is the finest group topology below $\lambda_A$.
Finite sums of $\lambda_A$-open sets are open, since translations
are homeomorphisms by
Proposition~\ref{prop-lambda-quasitopological}.
Thus every $O(U_1,U_2,\ldots)$ and each of its translates is
$\lambda_A$-open, giving $\rho\subseteq\lambda_A$.

Let $\theta\subseteq\lambda_A$ be any group topology, and let $W$
be a $\theta$-open zero-neighbourhood.  Continuity of addition and
inversion allows us to choose symmetric $\theta$-open
zero-neighbourhoods $(V_n)_{n\geq0}$ recursively with
\[
 V_0\subseteq W,\qquad V_{n+1}+V_{n+1}\subseteq V_n.
\]
For every $m\geq1$, backward induction gives
$V_k+\cdots+V_m\subseteq V_{k-1}$ for $1\leq k\leq m$:
the last term lies in $V_{m-1}$, and the induction step uses
$V_k+V_k\subseteq V_{k-1}$.
In particular,
\[
 V_1+\cdots+V_m\subseteq V_0\subseteq W
 \qquad(m\geq1).
\]
Each $V_i$ is $\lambda_A$-open, so
$O(V_1,V_2,\ldots)$ is a $\rho$-zero-neighbourhood contained in $W$.
Since group topologies are determined by their zero-neighbourhoods
and translations, $\theta\subseteq\rho$.
Hence $\rho=\gamma_A$.
\end{proof}

\begin{rem}
\label{rem-group-reflection-Tkachenko}
Lemma~\ref{lem-group-reflection-base} is the abelian special case of
the neighbourhood description of the group reflection in
\cite[Lemma~11]{Tkachenko2013}.  In the notation above, that result gives, 
for an arbitrary sequence $(U_i)_{i\geq1}$ of $\lambda_A$-neighbourhoods of $0$, 
basic neighbourhoods of the form
\[
 \bigcup_{n\geq1}\ \bigcup_{\pi\in S_{2n}}
 \bigl(
 U_{\pi(1)}-U_{\pi(2)}+\cdots+
 U_{\pi(2n-1)}-U_{\pi(2n)}
 \bigr),
\]
where $S_{2n}$ denotes the symmetric group on
$\{1,\ldots,2n\}$.

In our setting this formula simplifies to
\eqref{eq-group-reflection-base}.
Inversion is a homeomorphism for $\lambda_A$, so every
neighbourhood of $0$ contains a symmetric open one: take
$V\cap(-V)$ inside it.  We may therefore choose every $U_i$
symmetric.  Each difference in the formula then becomes a sum,
and commutativity gives, for every $\pi\in S_{2n}$,
\[
 U_{\pi(1)}-U_{\pi(2)}+\cdots+
 U_{\pi(2n-1)}-U_{\pi(2n)}
 =U_1+\cdots+U_{2n}.
\]
Thus the permutation union has just one value, and the union over
$n$ is precisely $O(U_1,U_2,\ldots)$.
\end{rem}

Lemma~\ref{lem-coordinate-identities}\ref{eq-bounded-coordinate-tails} 
shows that the bounded
coordinates are eventually $0$ in the positive direction and $u$ in
the negative direction.
Now we replace each negative-index coordinate $\pi_k(g)$ by
$u-\pi_k(g)$ and leave the nonnegative-index coordinates unchanged.

For $k\in\mathbb Z$ define the \emph{$k$th Mundici deviation} by
\begin{equation}
\label{eq-Mundici-deviation}
 \delta_k(g)=
 \begin{cases}
  \pi_k(g),&k\geq0,\\
  \neg\pi_k(g),&k<0.
 \end{cases}
\end{equation}
For each $g\in G$, only finitely many deviations $\delta_k(g)$ are
nonzero, and their signed sum recovers $g$, as the next lemma shows.

Before turning to that lemma, we record that this normalization
leaves the topologies induced by the coordinates unchanged.  Write
\[
 \pi,\delta\colon G\longrightarrow A^{\mathbb Z},
 \qquad
 \pi(g)=(\pi_k(g))_k,\quad \delta(g)=(\delta_k(g))_k,
\]
and define
\begin{equation}
\label{eq-coordinate-negation-homeomorphism}
 \Theta\colon A^{\mathbb Z}\longrightarrow A^{\mathbb Z},
 \qquad
 (\Theta(a))_k=
 \begin{cases}
  a_k,&k\geq0,\\
  \neg a_k,&k<0.
 \end{cases}
\end{equation}
Since MV-negation is an involutive homeomorphism, $\Theta$ is a
self-inverse homeomorphism for both the product and box topologies,
and $\delta=\Theta\pi$.  In particular,
Proposition~\ref{prop-kappa-lambda} shows that both $\pi$ and $\delta$
induce the initial topology $\kappa_A$ on $G$ when $A^{\mathbb Z}$
is endowed with the product topology.

For each family $\mathcal W=(W_k)_{k\in\mathbb Z}$ of open
zero-neighbourhoods in $A$, define the \emph{coordinate box}
\begin{equation}
\label{eq-coordinate-box-Q}
 Q(\mathcal W)=\bigcap_{k \in \mathbb Z } \delta_k^{-1}(W_k)
\end{equation}
and the \emph{symmetric coordinate box}
\begin{equation}
\label{eq-coordinate-box-B}
 B(\mathcal W)=Q(\mathcal W)\cap(-Q(\mathcal W))
\end{equation}
determined by $\mathcal W$.

\begin{lem}
\label{lem-deviation-decomposition}
For every $g\in G$, only finitely many $\delta_k(g)$ are nonzero, and
\begin{equation}
\label{eq-deviation-decomposition}
 g=\sum_{k\geq0}\delta_k(g)-\sum_{k<0}\delta_k(g).
\end{equation}
Moreover, for every family $\mathcal W$ as above,
$Q(\mathcal W)$ and $B(\mathcal W)$ are $\lambda_A$-open
zero-neighbourhoods, and $B(\mathcal W)$ is symmetric.
\end{lem}

\begin{proof}
Let $g\in G$.  By \eqref{eq-bounded-stage-exhaustion}, choose $m\geq1$
such that $g\in I_m$.  Then
Lemma~\ref{lem-coordinate-identities}\ref{eq-bounded-coordinate-tails} and
\eqref{eq-Mundici-deviation} give
\begin{equation}
\label{eq-deviation-stage-support}
 \delta_k(g)=0
 \qquad
 (k\geq m\text{ or }k\leq-m-1).
\end{equation}
Thus only finitely many deviations of $g$ are nonzero.  Moreover,
Lemma~\ref{lem-coordinate-identities}\ref{eq-coordinate-reconstruction} gives
\[
g
=\sum_{k=-m}^{-1}(u-\delta_k(g))
  +\sum_{k=0}^{m-1}\delta_k(g)-mu
 =\sum_{k=0}^{m-1}\delta_k(g)
  -\sum_{k=-m}^{-1}\delta_k(g),
\]
which proves \eqref{eq-deviation-decomposition}.

Now let $\mathcal W=(W_k)_{k\in\mathbb Z}$ be a family of open
zero-neighbourhoods in $A$.
Since $\delta_k(0)=0$ for every $k\in\mathbb Z$, we have
$0\in Q(\mathcal W)$.
By Definition~\ref{def-Mundici-final-topology}, $Q(\mathcal W)$ is
$\lambda_A$-open if and only if $Q(\mathcal W)\cap I_m$ is open
in $I_m$ for every $m\geq1$.
Fix $m\geq1$.  The identity \eqref{eq-deviation-stage-support}
gives the commutative square
\[
\begin{tikzcd}[column sep=large,row sep=large]
 (I_m,\tau_m) \arrow[r,hook] \arrow[d,"\Delta_m"'] &
 (G,\lambda_A) \arrow[d,"\delta"] \\
 (A^{2m},\tau^{2m}) \arrow[r,hook,"\Omega_m"'] &
 (A^{\mathbb Z},\tau_{\mathrm{box}}),
\end{tikzcd}
\]
where $\Delta_m(g)=(\delta_{-m}(g),\ldots,\delta_{m-1}(g))$ and
$\Omega_m$ extends each tuple by zero outside the indices
$-m,\ldots,m-1$, while $\tau_{\mathrm{box}}$ denotes the box topology
on $A^{\mathbb Z}$.
By \eqref{eq-Mundici-deviation}, $\Delta_m$ is obtained from the
stage map $s_m$ in \eqref{eq-stage-coordinate-map} by applying
MV-negation to its first $m$ coordinates.  Since $s_m$ is continuous
by \eqref{eq-stage-initial-map} and MV-negation is continuous,
$\Delta_m$ is continuous.

Since every $W_k$ contains $0$, the inverse image of the box
$\prod_{k\in\mathbb Z}W_k$ under $\Omega_m$ is
$\prod_{k=-m}^{m-1}W_k$.
By \eqref{eq-coordinate-box-Q}, $Q(\mathcal W)$ is the inverse image
of the same box under $\delta$, so the diagram gives
\[
 Q(\mathcal W)\cap I_m
 =
 \Delta_m^{-1}\left(\prod_{k=-m}^{m-1}W_k\right).
\]
The finite product $\prod_{k=-m}^{m-1}W_k$ is open in
$(A^{2m},\tau^{2m})$, so its inverse image under the continuous map
$\Delta_m$ is open in $I_m$.
Since $m$ was arbitrary, $Q(\mathcal W)$ is a
$\lambda_A$-open zero-neighbourhood.

Finally, inversion is continuous on $(G,\lambda_A)$ by
Proposition~\ref{prop-lambda-quasitopological}.  Therefore
$-Q(\mathcal W)$ is also a $\lambda_A$-open zero-neighbourhood, and
\[
 B(\mathcal W)=Q(\mathcal W)\cap(-Q(\mathcal W))
\]
is a symmetric $\lambda_A$-open zero-neighbourhood.
\end{proof}

The next lemma relates symmetric coordinate boxes to group
topologies on the Mundici group for which the interval inclusion
is continuous, including all admissible $\ell$-group topologies.

\begin{lem}
\label{lem-coordinate-box-refinement}
Let $\theta$ be a group topology on $G$ for which the inclusion
$(A,\tau)\longrightarrow(G,\theta)$ is continuous.
Every $\theta$-zero-neighbourhood contains a symmetric coordinate box.
\end{lem}

\begin{proof}
Let $W$ be a $\theta$-zero-neighbourhood.  Choose symmetric
$\theta$-open zero-neighbourhoods $(V_j)_{j\geq0}$ such that
\[
 V_0\subseteq W,\qquad V_{j+1}+V_{j+1}\subseteq V_j.
\]
The finite-tail argument in the proof of
Lemma~\ref{lem-group-reflection-base} gives
\[
 V_1+\cdots+V_M\subseteq V_0\subseteq W
 \qquad(M\geq1).
\]
By symmetry and zero padding, every finite signed sum of elements
chosen from the $V_j$ with distinct indices $j\geq1$ lies in $W$.
Commutativity makes the order of summation irrelevant.
Fix a bijection $e\colon\mathbb Z\longrightarrow\mathbb N_{\geq1}$
and set $W_k=A\cap V_{e(k)}$.  By continuity of the interval
inclusion, every $W_k$ is a $\tau$-open zero-neighbourhood.
If $g\in Q((W_k)_k)$, then
$\delta_k(g)\in V_{e(k)}$ for every $k$.
By \eqref{eq-deviation-decomposition}, $g$ is a finite signed sum
of these deviations.  Symmetry absorbs the signs, and the indices
$e(k)$ are distinct.  Padding unused indices with zero therefore gives
$g\in V_1+\cdots+V_M\subseteq W$ for some $M$.
Thus $B((W_k)_k)\subseteq Q((W_k)_k)\subseteq W$.
\end{proof}

We now restrict the finite-sum construction to coordinate boxes.
Let $\mathcal B_\square$ consist of all sets
\begin{equation}
\label{eq-coordinate-box-group-base}
 O(B_1,B_2,\ldots)=\bigcup_{n\geq1}(B_1+\cdots+B_{2n}),
\end{equation}
where each $B_i$ is a symmetric coordinate box of the form
\eqref{eq-coordinate-box-B}.

\begin{prop}
\label{prop-coordinate-box-topology}
The family $\mathcal B_\square$ is a zero-neighbourhood base for
$\gamma_A$.
\end{prop}

\begin{proof}
Every symmetric coordinate box is a symmetric $\lambda_A$-open
zero-neighbourhood by Lemma~\ref{lem-deviation-decomposition}.
Thus Lemma~\ref{lem-group-reflection-base} makes each member of
$\mathcal B_\square$ a $\gamma_A$-zero-neighbourhood.

Conversely, let $W$ be a $\gamma_A$-zero-neighbourhood.
Choose symmetric $\gamma_A$-open zero-neighbourhoods $V_i$ such that
$V_1+\cdots+V_M\subseteq W$ for every $M\geq1$, as in the preceding
proof.  The inclusion $(A,\tau)\longrightarrow(G,\gamma_A)$ is
continuous because $\gamma_A\subseteq\lambda_A$ and
$\lambda_A|_A=\tau$.  Lemma~\ref{lem-coordinate-box-refinement}
therefore supplies a symmetric coordinate box $B_i\subseteq V_i$
for every $i$.  Consequently,
$O(B_1,B_2,\ldots)\subseteq W$.
This proves the required cofinality.
\end{proof}

\subsection{The bounded-output lemma}
\label{subsec-group-reflection}

To show that $\gamma_A$ induces $\tau_m$ on each $I_m$ and makes
the lattice operations continuous, we shall need the following
estimate.  Fix $m\geq1$, a finite $\ell$-group term
$t(x_1,\ldots,x_r)$ with fixed parameters from $G$ and
$t(0,\ldots,0)=0$, and an open zero-neighbourhood $V$ in $(I_m,\tau_m)$.
We seek $\gamma_A$-zero-neighbourhoods $N_1,\ldots,N_r$ such that
\[
 t(N_1\times\cdots\times N_r)\cap I_m\subseteq V.
\]

We first choose open zero-neighbourhoods $W_q$ in $A$ so that
$t(g_1,\ldots,g_r)\in I_m$ implies $t(g_1,\ldots,g_r)\in V$
whenever the inputs are finite signed sums of elements
$a_q\in W_q$, with each index $q$ used at most once among all the sums.
We then use coordinate boxes to construct the neighbourhoods $N_i$
so that every tuple in $N_1\times\cdots\times N_r$ admits such
representations.  Applying the bounded-output lemma to these representations
gives the desired inclusion
(Lemma~\ref{prop-gamma-bounded-control}\ref{gamma-bounded-output}).

\begin{lem}[Bounded-output lemma]
\label{lem-strong-FBOT}
Let $r\geq1$, and let $t(x_1,\ldots,x_r)$ be a finite term in the language
of abelian $\ell$-groups, with finitely many fixed parameters from $G$,
such that
\begin{equation}
\label{eq-FBOT-zero}
 t(0,\ldots,0)=0.
\end{equation}
Fix $m\geq1$ and an open neighbourhood $V$ of $0$ in
$(I_m,\tau_m)$.  Then there is a sequence $(W_n)_{n\geq1}$ of 
open zero-neighbourhoods in $A$ with the following property:
\begin{itemize}
\item For any pairwise disjoint finite sets $E_1,\ldots,E_r$
of positive integers and any choices of $\epsilon_q\in\{-1,1\}$
and $a_q\in W_q$ for $q\in E_1\cup\cdots\cup E_r$, put
\[
 g_i=\sum_{q\in E_i}\epsilon_q a_q
 \qquad(1\leq i\leq r).
\]
Then
\[
 t(g_1,\ldots,g_r)\in I_m
 \implies
 t(g_1,\ldots,g_r)\in V.
\]
\end{itemize}
\end{lem}

\begin{proof}
\textbf{Step 1.} We first control truncated sums with uniformly
bounded integer coefficients.
Recall the truncation maps $r_M$ from \eqref{eq-truncation-map}.
We claim that, for all integers $M,C\geq1$ and every open
zero-neighbourhood $H$ in $(I_M,\tau_M)$, there is a sequence
$(D_q)_{q\geq1}$ of open zero-neighbourhoods in $A$ such that
\[
 r_M\left(\sum_{q\in E}k_q a_q\right)\in H
\]
for every finite set $E$ of positive integers, every choice of
$a_q\in D_q$, and all integers $k_q$ with $|k_q|\leq C$ $(q\in E)$.
In this claim, the sequence $(D_q)$ is chosen once $M$, $C$,
and $H$ are fixed.  It must then work for every finite set $E$ of positive
integers and every allowed choice of coefficients $k_q$ and elements $a_q$,
even though the sums before truncation need not lie in any fixed bounded stage.

To prove the claim, fix $M$, $C$, and $H$, put $S_0=H$, and set
\begin{equation}
\label{eq-FBOT-stage-radii}
 M_j=M+Cj
 \qquad(j\geq0).
\end{equation}
For $j\geq1$ and each integer $k$ with $|k|\leq C$, define
\begin{equation}
\label{eq-FBOT-transition-map}
 F_{k,j}\colon
 (A,\tau)\times(I_{M_j},\tau_{M_j})
 \longrightarrow (I_{M_{j-1}},\tau_{M_{j-1}}),
 \qquad
 F_{k,j}(a,z)=r_{M_{j-1}}(z+ka).
\end{equation}
These maps are continuous by Lemma~\ref{lem-stage-topology}.
Indeed, $a\mapsto ka$ sends $A$ continuously into $I_C$ by
repeated bounded-stage addition, inversion, and inclusion
(and is constant when $k=0$).  Addition sends
$I_{M_j}\times I_C$ continuously into $I_{M_j+C}$, and
truncation from that bounded stage to $I_{M_{j-1}}$ is continuous.
Moreover, $F_{k,j}(0,0)=0$ for every such $k$.

Suppose that an open zero-neighbourhood $S_{j-1}$ in
$I_{M_{j-1}}$ has been chosen.
Continuity at $(0,0)$ supplies open zero-neighbourhoods in $A$
and $I_{M_j}$ whose product is mapped into $S_{j-1}$.
Taking intersections over the finitely many integers $k$ with
$|k|\leq C$, we obtain open zero-neighbourhoods
$U_j\subseteq A$ and $S_j\subseteq I_{M_j}$ such that
\begin{equation}
\label{eq-FBOT-recursive-control}
 F_{k,j}(U_j\times S_j)\subseteq S_{j-1}
 \qquad(k\in\mathbb Z,\ |k|\leq C).
\end{equation}
This recursively defines $U_j$ and $S_j$ for every $j\geq1$.
Set
\[
 D_q=\bigcap_{j=1}^q U_j
 \qquad(q\geq1).
\]
Each $D_q$ is an open zero-neighbourhood in $A$.

We next verify the claimed estimate for every finite sum
with these coefficient bounds.
The empty sum is $0$, so assume that $E\neq\varnothing$ and write
\[
 E=\{q_1<\cdots<q_d\}.
\]
For $0\leq j\leq d$, put
\[
 y^j=\sum_{\ell=j+1}^d k_{q_\ell}a_{q_\ell},
 \qquad z^j=r_{M_j}(y^j).
\]
We shall show, working backwards from $j=d$ to $j=0$, that
$z^j\in S_j$.
The initial case is $z^d=0\in S_d$.
For $1\leq j\leq d$, we have
$y^{j-1}=y^j+k_{q_j}a_{q_j}$.
Since $k_{q_j}a_{q_j}\in I_C$ and $M_j=M_{j-1}+C$,
Lemma~\ref{lem-truncation-identity} gives
\begin{equation}
\label{eq-FBOT-backward-step}
 z^{j-1}
 =r_{M_{j-1}}\bigl(z^j+k_{q_j}a_{q_j}\bigr)
 =F_{k_{q_j},j}(a_{q_j},z^j).
\end{equation}
Since the labels are distinct positive integers, $q_j\geq j$.
Thus $a_{q_j}\in D_{q_j}\subseteq U_j$, so
\eqref{eq-FBOT-recursive-control} carries $z^j\in S_j$
to $z^{j-1}\in S_{j-1}$.
Backwards induction therefore yields
\[
 r_M\left(\sum_{q\in E}k_q a_q\right)=z^0\in S_0=H.
\]
This proves the claim.

\textbf{Step 2.} We deduce the assertion for the general term $t$
from Step~1.
We express $t$ as a lattice term in finitely many affine sums,
then use Step~1 to control their linear parts in one common
neighbourhood.

We first separate the group operations from the lattice operations.
We claim that every abelian $\ell$-group
term in the variables $x_1,\ldots,x_r$, with fixed parameters from $G$,
can be written in the form
\begin{equation}
\label{eq-FBOT-normal-form}
 t(x_1,\ldots,x_r)
 =
 P\bigl(c_1+v_1(x),\ldots,c_s+v_s(x)\bigr),
 \qquad
 v_h(x)=\sum_{i=1}^r n_{hi}x_i,
\end{equation}
where $P$ is a finite lattice term involving only $\vee$ and $\wedge$,
the $c_h$ are fixed elements of $G$, and $n_{hi}\in\mathbb Z$.

To prove the claim, let $\mathcal C$ be the class of terms that admit
such a representation.  Each variable belongs to $\mathcal C$, since
$x_i$ itself has the required form, and each fixed parameter belongs to
$\mathcal C$ as well.  The class $\mathcal C$ is plainly closed under
$\vee$ and $\wedge$.

It is also closed under negation.  Indeed,
\[
 -(a\vee b)=(-a)\wedge(-b),
 \qquad
 -(a\wedge b)=(-a)\vee(-b),
\]
so a negation can be pushed through a finite lattice term until it reaches
the expressions at its leaves, and
\[
 -\left(c+\sum_i n_i x_i\right)
 =
 -c+\sum_i(-n_i)x_i
\]
is again of the required form.

Finally, $\mathcal C$ is closed under addition.  Since translations in an
$\ell$-group are lattice automorphisms,
\[
 a+(b\vee c)=(a+b)\vee(a+c),
 \qquad
 a+(b\wedge c)=(a+b)\wedge(a+c),
\]
and, by commutativity, the analogous identities hold in the first
variable.  Repeatedly applying these identities expands the sum of two
finite lattice terms into a finite lattice term whose leaves are sums of
expressions of the form $c+\sum_i n_i x_i$.
Such sums are again of the same form.  

Hence, $\mathcal C$ is closed under
all the term-forming operations, and the claim follows by structural
induction on terms.

To apply Step~1 to all the linear parts, put
\[
 C=\max\bigl\{1,\max_{h,i}|n_{hi}|\bigr\}.
\]
Choose an integer $R\geq1$ such that $c_h\in I_R$ for every
$1\leq h\leq s$, and put $L=m+R$.
By Lemma~\ref{lem-truncation-identity}, this choice gives
\[
 r_m(c_h+v)=r_m\bigl(c_h+r_L(v)\bigr)
 \qquad(v\in G,\ 1\leq h\leq s).
\]
Thus the truncated sums $r_L(v)$ retain all the information
needed after adding the fixed parameters and truncating to $I_m$.

Truncation also preserves finite joins and meets: by distributivity,
both $z\mapsto z\wedge mu$ and $z\mapsto(-mu)\vee z$ do so,
and hence their composite $r_m$ does so as well.
Define
\[
 \Phi\colon(I_L^s,\tau_L^s)\longrightarrow(I_m,\tau_m),
 \qquad
 \Phi(z_1,\ldots,z_s)
 =r_mP(c_1+z_1,\ldots,c_s+z_s).
\]
This map is continuous by Lemma~\ref{lem-stage-topology}:
each translation by $c_h$ sends $I_L$ continuously into $I_{L+R}$,
the lattice term $P$ is continuous on $I_{L+R}^s$ with values in
$I_{L+R}$, and truncation to $I_m$ is continuous.
Since $P(c_1,\ldots,c_s)=t(0,\ldots,0)=0$, we have
\begin{equation}
\label{eq-FBOT-initial-control}
 \Phi(0,\ldots,0)
 =r_mP(c_1,\ldots,c_s)
 =0.
\end{equation}
Continuity at zero gives open zero-neighbourhoods
$H_1,\ldots,H_s$ in $I_L$ with
$\Phi(H_1\times\cdots\times H_s)\subseteq V$.
Their intersection $H=\bigcap_{h=1}^sH_h$ is an open
zero-neighbourhood in $I_L$ satisfying
\begin{equation}
\label{eq-FBOT-target-neighbourhood}
 \Phi(H^s)\subseteq V.
\end{equation}
Apply Step~1 with radius $L$, coefficient bound $C$, and target $H$,
and denote the resulting sequence in $A$ by $(W_q)_{q\geq1}$.

We verify that this sequence has the property in the statement.
Fix the sets $E_1,\ldots,E_r$, the pieces $a_q\in W_q$, and the
signs $\epsilon_q$ defining $g=(g_1,\ldots,g_r)$.
For each $1\leq h\leq s$,
\[
 v_h(g)=\sum_{i=1}^r\sum_{q\in E_i}n_{hi}\epsilon_q a_q.
\]
Because the $E_i$ are pairwise disjoint, each label $q$ occurs
only once in this sum, with an integer coefficient of absolute
value at most $C$.  Step~1 therefore gives
\[
 z_h:=r_L(v_h(g))\in H
 \qquad(1\leq h\leq s).
\]
The same sequence $(W_q)$ works for every $h$, since the estimate
in Step~1 is uniform over all coefficients with absolute value
at most $C$.

Finally, Lemma~\ref{lem-truncation-identity} gives
$r_m(c_h+z_h)=r_m(c_h+v_h(g))$.
Using \eqref{eq-FBOT-normal-form} and the fact that $r_m$
preserves finite joins and meets, we obtain
\[
\begin{aligned}
 \Phi(z_1,\ldots,z_s)
 &=P\bigl(r_m(c_1+z_1),\ldots,r_m(c_s+z_s)\bigr)\\
 &=P\bigl(r_m(c_1+v_1(g)),\ldots,r_m(c_s+v_s(g))\bigr)
 =r_mt(g).
\end{aligned}
\]
Since $(z_1,\ldots,z_s)\in H^s$, this value belongs to $V$ by
\eqref{eq-FBOT-target-neighbourhood}.
If $t(g)\in I_m$, then $r_mt(g)=t(g)$, and hence $t(g)\in V$,
as required.
\end{proof}

\begin{rem}
\label{rem-FBOT-box-continuity}
Lemma~\ref{lem-strong-FBOT} admits an interpretation in terms of
continuity at zero.  Fix $t$ and $m$ as in that lemma.
For each choice of the sets $E_1,\ldots,E_r$ and signs
$\epsilon_q$ allowed there, put $E=E_1\cup\cdots\cup E_r$
and define
\[
 \Upsilon\colon A^E\longrightarrow G,
 \qquad
 \Upsilon((a_q)_{q\in E})
 =t\left(\sum_{q\in E_1}\epsilon_q a_q,\ldots,
          \sum_{q\in E_r}\epsilon_q a_q\right).
\]
For every open zero-neighbourhood $V$ in $(I_m,\tau_m)$, the lemma
provides a basic open zero-neighbourhood
$W=\prod_{n\geq1}W_n$ in $A^{\mathbb N_{\geq1}}$ with the box topology
such that, simultaneously for all the above choices,
\[
 \left(\prod_{q\in E}W_q\right)\cap\Upsilon^{-1}(I_m)
 \subseteq\Upsilon^{-1}(V).
\]
Since $\Upsilon(0)=0$, this is the neighbourhood condition for continuity
at zero of the restriction
\[
 \left.\Upsilon\right|_{\Upsilon^{-1}(I_m)}\colon
 \Upsilon^{-1}(I_m)\longrightarrow(I_m,\tau_m),
\]
where the domain
carries the subspace topology from $(A^E,\tau^E)$.
For any fixed choice of $E_1,\ldots,E_r$ and the signs, this
continuity already follows from Lemma~\ref{lem-stage-topology}:
$\Upsilon$ has bounded range and is built from operations on bounded stages.
The additional assertion of Lemma~\ref{lem-strong-FBOT} is that the
same box $W$ works for all these choices.

Every particular sum uses only finitely many labels, but there is
no fixed finite set containing all the labels that may occur.
The box topology allows a separate zero-neighbourhood for every
potential summand.  Lemma~\ref{lem-strong-FBOT} chooses these
neighbourhoods so that, however many summands occur, the value of
$t$ belongs to $V$ whenever it lies in $I_m$; the input sums
themselves need not lie in a common bounded stage.
\end{rem}

\begin{samepage}
We now construct $\gamma_A$-zero-neighbourhoods from an arbitrary
sequence $(W_q)$ of open zero-neighbourhoods in $A$.

\begin{lem}
\label{prop-gamma-bounded-control}
Let $r\geq1$.
\begin{enumerate}[label=(\alph*),ref=(\alph*)]
\item\label{gamma-labelled-representations}
Given a sequence $(W_q)_{q\geq1}$ of open zero-neighbourhoods in $A$,
there are zero-neighbourhoods $N_1,\ldots,N_r$ in $(G,\gamma_A)$
such that, whenever $g_i\in N_i$ for $1\leq i\leq r$, there are representations
\[
 g_i=\sum_{q\in E_i}\epsilon_q a_q
 \qquad(1\leq i\leq r),
\]
where $E_1,\ldots,E_r$ are pairwise disjoint finite sets of positive
integers, $\epsilon_q\in\{-1,1\}$, and $a_q\in W_q$.
\item\label{gamma-bounded-output}
Let $t=t(x_1,\ldots,x_r)$, $m$, and $V$ be as in
Lemma~\ref{lem-strong-FBOT}.  Then there are
zero-neighbourhoods $N_1,\ldots,N_r$ in $(G,\gamma_A)$ such that
\[
 t(N_1\times\cdots\times N_r)\cap I_m\subseteq V.
\]
\end{enumerate}
In both parts, the $N_i$ may be chosen in the form
\[
 N_i=O(C_{i,1},C_{i,2},\ldots),
\]
where every $C_{i,j}$ is a symmetric coordinate box
$B(\mathcal W^{i,j})$.
\end{lem}
\end{samepage}

\begin{proof}
\ref{gamma-labelled-representations}
Fix a sequence $(W_q)_{q\geq1}$ of open zero-neighbourhoods in $A$.
We construct the $N_i$ and then verify the required representations.
Choose an injection
\[
 e\colon\{1,\ldots,r\}\times\mathbb N_{\geq1}\times\mathbb Z
 \longrightarrow\mathbb N_{\geq1}.
\]
Thus $e(i,j,k)$ provides a distinct label for every triple $(i,j,k)$;
below, these indices will record, respectively, the variable, the summand
occurring in that variable, and the Mundici deviation of that summand.

For $1\leq i\leq r$ and $j\geq1$, define
\[
 \mathcal W^{i,j}
 =\bigl(W_{e(i,j,k)}\bigr)_{k\in\mathbb Z},\qquad
 C_{i,j}=B(\mathcal W^{i,j}),\qquad
 N_i=O(C_{i,1},C_{i,2},\ldots).
\]
By Lemma~\ref{lem-deviation-decomposition}, each $C_{i,j}$ is a
symmetric $\lambda_A$-open zero-neighbourhood.  Hence
Lemma~\ref{lem-group-reflection-base} shows that every $N_i$ is a
$\gamma_A$-zero-neighbourhood.

We next express arbitrary elements of the $N_i$ as labelled
signed interval sums.  Choose
$g_i\in N_i$ for $1\leq i\leq r$.  By
\eqref{eq-group-reflection-base}, for each $i$ there are $n_i\geq1$ and
elements $c_{i,j}\in C_{i,j}$ $(1\leq j\leq2n_i)$ such that
\begin{equation}
\label{eq-gamma-bounded-sum}
 g_i=\sum_{j=1}^{2n_i}c_{i,j}.
\end{equation}
For every pair $(i,j)$ occurring here, put
\[
 K_{i,j}=\{k\in\mathbb Z:\delta_k(c_{i,j})\neq0\}.
\]
This set is finite by Lemma~\ref{lem-deviation-decomposition}.  Define
\[
 E_i
 =\{e(i,j,k):1\leq j\leq2n_i,\ k\in K_{i,j}\}.
\]
The sets $E_i$ are finite, and they are pairwise disjoint because
$e$ is injective.

For $q=e(i,j,k)\in E_i$, set
\[
 a_q=\delta_k(c_{i,j}),
 \qquad
 \epsilon_q=
 \begin{cases}
  1,&k\geq0,\\
  -1,&k<0.
 \end{cases}
\]
Since
$c_{i,j}\in C_{i,j}=B(\mathcal W^{i,j})\subseteq
Q(\mathcal W^{i,j})$, the definition
\eqref{eq-coordinate-box-Q} gives
\[
 a_q=\delta_k(c_{i,j})
 \in W_{e(i,j,k)}=W_q.
\]
Using \eqref{eq-deviation-decomposition} for each $c_{i,j}$ and then
\eqref{eq-gamma-bounded-sum}, we obtain
\[
 g_i=\sum_{q\in E_i}\epsilon_q a_q
 \qquad(1\leq i\leq r).
\]
This proves \ref{gamma-labelled-representations}.

\ref{gamma-bounded-output}
Take the sequence $(W_q)$ from Lemma~\ref{lem-strong-FBOT}, and choose
$N_1,\ldots,N_r$ as in \ref{gamma-labelled-representations}.
For any tuple $(g_1,\ldots,g_r)\in N_1\times\cdots\times N_r$,
part~\ref{gamma-labelled-representations} provides representations
with pairwise disjoint label sets and $a_q\in W_q$.
If $t(g_1,\ldots,g_r)\in I_m$, Lemma~\ref{lem-strong-FBOT}
therefore gives $t(g_1,\ldots,g_r)\in V$.  Thus
\[
 t(N_1\times\cdots\times N_r)\cap I_m\subseteq V,
\]
as required.
\end{proof}

Write $g^+=g\vee0$ for the positive part of $g\in G$.
We now control how the positive part changes when a finite signed sum
$h$ is added to a fixed element $y\in G$, that is, the difference
\[
 (y+h)^+-y^+.
\]

\begin{lem}
\label{lem-tail-telescoping}
Fix $y\in G$, and let $(V_r)_{r\geq1}$ be a sequence of open
zero-neighbourhoods in $(I_1,\tau_1)$.  
There is a sequence $(W_r)_{r\geq1}$ of open zero-neighbourhoods
in $(A,\tau)$ such that, for every finite subset
$E\subseteq\mathbb N_{\geq1}$ and every choice of
$a_r\in W_r$ and $\epsilon_r\in\{-1,1\}$ $(r\in E)$, if
\[
 h=\sum_{r\in E}\epsilon_r a_r,
\]
then there are elements $d_r\in V_r$ $(r\geq1)$, only finitely many of
them nonzero, such that
\[
 (y+h)^+-y^+=\sum_{r\geq1}d_r.
\]
\end{lem}

\begin{proof}
We first isolate the change in positive part produced by adding one piece
while the pieces with larger labels are kept fixed.  For this purpose,
fix the finite two-variable $\ell$-group term with parameter $y$:
\[
 \Psi_y\colon G^2\longrightarrow G,
 \qquad
 \Psi_y(x,z)=(y+x+z)^+-(y+z)^+.
\]
Here $x$ represents the piece currently being added, while $z$ represents
the sum of the pieces with larger labels.  Thus $\Psi_y(x,z)$ measures the
increment in positive part caused by adding $x$ after the contribution
$z$ has already been included.  Clearly, $\Psi_y(0,0)=0$. 

In any abelian $\ell$-group, the positive-part map satisfies
\begin{equation}
\label{eq-positive-part-contraction}
 |a^+-b^+|\leq|a-b|
\end{equation}
for all $a,b\in G$.  Indeed,
\[
 a=b+(a-b)\leq b^++(a-b)^+.
\]
Since the right-hand side is nonnegative, the definition
$a^+=a\vee0$ gives $a^+\leq b^++(a-b)^+$, and hence 
\[
 a^+-b^+\leq(a-b)^+.
\]
Interchanging $a$ and $b$ yields $b^+-a^+\leq(b-a)^+$.
Using the standard facts for $\ell$-groups that
\[
 x^+\leq |x|,
 \quad
 |-x|=|x|,
 \quad\text{and}\quad
 |x|=x\vee(-x),
\]
we obtain
\[
 |a^+-b^+|
 =
 (a^+-b^+)\vee(b^+-a^+)
 \leq(a-b)^+\vee(b-a)^+
 \leq |a-b|,
\]
which proves \eqref{eq-positive-part-contraction}.

In particular, by \eqref{eq-positive-part-contraction}, for
$a\in A=[0,u]$, $\epsilon\in\{-1,1\}$, and $z\in G$, we have
\[
 |\Psi_y(\epsilon a,z)|
 \leq |\epsilon a|
 =a
 \leq u.
\]
Thus $\Psi_y(\epsilon a,z)\in I_1$.
Hence, although $\Psi_y$ need not have bounded range on all of $G^2$, its
value lies in $I_1$ whenever its first variable is a signed element
$\epsilon a$ with $a\in A$.  This is exactly the bounded-output situation
needed below to apply
Lemma~\ref{prop-gamma-bounded-control}\ref{gamma-bounded-output}.

To control both the current piece and its tail, we choose
symmetric open zero-neighbourhoods $H_q$ in $(G,\gamma_A)$ such that
\[
 \Psi_y(H_q\times H_q)\cap I_1\subseteq V_q,
 \qquad H_{q+1}+H_{q+1}\subseteq H_q
 \qquad(q\geq1).
\]
For each $q\geq1$,
Lemma~\ref{prop-gamma-bounded-control}\ref{gamma-bounded-output}, applied
to $\Psi_y$ with target stage $I_1$ and neighbourhood $V_q$,
provides two input zero-neighbourhoods.
Choose $H_q$ recursively inside both of them; for $q>1$,
shrink it further so that $H_q+H_q\subseteq H_{q-1}$.
This is possible because $\gamma_A$ is a group topology, and
shrinking preserves the first inclusion.
Put $W_q=A\cap H_q$.  These are open zero-neighbourhoods in
$(A,\tau)$ because $\gamma_A\subseteq\lambda_A$ and
$\lambda_A|_A=\tau$.

If $E=\varnothing$, then $h=0$, and the conclusion follows by taking
$d_r=0$ for every $r\geq1$.  Assume therefore that $E\neq\varnothing$,
and let $N=\max E$. 
It is convenient to fill the gaps in the set of labels.  
For every $1\leq r\leq N$ with $r\notin E$, set $a_r=0$ and $\epsilon_r=1$.
Since $0\in W_r$, the inserted zero pieces also satisfy $a_r\in W_r$.  
We may therefore write
\[
 h=\sum_{r=1}^N\epsilon_r a_r.
\]

We now control the tails before adding the pieces in decreasing
order of their labels.  For
$1\leq q\leq N$, set
\[
 \beta_q
 =
 \sum_{r=q}^N\epsilon_r a_r,
\]
and put $\beta_{N+1}=0$.  We claim that $\beta_{q+1}\in H_q$ for every
$1\leq q\leq N$.  For $q=N$, this holds because $\beta_{N+1}=0$.
If $q<N$ and $\beta_{q+2}\in H_{q+1}$, then symmetry of $H_{q+1}$
and $a_{q+1}\in W_{q+1}\subseteq H_{q+1}$ give
\[
 \beta_{q+1}=\epsilon_{q+1}a_{q+1}+\beta_{q+2}
 \in H_{q+1}+H_{q+1}\subseteq H_q.
\]
Descending induction proves the claim.

The step from $y+\beta_{q+1}$ to $y+\beta_q$ adds the piece
$\epsilon_q a_q$.  Define its positive-part increment by
\[
 d_q
 =(y+\beta_q)^+-(y+\beta_{q+1})^+
 =\Psi_y(\epsilon_q a_q,\beta_{q+1})
 \qquad(1\leq q\leq N),
\]
and put $d_q=0$ for $q>N$.

For $1\leq q\leq N$, symmetry and $a_q\in W_q\subseteq H_q$
give $\epsilon_q a_q\in H_q$, while the tail estimate gives
$\beta_{q+1}\in H_q$.
The contraction estimate above also gives $d_q\in I_1$, so
\[
 d_q\in\Psi_y(H_q\times H_q)\cap I_1\subseteq V_q.
\]
For $q>N$, we also have $d_q=0\in V_q$.

Only $d_1,\ldots,d_N$ can be nonzero.
Since $\beta_1=h$ and $\beta_{N+1}=0$, cancelling the intermediate terms gives
\[
 \sum_{q\geq1}d_q
 =
 \sum_{q=1}^N\bigl((y+\beta_q)^+-(y+\beta_{q+1})^+\bigr)
 =
 (y+h)^+-y^+.
\]
Thus the sequence $(d_q)_{q\geq1}$ has all the required properties.
\end{proof}

\subsection{The topological Mundici representation theorem}
\label{subsec-exact-lgroup}

We now apply the estimates to the group reflection $\gamma_A$.
First we prove $\kappa_A\subseteq\gamma_A$ using the initial
description of $\kappa_A$ by the bounded coordinates.
The comparison $\kappa_A\subseteq\gamma_A\subseteq\lambda_A$
then shows that $\gamma_A$ induces $\tau_n$ on every bounded stage $I_n$.
Next the reverse-tail estimate gives continuity of positive part,
and hence joint lattice continuity.  The resulting $\ell$-group
extension is identified with the Mundici topology $\sigma_A$.

\begin{prop}
\label{prop-kappa-rho-gamma}
One has
\[
 \kappa_A\subseteq\gamma_A.
\]
Moreover, the group reflection preserves every bounded-stage topology:
\begin{equation}
\label{eq-gamma-stage-exactness}
 \gamma_A|_{I_n}=\tau_n
 \qquad(n\geq1).
\end{equation}
In particular,
\[
 \kappa_A|_A=\gamma_A|_A=\tau,
\]
so $\gamma_A$ is a group topology inducing $\tau$ on $A$.
\end{prop}

\begin{proof}
By definition, $\gamma_A\subseteq\lambda_A$.
To prove $\kappa_A\subseteq\gamma_A$, we use the coordinate-initial
description in Proposition~\ref{prop-kappa-lambda}.

Fix $k\in\mathbb Z$, $x\in G$, and an open neighbourhood $V$
of $a=\pi_k(x)$ in $(A,\tau)$.
Since $\lambda_A|_A=\tau$ and translations are
$\lambda_A$-homeomorphisms, $V-a$ is open in the translated
subspace $A-a\subseteq I_1$.
By $\lambda_A|_{I_1}=\tau_1$, choose an open zero-neighbourhood
$W$ in $(I_1,\tau_1)$ such that
\[
 W\cap(A-a)=V-a.
\]
The centred $\ell$-group term
$\psi_{k,x}(h)=\pi_k(x+h)-a$ vanishes at zero and has range in
$A-a\subseteq I_1$.
Lemma~\ref{prop-gamma-bounded-control}\ref{gamma-bounded-output} therefore supplies a
$\gamma_A$-zero-neighbourhood $N$ with $\psi_{k,x}(N)\subseteq W$.
Since its values also lie in $A-a$, we obtain
$\pi_k(x+N)\subseteq V$.
Thus each $\pi_k$ is $\gamma_A$-continuous, proving
$\kappa_A\subseteq\gamma_A$ by initiality.

For every $n\geq1$, the resulting sandwich
\[
 \kappa_A|_{I_n}\subseteq\gamma_A|_{I_n}
 \subseteq\lambda_A|_{I_n}
\]
has both endpoints equal to $\tau_n$ by
\eqref{eq-lambda-stage-restriction}.
This proves \eqref{eq-gamma-stage-exactness}; restricting to
$A\subseteq I_1$ gives $\gamma_A|_A=\tau$.
\end{proof}

The group topology $\gamma_A$ now induces $\tau$ on $A$.
It remains to prove joint continuity of the lattice operations.
For each fixed $y\in G$ and each basic zero-neighbourhood
$O=O(B_1,B_2,\ldots)$, we seek a $\gamma_A$-zero-neighbourhood $N$ such that
\[
 (y+N)^+-y^+\subseteq O.
\]
We will use Lemma~\ref{lem-tail-telescoping} to express each difference
as a finite sum whose summands lie in the prescribed boxes $B_r$.
This will prove continuity of positive part; the standard $\ell$-group
formulas will then give joint continuity of joins and meets.

\begin{thm}
\label{thm-rho-lgroup}
The positive-part map
\begin{equation}
\label{eq-positive-part-map-rho-square}
 (G,\gamma_A)\longrightarrow(G,\gamma_A),
 \qquad g\mapsto g^+
\end{equation}
is continuous.  Consequently, $\gamma_A$ is an $\ell$-group
extension of $\tau$ on $G$.
\end{thm}

\begin{proof}
We prove continuity of the positive-part map at an arbitrary $y\in G$.
Let $O(B_1,B_2,\ldots)$ be a basic $\gamma_A$-zero-neighbourhood as in
\eqref{eq-coordinate-box-group-base}, and put
\begin{equation}
\label{eq-positive-part-control-neighbourhood}
V_r=B_r\cap I_1
\qquad(r\geq1).
\end{equation}
Each $B_r$ is $\lambda_A$-open by
Lemma~\ref{lem-deviation-decomposition}, while
Proposition~\ref{prop-lambda-exact-space} gives
$\lambda_A|_{I_1}=\tau_1$.  Hence $(V_r)_{r\geq1}$ is a sequence of
open zero-neighbourhoods in $(I_1,\tau_1)$.  
Choose a sequence $(W_r)_{r\geq1}$ of open zero-neighbourhoods in
$(A,\tau)$ as supplied by Lemma~\ref{lem-tail-telescoping}.

Apply Lemma~\ref{prop-gamma-bounded-control}\ref{gamma-labelled-representations}
with one variable and this
sequence $(W_r)_{r\geq1}$.
It gives a basic $\gamma_A$-zero-neighbourhood
$N=O(C_1,C_2,\ldots)$ such that every $h\in N$ admits a representation
\begin{equation}
\label{eq-source-neighbourhood-labeling}
 h=\sum_{r\in E}\epsilon_r a_r,
 \qquad a_r\in W_r,\quad\epsilon_r\in\{-1,1\},
\end{equation}
where $E$ is a finite set of positive integers with distinct labels.
Lemma~\ref{lem-tail-telescoping} supplies $d_r\in V_r\subseteq B_r$,
only finitely many nonzero, such that
\[
 (y+h)^+-y^+=\sum_{r\geq1}d_r.
\]
Choose $p\geq1$ so that $d_r=0$ for $r>2p$.
Padding with these zero terms gives
\[
 (y+h)^+-y^+
 =d_1+\cdots+d_{2p}
 \in B_1+\cdots+B_{2p}
 \subseteq O(B_1,B_2,\ldots).
\]
Consequently,
\[
 (y+N)^+\subseteq y^++O(B_1,B_2,\ldots),
\]
which proves continuity at $y$.
Since $y$ was arbitrary, positive part is continuous on $G$.

It remains to recover the binary lattice operations.  The topology $\gamma_A$ is a group topology by construction,
so addition and inversion are continuous.  The standard
$\ell$-group identities
\begin{equation}
\label{eq-lattice-operations-positive-part}
 x\vee y=y+(x-y)^+,
 \qquad
 x\wedge y=x-(x-y)^+
\end{equation}
therefore show that both
\[
 \vee,\wedge\colon
 (G^2,\gamma_A^2)
 \longrightarrow
 (G,\gamma_A)
\]
are continuous.  Hence $\gamma_A$ is an $\ell$-group topology.
Finally, Proposition~\ref{prop-kappa-rho-gamma} gives
$\gamma_A|_A=\tau$, so $\gamma_A$ is an $\ell$-group
extension of $\tau$.
\end{proof}

\begin{thm}[Topological Mundici representation theorem]
\label{thm-topological-representation}
\label{thm-intro-representation}
For every topological MV-algebra $(A,\tau)$, write $\Xi A=(G_A,u_A)$.
Then the canonical Mundici isomorphism
\[
 \eta_A\colon(A,\tau)\longrightarrow
 \Gamma_t\Xi_t(A,\tau)
\]
is a homeomorphism, and hence an isomorphism of topological
MV-algebras.  In particular, the Mundici topology $\sigma_A$ is an
$\ell$-group topology satisfying $\sigma_A|_{[0,u_A]}=\tau$.
\end{thm}

\begin{proof}
By Theorem~\ref{thm-rho-lgroup} and
Proposition~\ref{prop-kappa-rho-gamma}, $\gamma_A$ is an
$\ell$-group topology inducing $\tau$ on $A$.
It is therefore admissible, so maximality of $\sigma_A$ gives
$\gamma_A\subseteq\sigma_A$.  Hence
\[
 \tau=\gamma_A|_A\subseteq\sigma_A|_A\subseteq\tau,
\]
where the last inclusion follows from admissibility of $\sigma_A$.
Thus $\sigma_A|_{[0,u_A]}=\tau$, so $\eta_A$ is a
homeomorphism and hence an isomorphism in $\TopMV$.

Moreover, Theorem~\ref{thm-lambda-envelope} gives
$\sigma_A\subseteq\lambda_A$.
Since $\sigma_A$ is a group topology by Lemma~\ref{lem-finest-admissible},
the defining maximality of $\gamma_A$ yields
$\sigma_A\subseteq\gamma_A$.
Thus
\begin{equation}
\label{eq-sigma-gamma}
 \sigma_A=\gamma_A.
\end{equation}
This completes the proof.
\end{proof}

\begin{cor}[Structure of the Mundici topology]
\label{cor-reconstruction-topology-chain}
For every $(A,\tau)\in\TopMV$, identify $A=[0,u]$ in $\Xi A=(G,u)$.
The Mundici topology $\sigma_A$ is the largest $\ell$-group extension
of $\tau$, and
\begin{equation}
\label{eq-topological-reconstruction-chain}
 \kappa_A\subseteq\sigma_A=\gamma_A
 =R_{\mathrm{grp}}(\lambda_A)\subseteq\lambda_A.
\end{equation}
All these topologies induce $\tau_n$ on every bounded stage $I_n$,
and hence induce $\tau$ on $A$.
Moreover, the family $\mathcal B_\square$ in
\eqref{eq-coordinate-box-group-base}, consisting of unions of finite
sums of symmetric coordinate boxes, is a zero-neighbourhood base
for $\sigma_A$.
\end{cor}

\begin{proof}
By Theorem~\ref{thm-topological-representation}, $\sigma_A$ is an
$\ell$-group extension of $\tau$.  Every other such extension is
admissible and hence contained in $\sigma_A$ by
Lemma~\ref{lem-finest-admissible}.
Equation~\eqref{eq-sigma-gamma} and the definition of $\gamma_A$ give
$\sigma_A=\gamma_A=R_{\mathrm{grp}}(\lambda_A)\subseteq\lambda_A$.
Propositions~\ref{prop-kappa-rho-gamma} and~\ref{prop-lambda-exact-space}
give the remaining inclusion and the bounded-stage restrictions.
Finally, the assertion about $\mathcal B_\square$ follows from
\eqref{eq-sigma-gamma} and
Proposition~\ref{prop-coordinate-box-topology}.
\end{proof}

\begin{rem}[Strictness of both inclusions]
\label{rem-reconstruction-strictness}
Both inclusions in \eqref{eq-topological-reconstruction-chain} can
be strict simultaneously.  Take the rational product cube of
Example~\ref{exmp-lambda-not-group}.  Its topology $\lambda_A$ is
not a group topology, whereas $\gamma_A=\sigma_A$ is one, so
$\sigma_A\subsetneq\lambda_A$.
For the other inclusion, put
\[
 p(x)=\sup_{j\geq1}\frac{|x_j|}{j}\qquad(x\in G).
\]
The metric $d_p(x,y)=p(x-y)$ defines an $\ell$-group topology
$\theta_p$: addition and inversion are continuous, and coordinatewise
join and meet satisfy
\[
 p\bigl((x\vee y)-(x'\vee y')\bigr)
 \leq\max\{p(x-x'),p(y-y')\},
\]
with the analogous inequality for meet.
On $A$, this topology is exactly the product topology.
Indeed, each evaluation is $p$-continuous, while
$\frac{|a_j-b_j|}{j}\leq1/j$ uniformly for $a,b\in A$, so a $p$-neighbourhood
can be obtained by restricting only finitely many coordinates.
Thus $\theta_p$ is an $\ell$-group extension of $\tau$ and
$\theta_p\subseteq\sigma_A$.
But $j\mathbf e_j\to0$ pointwise, hence in $\kappa_A$, whereas
$p(j\mathbf e_j)=1$.  Therefore
\begin{equation}
\label{eq-topological-chain-simultaneously-strict}
 \kappa_A\subsetneq\sigma_A=\gamma_A\subsetneq\lambda_A.
\end{equation}
The assertion is that strictness can occur, not that it occurs for
every input.  For the usual interval $A=[0,1]$, all these topologies
are the usual topology on $G=\mathbb R$: the usual bounded stages
cover $\mathbb R$ by their interiors, so their final topology is
usual, as is the coordinate topology.
\end{rem}

% The representation theorem identifies the formally defined maximal
% admissible topology with an explicit group reflection.
% The bounded-output and reverse-tail estimates prove that this
% reflection preserves the prescribed bounded stages and has joint
% lattice continuity.  The next section studies the counit and
% identifies the objects recovered by this reconstruction.

\section{The topological Mundici correspondence}
\label{sec-correspondence}

By Theorem~\ref{thm-topological-representation}, the unit of
$\Xi_t\dashv\Gamma_t$ is an isomorphism.  We now identify the
topological unital abelian $\ell$-groups for which the counit is
an isomorphism and obtain the topological Mundici correspondence.
The relevant condition is that the topology on the whole
group be the largest $\ell$-group topology inducing its given topology
on the unit interval.  
We conclude with retraction and separation-axiom consequences
of the representation theorem, together with the restriction
of the correspondence to Hausdorff objects.

Fix $(G,u,\rho)\in\Topulg$.  The counit of the topological adjunction
$\Xi_t\dashv\Gamma_t$ at $(G,u,\rho)$ is
\[
 \varepsilon_{(G,u,\rho)}\colon
 \Xi_t\Gamma_t(G,u,\rho)\longrightarrow(G,u,\rho).
\]
Its underlying unital $\ell$-group homomorphism is the algebraic
Mundici counit
\begin{equation}
\label{eq-transported-counit}
 \varepsilon_{(G,u)}\colon
 \Xi\Gamma(G,u)\longrightarrow(G,u),
\end{equation}
which is an isomorphism by \eqref{eq-Mundici-unit-counit}.
Put $\tau=\rho|_{[0,u]}$.  Transport the Mundici topology
$\sigma_{\Gamma_t(G,u,\rho)}$ along $\varepsilon_{(G,u)}$, and denote
the resulting $\ell$-group topology on $G$ by $\widetilde\sigma_G$.
Explicitly, for $U\subseteq G$,
\begin{equation}
\label{eq-transported-maximal-topology}
 U\in\widetilde\sigma_G
 \iff
 \varepsilon_{(G,u)}^{-1}(U)
 \in\sigma_{\Gamma_t(G,u,\rho)}.
\end{equation}
By construction,
\begin{equation}
\label{eq-transported-counit-homeomorphism}
 \varepsilon_{(G,u)}\colon
 \Xi_t\Gamma_t(G,u,\rho)\longrightarrow
 (G,u,\widetilde\sigma_G)
\end{equation}
is a homeomorphism, and hence an isomorphism in $\Topulg$.
The algebraic counit depends only on $(G,u)$, and
$\Gamma_t(G,u,\rho)=([0,u],\tau)$, so $\widetilde\sigma_G$ depends
on $\rho$ only through $\tau$.

\begin{prop}[Counit criterion]
\label{prop-transported-maximality}
Let $(G,u,\rho)\in\Topulg$, and put
$\tau=\rho|_{[0,u]}$.  Then
\begin{equation}
\label{eq-transported-interval-exactness}
 \widetilde\sigma_G|_{[0,u]}=\tau.
\end{equation}
Moreover, if $\theta$ is any $\ell$-group topology on $G$ such that
$\theta|_{[0,u]}\subseteq\tau$, then $\theta\subseteq\widetilde\sigma_G$.
Consequently, $\widetilde\sigma_G$ is the largest $\ell$-group topology
on $G$ inducing $\tau$ on $[0,u]$.
The counit is a homeomorphism precisely when $\rho$ is this largest
extension:
\begin{equation}
\label{eq-counit-Mundici-maximality}
 \varepsilon_{(G,u,\rho)}
 \text{ is a homeomorphism}
 \iff \rho=\widetilde\sigma_G.
\end{equation}
\end{prop}

\begin{proof}
Put $A=\Gamma(G,u)=[0,u]$, write $\Xi A=(G_A,u_A)$, and let
$j_G\colon A\longhookrightarrow G$ be the inclusion.
The algebraic triangle identity gives the commutative triangle
\[
\begin{tikzcd}[column sep=4.5em,row sep=3em]
A
  \arrow[r,hook,"\iota_A"]
  \arrow[dr,hook,"j_G"']
&
G_A
  \arrow[d,"\varepsilon_{(G,u)}"]
\\
&
G.
\end{tikzcd}
\]
Since $\varepsilon_{(G,u)}$ is a unital $\ell$-group isomorphism,
transport along it identifies the admissible $\ell$-group topologies
for $(A,\tau)$ with the $\ell$-group topologies $\theta$ on $G$
satisfying $\theta|_A\subseteq\tau$.
By Lemma~\ref{lem-finest-admissible}, $\widetilde\sigma_G$ is the
largest of these topologies.
In particular, $\widetilde\sigma_G|_A\subseteq\tau$, while
$\rho\subseteq\widetilde\sigma_G$ because $\rho|_A=\tau$.
Hence
\[
 \tau=\rho|_A\subseteq\widetilde\sigma_G|_A\subseteq\tau,
\]
which proves \eqref{eq-transported-interval-exactness} and the
largest-extension assertion.

The topological counit fits into the commutative triangle
\begin{equation}
\label{eq-transported-counit-factorization}
\begin{tikzcd}[column sep=5.5em,row sep=4em]
\Xi_t\Gamma_t(G,u,\rho)
  \arrow[r,"\varepsilon_{(G,u)}"]
  \arrow[dr,"\varepsilon_{(G,u,\rho)}"']
&
(G,u,\widetilde\sigma_G)
  \arrow[d,"1_G"]
\\
&
(G,u,\rho).
\end{tikzcd}
\end{equation}
The horizontal arrow is a homeomorphism by
\eqref{eq-transported-counit-homeomorphism}, and the vertical arrow is a
homeomorphism if and only if $\rho=\widetilde\sigma_G$.
This proves \eqref{eq-counit-Mundici-maximality}.
\end{proof}

We give a name to the objects satisfying this counit criterion.

\begin{defn}
An object $(G,u,\rho)\in\Topulg$ is \emph{Mundici-fine} if
\[
 \rho=\widetilde\sigma_G.
\]
Equivalently, $\rho$ is the finest $\ell$-group topology on $G$
inducing $\rho|_{[0,u]}$ on the unit interval.
\end{defn}

By \eqref{eq-counit-Mundici-maximality}, the Mundici-fine objects are
precisely those at which the counit is a homeomorphism.
The full subcategory of Mundici-fine objects is denoted by
\[
 \TopulgM.
\]

\begin{thm}[Topological Mundici correspondence]
\label{thm-topological-correspondence}
\label{thm-intro-correspondence}
The left adjoint
\[
 \Xi_t\colon\TopMV\longrightarrow\Topulg
\]
is fully faithful, and its essential image is precisely
$\TopulgM$.  Consequently, the adjunction of
Theorem~\ref{thm-topological-adjunction} restricts to an equivalence
\[
 \TopMV\simeq\TopulgM.
\]
Moreover, $\TopulgM$ is a replete coreflective full subcategory of
$\Topulg$, with coreflector given by the corestriction of
$\Xi_t\Gamma_t$ to $\TopulgM$.
\end{thm}

\begin{proof}
By Theorem~\ref{thm-topological-representation}, the unit
$\eta\colon1_{\TopMV}\Longrightarrow\Gamma_t\Xi_t$ is a natural
isomorphism.  Hence $\Xi_t$ is fully faithful by the standard
left-adjoint criterion
(see, e.g., \cite[Proposition~II.2.6.1]{Hofmann2014}).

By Corollary~\ref{cor-reconstruction-topology-chain}, every object
$\Xi_t(A,\tau)$ is Mundici-fine.
Transport along a topological unital $\ell$-group isomorphism
preserves this largest-extension property, so $\TopulgM$ is replete
and contains the essential image of $\Xi_t$.
Conversely, at a Mundici-fine object the counit is an isomorphism by
Proposition~\ref{prop-transported-maximality}, giving
\[
 (G,u,\rho)\cong\Xi_t\Gamma_t(G,u,\rho).
\]
Hence $\TopulgM$ is exactly the essential image of the fully
faithful functor $\Xi_t$, and its corestriction is the stated
equivalence, with the restricted unit and counit.

The standard coreflection consequence for fully faithful left
adjoints (see, e.g., \cite[Corollary~II.2.6.3]{Hofmann2014} for the
dual statement) now makes $\TopulgM$ coreflective.
The coreflector is the corestriction of $\Xi_t\Gamma_t$, with
coreflection morphism
\[
 \varepsilon_{(G,u,\rho)}\colon
 \Xi_t\Gamma_t(G,u,\rho)\longrightarrow(G,u,\rho).
 \qedhere
\]
\end{proof}

\begin{rem}
Theorem~\ref{thm-topological-representation} shows that the unit
of the topological Mundici adjunction is always an isomorphism.
By Proposition~\ref{prop-transported-maximality}, the counit is an
isomorphism exactly at the Mundici-fine objects, whose topology is
the largest $\ell$-group extension of the unit-interval topology.
Theorem~\ref{thm-topological-correspondence} therefore lifts the
classical Mundici equivalence topologically not to all
of $\Topulg$, but exactly to the replete coreflective full subcategory
$\TopulgM$.
The lifting diagram \eqref{eq-topological-lifting-square} now restricts
to the commutative diagram of equivalences
\begin{equation}
\label{eq-topological-equivalence-square}
\begin{tikzcd}[column sep=4.5em,row sep=3em]
\TopMV
  \arrow[r,shift left=1.4ex,"\Xi_t"]
  \arrow[r,phantom,"\simeq" description]
  \arrow[d,"U_{\mathrm{MV}}"']
&
\TopulgM
  \arrow[l,shift left=1.4ex,"\Gamma_t"]
  \arrow[d,"U_{\ell}"]
\\
\MV
  \arrow[r,shift left=1.4ex,"\Xi"]
  \arrow[r,phantom,"\simeq" description]
&
\uLG
  \arrow[l,shift left=1.4ex,"\Gamma"]
\end{tikzcd}
\end{equation}
Here $\Xi_t$ is corestricted to $\TopulgM$, while $\Gamma_t$ and
$U_{\ell}$ are restricted to this subcategory.

Section~\ref{sec-quasi-reflection} gives the parallel quasitopological
correspondence and shows how the topological reconstruction is obtained
from the quasitopological one through reflection into $\Topulg$.
\end{rem}

We conclude with topological consequences of the representation
theorem and the Hausdorff restriction of the correspondence.
The canonical topological Mundici group gives a retraction onto the
unit interval and transfers separation properties from the MV-algebra
to that group.

\begin{cor}
\label{thm-retract-Tychonoff}
\label{thm-intro-group-envelope}
Every topological MV-algebra $(A,\tau)$ is a retract of its canonical
topological Mundici group $(G_A,\sigma_A)$.
If $(A,\tau)$ is $T_0$, then $(G_A,\sigma_A)$ is Tychonoff,
and consequently so is $(A,\tau)$.
\end{cor}

\begin{proof}
Identify $A=[0,u_A]$.  By
Theorem~\ref{thm-topological-representation}, the inclusion
$j_A\colon(A,\tau)\longhookrightarrow(G_A,\sigma_A)$ is a
topological embedding.  Since $\sigma_A$ is an $\ell$-group topology,
the interval truncation
\[
 \pi_0\colon(G_A,\sigma_A)\longrightarrow(A,\tau),
 \qquad \pi_0(g)=(g\wedge u_A)\vee0,
\]
is continuous and satisfies $\pi_0j_A=1_A$.  Thus $A$ is a retract.

Now assume that $(A,\tau)$ is $T_0$.  The bounded coordinates
separate points of $G_A$: any two elements lie in a common $I_n$,
and $q_ns_n=1_{I_n}$ implies that $s_n$ is injective.
By the coordinate-initial characterization in
Proposition~\ref{prop-kappa-lambda}, the map
\[
 (G_A,\kappa_A)\longrightarrow(A,\tau)^{\mathbb Z},
 \qquad g\mapsto(\pi_k(g))_{k\in\mathbb Z},
\]
is therefore a topological embedding.  A product of $T_0$ spaces
is $T_0$, so $(G_A,\kappa_A)$ is $T_0$.
Since $\kappa_A\subseteq\sigma_A$ by
\eqref{eq-topological-reconstruction-chain}, $(G_A,\sigma_A)$ is
$T_0$ as well.  Every $T_0$ topological group is Tychonoff;
hence $(G_A,\sigma_A)$ and its subspace $(A,\tau)$ are Tychonoff.
\end{proof}

Hoo established the equivalence of $T_0,T_1$, and $T_2$ for
topological MV-algebras \cite{Hoo1997}.
Gan et al.\ extended the equivalence through $T_3$ \cite{Gan2024}.
Since Corollary~\ref{thm-retract-Tychonoff} makes every $T_0$
topological MV-algebra Tychonoff, the standard implications between
separation axioms give the following equivalence.
Here $T_3$ means regular and $T_1$, while $T_{3\frac12}$ means
completely regular and $T_1$, equivalently, Tychonoff.

\begin{cor}
\label{cor-separation-axioms}
For topological MV-algebras the separation axioms
$T_0,T_1,T_2,T_3$, and $T_{3\frac12}$ are equivalent.
\end{cor}

\begin{cor}
\label{cor-Hausdorff-adjunction}
If $(A,\tau)$ is Hausdorff and $\Xi A=(G_A,u_A)$, then
$(G_A,u_A,\sigma_A)$ is Hausdorff.
Consequently the adjunction of Theorem~\ref{thm-topological-adjunction}
and the equivalence of Theorem~\ref{thm-topological-correspondence}
restrict to the corresponding full Hausdorff subcategories.
\end{cor}

\begin{proof}
If $(A,\tau)$ is Hausdorff, it is $T_0$, so
Corollary~\ref{thm-retract-Tychonoff} shows that
$(G_A,\sigma_A)$ is Tychonoff and hence Hausdorff.

Finally, $\Gamma_t$ preserves Hausdorffness, since it equips the unit
interval with the subspace topology, while the first part of the proof
shows that $\Xi_t$ sends every Hausdorff topological MV-algebra to a
Hausdorff topological unital $\ell$-group.  Thus both functors of the
adjunction preserve Hausdorff objects, so the adjunction of
Theorem~\ref{thm-topological-adjunction} restricts to the corresponding
full Hausdorff subcategories.  Restricting further to the
Mundici-fine objects, the restricted unit and counit appearing in
Theorem~\ref{thm-topological-correspondence} remain isomorphisms.
Hence the equivalence in that theorem restricts to the corresponding
full Hausdorff subcategories as well.
\end{proof}

\section{Extremal \texorpdfstring{$\ell$}{l}-group extensions
under additional hypotheses}
\label{sec-restricted-correspondence}

Sections~\ref{sec-group-envelope} and~\ref{sec-correspondence} show that
the Mundici topology $\sigma_A$ is the largest $\ell$-group extension
of every topological MV-algebra $(A,\tau)$ and establish the
topological Mundici correspondence.
The coordinate initial topology $\kappa_A$ and the bounded-stage final
topology $\lambda_A$ of Section~\ref{sec-coordinate-bounded-reconstruction}
both induce $\tau$ on $A$, with $\kappa_A\subseteq\lambda_A$.
Every $\ell$-group extension $\rho$ of $\tau$ satisfies
$\kappa_A\subseteq\rho\subseteq\lambda_A$, but the definitions of
these bounds alone do not ensure joint continuity of the
$\ell$-group operations on the whole group.
The $k_\omega$ and local convexity hypotheses make the upper and lower
bounds into $\ell$-group extensions, respectively.
In the first case, $\lambda_A=\sigma_A$ is the largest extension;
in the second, $\kappa_A$ is the least extension.
The first realization restricts the topological Mundici correspondence,
while the second gives an alternative equivalence through the least
extension.

\subsection{\texorpdfstring{$k_\omega$}{k-omega} topological MV-algebras:
\texorpdfstring{$\lambda_A=\sigma_A$}{lambda\_A=sigma\_A}}
\label{sec-komega}

For a $k_\omega$ topological MV-algebra, the bounded-stage final
topology $\lambda_A$ is itself an $\ell$-group topology and hence
coincides with the Mundici topology $\sigma_A$.

Recall that a sequence $(K_n)_{n\geq1}$ of compact Hausdorff
subspaces of a space $X$ is called a \emph{$k_\omega$-decomposition}
of $X$ if
\[
 K_1\subseteq K_2\subseteq\cdots,\qquad
 X=\bigcup_{n\geq1}K_n,
\]
and a subset $U\subseteq X$ is open precisely when
$U\cap K_n$ is open in $K_n$ for every $n$.
A space admitting a $k_\omega$-decomposition is called a
\emph{$k_\omega$-space} \cite{Franklin1977}.  Recall also that a space $X$ is a
\emph{$k$-space} if a subset $U\subseteq X$ is open precisely when
$U\cap K$ is open in $K$ for every compact subset $K\subseteq X$ \cite{Willard1970}.
We shall use the following standard facts about $k_\omega$-spaces
and $k_\omega$-decompositions
(see, e.g., \cite[items~1, 3, 4, and~14]{Franklin1977}):
\begin{enumerate}[label=(K\arabic*),ref=(K\arabic*)]
\item\label{komega-Hausdorff-k-space}
every $k_\omega$-space is a Hausdorff $k$-space;
\item\label{komega-compact-stage}
if $X=\bigcup_{n\geq1}X_n$ is a $k_\omega$-decomposition, then every
compact subset of $X$ is contained in some $X_n$;
\item\label{komega-finite-products}
if $X=\bigcup_{n\geq1}X_n$ and $Y=\bigcup_{n\geq1}Y_n$ are
$k_\omega$-decompositions, then
\[
 X\times Y=\bigcup_{n\geq1}(X_n\times Y_n)
\]
is a $k_\omega$-decomposition; in particular, by iteration,
$(X_n^m)_n$ is a $k_\omega$-decomposition of $X^m$ for every $m\geq1$;
\item\label{komega-closed-subspaces}
if $X=\bigcup_{n\geq1}X_n$ is a $k_\omega$-decomposition and
$Y\subseteq X$ is closed, then
$Y=\bigcup_{n\geq1}(Y\cap X_n)$ is a $k_\omega$-decomposition.
\end{enumerate}

The following example provides a noncompact $k_\omega$ topological MV-algebra:

\begin{exmp}
Let
\[
 A_{\mathrm{ec}}
 =
 \bigl\{(a_m)_{m\geq1}\in[0,1]^{\mathbb N}:
        (a_m)_{m\geq1}\text{ is eventually constant}\bigr\},
\]
with the coordinatewise MV-operations, and for $n\geq1$ set
\[
 K_n=
 \bigl\{a\in A_{\mathrm{ec}}:a_m=a_n\text{ for every }m\geq n\bigr\}.
\]
With the subspace topology inherited from $[0,1]^{\mathbb N}$,
$K_n$ is homeomorphic to $[0,1]^n$, hence is compact Hausdorff, and
$K_n\subseteq K_{n+1}$ with
$A_{\mathrm{ec}}=\bigcup_{n\geq1}K_n$.  Equip $A_{\mathrm{ec}}$
with the final topology $\tau_{\mathrm{ec}}$ of the family of inclusions
\[
 \bigl(i_n\colon K_n\longhookrightarrow A_{\mathrm{ec}}\bigr)_{n\geq1}.
\]

To check that these are indeed the subspace topologies, define
$\chi_n\colon A_{\mathrm{ec}}\longrightarrow K_n$ by
\[
 (\chi_n(a))_m=
 \begin{cases}
  a_m,&m\leq n,\\
  a_n,&m>n.
 \end{cases}
\]
For every $m$, the restriction $\chi_n|_{K_m}$ is continuous, so $\chi_n$ is
continuous by finality.  The inclusion $i_n$ is continuous by definition, and
$\chi_ni_n=1_{K_n}$; hence $i_n$ is a topological embedding.  Thus
$(K_n)_n$ is a $k_\omega$-decomposition of
$(A_{\mathrm{ec}},\tau_{\mathrm{ec}})$.

The MV-negation restricts to a continuous map $K_n\longrightarrow K_n$ for every
$n$, so it is continuous on $A_{\mathrm{ec}}$ by finality.  Likewise,
$\oplus$ restricts to a continuous map $K_n\times K_n\longrightarrow K_n$.
By \ref{komega-finite-products}, $(K_n\times K_n)_n$ is a
$k_\omega$-decomposition of $A_{\mathrm{ec}}\times A_{\mathrm{ec}}$,
and hence $\oplus$ is jointly continuous.  Therefore
$(A_{\mathrm{ec}},\tau_{\mathrm{ec}})$ is a $k_\omega$ topological
MV-algebra.  

Moreover, $A_{\mathrm{ec}}$ is nowhere locally compact.
Indeed, each $K_n$ has empty interior in $K_{n+1}$,
where it is defined by $a_{n+1}=a_n$.
Consequently, $K_n$ has empty interior in $A_{\mathrm{ec}}$.
By \ref{komega-compact-stage}, every compact subset of
$A_{\mathrm{ec}}$ is contained in some $K_n$.
Thus no point has a compact neighbourhood.
In particular, $A_{\mathrm{ec}}$ is not compact.
\end{exmp}

\begin{thm}
\label{thm-komega-reconstruction}
Let $(A,\tau)$ be a $k_\omega$ topological MV-algebra and write
$\Xi A=(G,u)$.  Then:
\begin{enumerate}[label=(\alph*),ref=(\alph*)]
\item\label{thm-komega-reconstruction-topology}
$(G,\lambda_A)$ is a $k_\omega$ topological unital abelian
$\ell$-group;
\item\label{thm-komega-reconstruction-compact}
every compact subset of $(G,\lambda_A)$ is order bounded;
\item\label{thm-komega-reconstruction-maximality}
the bounded-stage final topology $\lambda_A$ coincides with the
Mundici topology $\sigma_A$ and is therefore the largest
$\ell$-group extension of $\tau$.
\end{enumerate}
\end{thm}

\begin{proof}
\ref{thm-komega-reconstruction-topology}
Choose a $k_\omega$-decomposition
\[
 K_1\subseteq K_2\subseteq\cdots,\qquad
 A=\bigcup_{j\geq1}K_j.
\]
Since $0,u\in A$ and the $K_j$ are increasing, there is $N$ such that
$0,u\in K_N$.  Passing to the tail $(K_j)_{j\geq N}$ and reindexing
preserves the $k_\omega$-decomposition, so we may assume that
$0,u\in K_j$ for every $j$.

By \ref{komega-Hausdorff-k-space}, $A$ is Hausdorff.
Since $s_j$ embeds $I_j$ into $A^{2j}$ by
Lemma~\ref{lem-stage-topology}\ref{lem-stage-topology-quotient},
each $I_j$ is Hausdorff as well.
Set
\[
 C_j=q_j(K_j^{2j})\subseteq I_j\qquad(j\geq1).
\]
Each $C_j$ is compact Hausdorff by continuity of $q_j$.
By Proposition~\ref{prop-lambda-exact-space}, $\lambda_A|_{I_j}=\tau_j$,
so each $C_j$ carries its $\lambda_A$-subspace topology.

For $j\geq n$, padding a tuple in $K_j^{2n}$ with $j-n$ copies
of $u$ and $j-n$ copies of $0$ leaves its value under the
corresponding $q$ map unchanged.  Hence
\[
 q_n(K_j^{2n})\subseteq C_j\qquad(j\geq n).
\]
Together with $K_j\subseteq K_{j+1}$, this gives
$C_j\subseteq C_{j+1}$.  These sets cover $G$: write any $g\in G$
as $q_n(a)$ with $a\in A^{2n}$, and choose $j\geq n$ large enough
that all coordinates of $a$ lie in $K_j$.  Then $g\in C_j$.

We claim that $(C_j)_j$ determines $\lambda_A$.
If $U\in\lambda_A$, then each $U\cap C_j$ is open in $C_j$.
Conversely, suppose that these intersections are open, and fix $n\geq1$.
For $j\geq n$, the restriction
\[
 q_n|_{K_j^{2n}}\colon K_j^{2n}\longrightarrow C_j
\]
is continuous: it factors through the continuous padding map
$K_j^{2n}\longrightarrow K_j^{2j}$ and
$q_j\colon K_j^{2j}\longrightarrow C_j$.  Thus
\[
 q_n^{-1}(U\cap I_n)\cap K_j^{2n}
 =\bigl(q_n|_{K_j^{2n}}\bigr)^{-1}(U\cap C_j)
\]
is open in $K_j^{2n}$ for every $j\geq n$.
By \ref{komega-finite-products}, the tail $(K_j^{2n})_{j\geq n}$
is a $k_\omega$-decomposition of $A^{2n}$, so
$q_n^{-1}(U\cap I_n)$ is open in $A^{2n}$.
As this holds for every $n$, \eqref{eq-lambda-quotient-test}
gives $U\in\lambda_A$.  This proves the claim, and hence
$(C_j)_j$ is a $k_\omega$-decomposition of $(G,\lambda_A)$.

It remains to prove continuity of the $\ell$-group operations.
By \ref{komega-finite-products},
\[
 G\times G=\bigcup_{j\geq1}(C_j\times C_j)
\]
is a $k_\omega$-decomposition.  Since $C_j\subseteq I_j$,
Lemma~\ref{lem-stage-topology}%
\ref{lem-stage-topology-operations} gives continuous restrictions
\[
 +\colon C_j\times C_j\longrightarrow I_{2j},
 \qquad
 \vee,\wedge\colon C_j\times C_j\longrightarrow I_j.
\]
Composing with the corresponding stage inclusions into
$(G,\lambda_A)$ shows that the restrictions of $+$, $\vee$, and
$\wedge$ to every $C_j\times C_j$ are continuous.  Since
$(C_j\times C_j)_j$ determines the topology of $G\times G$, these
three operations are continuous on $G\times G$.

Inversion is already continuous by
Proposition~\ref{prop-lambda-quasitopological}.
Thus $(G,\lambda_A)$ is a $k_\omega$ topological unital abelian
$\ell$-group, proving \ref{thm-komega-reconstruction-topology}.

\ref{thm-komega-reconstruction-compact}
Let $K\subseteq G$ be compact.  By \ref{komega-compact-stage}, $K$ is
contained in some $C_j$.  Since
\[
 C_j\subseteq I_j=[-ju,ju],
\]
the set $K$ is order bounded.

\ref{thm-komega-reconstruction-maximality}
By \ref{thm-komega-reconstruction-topology}, $\lambda_A$ is itself an
$\ell$-group topology.  Hence
Theorem~\ref{thm-lambda-envelope}%
\ref{thm-lambda-envelope-maximality} gives
$\sigma_A=\lambda_A$.  Since $\sigma_A$ is the largest
$\ell$-group extension of $\tau$ by
Corollary~\ref{cor-reconstruction-topology-chain}, the same is true of
$\lambda_A$.
\end{proof}

We next use the preceding reconstruction to characterize
Mundici-fine objects intrinsically on the $\ell$-group side
when the unit interval is a $k_\omega$-space.
By \eqref{eq-counit-Mundici-maximality}, this amounts to
characterizing when the counit is a homeomorphism.

\begin{thm}
\label{thm-intrinsic-komega}
Let $(G,u,\rho)$ be a topological unital abelian $\ell$-group,
and write
\[
 (A,\tau)=([0,u],\rho|_{[0,u]}).
\]
Assume that $A$ is a $k_\omega$-space.  The following are equivalent:
\begin{enumerate}[label=(\roman*),ref=(\roman*)]
\item\label{thm-intrinsic-komega-counit}
$(G,u,\rho)$ is Mundici-fine.  Equivalently, the counit
\[
 \varepsilon_{(G,u,\rho)}\colon
 \Xi_t\Gamma_t(G,u,\rho)\longrightarrow(G,u,\rho)
\]
is a homeomorphism.
\item\label{thm-intrinsic-komega-final}
The topology $\rho$ is the final topology of the family
of inclusions
\[
 \bigl(j_n\colon(I_n,\rho|_{I_n})
 \longhookrightarrow G\bigr)_{n\geq1},
 \qquad I_n=[-nu,nu].
\]
Equivalently, a subset $U\subseteq G$ is $\rho$-open
if and only if $U\cap I_n$ is open in $I_n$ for every $n$.
\item\label{thm-intrinsic-komega-k-space}
$G$ is a $k$-space and every compact subset of $G$ is order bounded.
\item\label{thm-intrinsic-komega-komega}
$G$ is a $k_\omega$-space and every compact subset of $G$ is order
bounded.
\end{enumerate}
\end{thm}

\begin{proof}
Write $\Xi A=(G_A,u_A)$, and abbreviate the algebraic
Mundici counit \eqref{eq-transported-counit} by
\[
 \varepsilon=\varepsilon_{(G,u)}\colon
 (G_A,u_A)\longrightarrow(G,u).
\]
For every $n\geq1$, put $I_n^A=[-nu_A,nu_A]$.
Then the restriction 
\[
 \varepsilon|_{I_n^A}\colon I_n^A\longrightarrow I_n
\]
is an order isomorphism.

Recall from \eqref{eq-transported-maximal-topology}
that $\widetilde\sigma_G$ is the topology obtained by
transporting $\sigma_A$ along $\varepsilon$.
By Proposition~\ref{prop-transported-maximality},
$\widetilde\sigma_G$ and $\rho$ induce the same topology
$\tau$ on $[0,u]$.  Hence
Proposition~\ref{prop-lgroup-extension-bounded-stages} gives
\begin{equation}
\label{eq-intrinsic-komega-stagewise}
 \widetilde\sigma_G|_{I_n}=\rho|_{I_n}
 \qquad(n\geq1).
\end{equation}

Since $A$ is a $k_\omega$-space,
Theorem~\ref{thm-komega-reconstruction}%
\ref{thm-komega-reconstruction-maximality}
gives $\sigma_A=\lambda_A$.
Thus, for every $U\subseteq G$,
\begin{align*}
 U\in\widetilde\sigma_G
 &\iff \varepsilon^{-1}(U)\in\sigma_A=\lambda_A
 &&\text{(by \eqref{eq-transported-maximal-topology})}\\
 &\iff \varepsilon^{-1}(U)\cap I_n^A
 \in\lambda_A|_{I_n^A}\quad(n\geq1)
 &&\text{(by \eqref{eq-lambda-stage-test}
 and \eqref{eq-lambda-stage-restriction})}\\
 &\iff U\cap I_n\in\widetilde\sigma_G|_{I_n}=\rho|_{I_n}
 \quad(n\geq1)
 &&\text{(by \eqref{eq-transported-maximal-topology}
 and \eqref{eq-intrinsic-komega-stagewise})}.
\end{align*}
By \eqref{eq-counit-Mundici-maximality},
\ref{thm-intrinsic-komega-counit} holds precisely when
$\rho=\widetilde\sigma_G$, which by the displayed
equivalence is exactly \ref{thm-intrinsic-komega-final}.

If \ref{thm-intrinsic-komega-final} holds, the preceding
argument shows that
\[
 \varepsilon\colon
 (G_A,\lambda_A)\longrightarrow(G,\rho)
\]
is a homeomorphism.  Thus
Theorem~\ref{thm-komega-reconstruction}%
\ref{thm-komega-reconstruction-topology}%
\ref{thm-komega-reconstruction-compact}
implies that $(G,\rho)$ is a $k_\omega$-space and that
every compact subset of $G$ is order bounded.  
Thus \ref{thm-intrinsic-komega-final} implies
\ref{thm-intrinsic-komega-komega}, and
\ref{komega-Hausdorff-k-space} shows that the latter implies
\ref{thm-intrinsic-komega-k-space}.

Finally assume \ref{thm-intrinsic-komega-k-space}, and let $U\subseteq G$ have the property that
$U\cap I_n$ is open in $I_n$ for every $n$.  If $K\subseteq G$ is compact,
then $K$ is order bounded.  Since $u$ is a strong unit, there is $n$ with
$K\subseteq I_n$.  Therefore $U\cap K$ is open in $K$.  Since $G$ is a
$k$-space, $U$ is open.  Hence $\rho$ is the final topology of the family
$(j_n)_{n\geq1}$ in \ref{thm-intrinsic-komega-final}.
\end{proof}

Write $\TopMVkw$ for the full subcategory of $k_\omega$ objects of
$\TopMV$, and $\Topulgkw$ for the full subcategory of $k_\omega$
objects of $\Topulg$ whose compact subsets are all order bounded.
Here $\mathrm{cob}$ stands for compact-order-boundedness.

\begin{thm}
\label{thm-komega-equivalence}
The topological Mundici correspondence of
Theorem~\ref{thm-topological-correspondence} restricts to an equivalence
of full subcategories
\[
 \TopMVkw\simeq\Topulgkw.
\]
The restrictions of $\Xi_t$ and $\Gamma_t$ are mutually quasi-inverse.
\end{thm}

\begin{proof}
By Theorem~\ref{thm-komega-reconstruction}, $\Xi_t$ sends
$\TopMVkw$ to $\Topulgkw$.  Its values are Mundici-fine by
Theorem~\ref{thm-topological-correspondence}.
Conversely, let $(G,u,\rho)\in\Topulgkw$.  By
\ref{komega-Hausdorff-k-space}, $G$ is Hausdorff.  
We claim that the interval $[0,u]$ is closed in $G$. 
Indeed, the positive cone $G_+$ is closed in a Hausdorff
topological $\ell$-group, since
$G_+=\{x:(-x)^+=0\}$ and the map $x\mapsto(-x)^+$
is continuous.  Moreover, $u-G_+$ is also closed,
as $x\mapsto u-x$ is a homeomorphism.  Hence
\[
 [0,u]=G_+\cap(u-G_+)
\]
is closed in $G$.
By \ref{komega-closed-subspaces}, the closed subspace $[0,u]$ is
$k_\omega$; hence $\Gamma_t(G,u,\rho)\in\TopMVkw$.
Again by \ref{komega-Hausdorff-k-space}, $G$ is a $k$-space, and by
assumption every compact subset is order bounded.  Hence
Theorem~\ref{thm-intrinsic-komega}%
\ref{thm-intrinsic-komega-k-space} implies
Theorem~\ref{thm-intrinsic-komega}%
\ref{thm-intrinsic-komega-counit}: the counit is a homeomorphism.
Thus $(G,u,\rho)$ is Mundici-fine, and
$\Topulgkw$ is a full subcategory of $\TopulgM$.
Both functors $\Xi_t$ and $\Gamma_t$ therefore restrict to the stated
subcategories.  The unit and counit of
Theorem~\ref{thm-topological-correspondence} restrict to natural
isomorphisms between them.
\end{proof}

\begin{cor}
\label{cor-compact-reconstruction}
If $A$ is a compact Hausdorff topological MV-algebra and
$\Xi A=(G,u)$, then the bounded stages
$I_n=[-nu,nu]$, equipped with their canonical
compact Hausdorff topologies, form a
$k_\omega$-decomposition of $\Xi_tA$.

On the group side, if $(G,u,\rho)\in\Topulg$ and $[0,u]$ is compact
Hausdorff, then the counit is a homeomorphism if and only if
$G$ is a $k$-space and every compact subset of $G$ is order bounded;
equivalently, $(I_n)_{n\geq1}$ is a $k_\omega$-decomposition of $(G,\rho)$,
where $I_n=[-nu,nu]$.
\end{cor}

\begin{proof}
Since $A$ is compact Hausdorff, it is a $k_\omega$-space.
Thus Theorem~\ref{thm-komega-reconstruction}%
\ref{thm-komega-reconstruction-maximality}
gives $\sigma_A=\lambda_A$.  For every $n\geq1$,
Lemma~\ref{lem-stage-topology}\ref{lem-stage-topology-quotient}
shows that $q_n\colon A^{2n}\longrightarrow(I_n,\tau_n)$ is a
continuous surjection and that $s_n$ embeds $(I_n,\tau_n)$ into
$A^{2n}$.  Since $A^{2n}$ is compact Hausdorff, so is
$(I_n,\tau_n)$.  By \eqref{eq-lambda-stage-restriction}, $\tau_n$
is the subspace topology induced by $\lambda_A$,
and hence by $\sigma_A$.
Moreover, \eqref{eq-lambda-stage-test} says that $\lambda_A$ is
final with respect to the stage inclusions.  Together with
\eqref{eq-bounded-stage-exhaustion}, these facts show that
$(I_n)_{n\geq1}$ is a $k_\omega$-decomposition of $\Xi_tA$.

For the group-side assertion, put
$(A,\tau)=([0,u],\rho|_{[0,u]})$.
By Proposition~\ref{prop-lgroup-extension-bounded-stages},
the map $s_n$ embeds $(I_n,\rho|_{I_n})$ into
$(A^{2n},\tau^{2n})$, so $I_n$ is Hausdorff.
Moreover, $q_n\colon A^{2n}\longrightarrow I_n$
is a continuous surjection, so $I_n$ is compact.
Theorem~\ref{thm-intrinsic-komega} identifies the
counit-homeomorphism condition with the stated $k$-space
condition and with finality for the inclusions
$I_n\longhookrightarrow G$.
Since the $I_n$ are increasing compact Hausdorff
subspaces whose union is $G$ by
\eqref{eq-bounded-stage-exhaustion},
the latter final-topology condition is exactly the
$k_\omega$-decomposition condition.
\end{proof}

\subsection{Locally convex topological MV-algebras:
the least extension \texorpdfstring{$\kappa_A$}{kappa\_A}}
\label{sec-locally-convex}

We turn to the other canonical topology from
Section~\ref{sec-coordinate-bounded-reconstruction}.  Local convexity
makes $\kappa_A$ an unbounded locally solid $\ell$-group topology.
Its initial definition then identifies it as the least $\ell$-group
extension of the interval topology, complementing the maximal
reconstruction by $\lambda_A$ in the $k_\omega$ case.

A subset $C$ of an MV-algebra $A$ is \emph{order convex} if
$x,z\in C$ and $x\leq y\leq z$ imply $y\in C$.  A topological
MV-algebra is \emph{locally convex} if every point has a neighbourhood
base of order-convex sets.  Since $0$ is the least element, an
order-convex set containing $0$ is precisely a downward-closed set.
By \cite[Proposition~2.3]{Weber2012}, we obtain the following
characterization, used repeatedly below:
\begin{enumerate}[label=(LC),ref=(LC)]
\item\label{lc-zero-neighbourhood-base}
A topological MV-algebra $(A,\tau)$ is locally convex if and only if
$0$ has a neighbourhood base consisting of downward-closed sets.
\end{enumerate}

A subset $V$ of an $\ell$-group $G$ is \emph{solid} if
$|x|\leq|y|$ and $y\in V$ imply $x\in V$.  A group topology $\tau$ on
$G$ is \emph{locally solid} if $0$ has a neighbourhood base consisting of
solid sets.  For background on locally solid topological $\ell$-groups,
see \cite{Khan2007,Hong2015}.

Let $(G,\tau)$ be a locally solid abelian $\ell$-group.  Following
\cite{Zabeti2020}, the \emph{unbounded topology} $\mathfrak u\tau$ has a
zero-neighbourhood base
\begin{equation}
\label{eq-unbounded-topology-base}
 V_a=\{x\in G:|x|\wedge a\in V\},
 \qquad a\in G_+,
\end{equation}
where $V$ ranges through a solid zero-neighbourhood base for $\tau$.
It is again a locally solid $\ell$-group topology, and Hausdorffness is
preserved \cite[Theorem~2.11]{Zabeti2020}.\footnote{The paper
\cite{Zabeti2020} assumes throughout that all singletons are bounded;
see Remark~2.1 therein.  This assumption is not needed for the
neighbourhood construction in the proof of \cite[Theorem~2.11]{Zabeti2020}, 
which remains valid here and gives the stated conclusion.}
We call $\tau$ \emph{unbounded}
if $\tau=\mathfrak u\tau$; see also \cite{Abela2023}.

When a strong unit $u$ is available, the neighbourhood base in
\eqref{eq-unbounded-topology-base} can be described using only
truncation by $u$.
To prove this reduction, we first establish the following elementary
lattice estimate in an abelian $\ell$-group:
\begin{equation}
\label{eq-lattice-meet-sum}
 b\wedge(c+d)\leq(b\wedge c)+(b\wedge d)
 \qquad(b,c,d\geq0).
\end{equation}
Indeed, put $a=b\wedge(c+d)$, $a_1=a\wedge c$, and
$a_2=a-a_1=(a-c)^+$.  Then $a=a_1+a_2$ and
\[
 0\leq a_1\leq b\wedge c,
 \qquad
 0\leq a_2\leq b\wedge d.
\]
For the second inequality, $a\leq c+d$ gives $(a-c)^+\leq d$,
while $a_1\geq0$ gives $a_2\leq a\leq b$.
Adding the two inequalities proves \eqref{eq-lattice-meet-sum}.

\begin{lem}
\label{lem-strong-unit-unbounded}
Let $(G,u)$ be a unital abelian $\ell$-group with strong unit and let
$\tau$ be a locally solid topology.  Then, as $V$ ranges through a solid
zero-neighbourhood base for $\tau$, the sets
\begin{equation}
\label{eq-strong-unit-unbounded-base}
 V_u=\{x\in G:|x|\wedge u\in V\},
\end{equation}
form a zero-neighbourhood base for $\mathfrak u\tau$.
\end{lem}

\begin{proof}
By \eqref{eq-unbounded-topology-base}, each $V_u$ is already
a basic zero-neighbourhood of $\mathfrak u\tau$.
Thus it suffices to show that every basic zero-neighbourhood
$V_a$ contains some $W_u$, where $W$ is a solid
zero-neighbourhood.

Fix $a\in G_+$ and a solid zero-neighbourhood $V$.
Since $u$ is a strong unit, choose $n\geq1$ with $a\leq nu$.
By continuity of addition and local solidity, choose a solid
zero-neighbourhood $W$ with $nW\subseteq V$.

If $x\in W_u$, then $|x|\wedge u\in W$.
By iterating \eqref{eq-lattice-meet-sum}, we obtain
\[
 |x|\wedge a
 \leq |x|\wedge(nu)
 \leq n(|x|\wedge u)
 \in nW\subseteq V.
\]
Since $V$ is solid,
we obtain $|x|\wedge a\in V$, hence $x\in V_a$.
Thus $W_u\subseteq V_a$, as required.
\end{proof}

\begin{prop}
\label{prop-interval-locally-convex}
Let $(G,u,\tau)$ be a locally solid topological unital abelian
$\ell$-group.  Then the subspace topology on $[0,u]$ is locally convex as
a topological MV-algebra.
\end{prop}

\begin{proof}
Let $\mathcal V$ be a solid zero-neighbourhood base for $\tau$.
Then
\[
 \{V\cap[0,u]:V\in\mathcal V\}
\]
is a zero-neighbourhood base for the subspace topology on $[0,u]$.
Each $V\cap[0,u]$ is downward closed: if
$0\leq x\leq y\in V\cap[0,u]$, then
$|x|\leq|y|$, so solidity gives $x\in V\cap[0,u]$.
Local convexity now follows from \ref{lc-zero-neighbourhood-base}.
\end{proof}

For a locally convex $(A,\tau)\in\TopMV$, write
$\Xi A=(G_A,u_A)$ and identify $A=[0,u_A]\subseteq G_A$.
For every downward-closed zero-neighbourhood $U$ of $A$, put
\begin{equation}
\label{eq-unbounded-neighbourhood}
 W_U=\{g\in G_A:|g|\wedge u_A\in U\}.
\end{equation}
Weber already observed that these sets define a locally solid
extension of $\tau$ \cite[p.~69]{Weber2023}.
The next theorem identifies this extension with the coordinate
topology $\kappa_A$ and proves its unboundedness; its extremal and
categorical properties follow below.  Without local convexity,
$\kappa_A$ need not be a group topology, as
Example~\ref{exmp-kappa-not-group} shows.

\begin{thm}[Locally convex coordinate reconstruction]
\label{lem-unbounded-lift}
\label{thm-locally-convex-reconstruction}
Let $(A,\tau)$ be a locally convex topological MV-algebra and write
$\Xi A=(G_A,u_A)$.  The sets $W_U$ in
\eqref{eq-unbounded-neighbourhood}, as $U$ ranges through the
downward-closed zero-neighbourhoods of $A$, form a zero-neighbourhood
base for $\kappa_A$.  Moreover, $\kappa_A$ is an unbounded locally solid
$\ell$-group topology and
\[
 \kappa_A|_{[0,u_A]}=\tau.
\]
\end{thm}

\begin{proof}
Write $(G,u)=(G_A,u_A)$ and identify $A=[0,u]$.

\textbf{Step 1.} We show that the sets $W_U$ form a zero-neighbourhood
base for $\kappa_A$.

Put $F(g)=|g|\wedge u$.  By \eqref{eq-pi-k} and \eqref{eq-Gamma-def}, we have
\[
 \pi_0(g)=g^+\wedge u,
 \qquad \neg\pi_{-1}(g)=(-g)^+\wedge u.
\]
Since $|g|=g^+\vee(-g)^+$, distributivity yields
\begin{equation}
\label{eq-modulus-coordinate-truncation}
 F(g)=\pi_0(g)\vee\neg\pi_{-1}(g).
\end{equation}
By Proposition~\ref{prop-kappa-lambda}, both coordinate maps
are continuous from $(G,\kappa_A)$ to $(A,\tau)$;
since $(A,\tau)$ is a topological MV-algebra,
MV-negation and join are also continuous.
Thus $F$ is continuous. 
Since $F(0)=0$, every
$W_U=F^{-1}(U)$ is a $\kappa_A$-neighbourhood of zero, even when
$U$ is not open.

To show that the sets $W_U$ form a zero-neighbourhood
base, let $N$ be a $\kappa_A$-neighbourhood of zero.
The coordinate initial description of $\kappa_A$ in
Proposition~\ref{prop-kappa-lambda} gives a finite set
$E\subseteq\mathbb Z$ and open neighbourhoods $O_k$ of $\pi_k(0)$
in $A$ such that
\begin{equation}
\label{eq-coordinate-neighbourhood-control}
 \bigcap_{k\in E}\pi_k^{-1}(O_k)\subseteq N.
\end{equation}
By \eqref{eq-pi-k}, the coordinates are decreasing in $k$;
that is, $\pi_{k+1}(g)\leq\pi_k(g)$ for all
$g\in G$ and $k\in\mathbb Z$.
Together with \eqref{eq-modulus-coordinate-truncation},
this yields
\begin{equation}
\label{eq-coordinate-modulus-bounds}
 \begin{cases}
  \pi_k(g)\leq\pi_0(g)\leq F(g),&k\geq0,\\
  \neg\pi_k(g)\leq\neg\pi_{-1}(g)\leq F(g),&k<0.
 \end{cases}
\end{equation}
By \eqref{eq-zero-coordinate-values}, $O_k$ is an
open zero-neighbourhood for $k\geq0$.
For $k<0$, the set
$\neg O_k=\{\neg a:a\in O_k\}$ is also an
open zero-neighbourhood, since MV-negation is a
homeomorphism sending $u$ to $0$.
By \ref{lc-zero-neighbourhood-base}, there is a downward-closed
zero-neighbourhood $U$ contained in
\[
 \bigcap_{\substack{k\in E\\k\geq0}}O_k
 \ \cap\!
 \bigcap_{\substack{k\in E\\k<0}}\neg O_k.
\]
If $g\in W_U$, downward closure and
\eqref{eq-coordinate-modulus-bounds} give
$\pi_k(g)\in O_k$ for every $k\in E$.
Hence $W_U\subseteq N$ by
\eqref{eq-coordinate-neighbourhood-control}, as required.

The topology $\kappa_A$ is translation invariant: for every
$h\in G$, the identity
$t_a(x+h)=t_{a+h}(x)$ and the defining initial topology
make translation by $h$ continuous.  Translation by $-h$
is its continuous inverse, so all translations are
$\kappa_A$-homeomorphisms.  Consequently, the sets
$h+W_U$, as $U$ varies, form a neighbourhood base at $h$.

\textbf{Step 2.} We prove that $\kappa_A$ is a locally solid
$\ell$-group topology.
Given $U$, continuity of $\oplus$ at $(0,0)$, 
together with \ref{lc-zero-neighbourhood-base}, provides a
downward-closed zero-neighbourhood $V$ with $V\oplus V\subseteq U$.  
For $g,h\in W_V$, put
$a=|g|\wedge u$ and $b=|h|\wedge u$.  The inequality
$|g+h|\leq|g|+|h|$, \eqref{eq-lattice-meet-sum}, and \eqref{eq-Gamma-def} give
\[
 |g+h|\wedge u
 \leq (|g|+|h|)\wedge u
 \leq \bigl((|g|\wedge u)+(|h|\wedge u)\bigr)\wedge u
 =a\oplus b\in U.
\]
Since $U$ is downward closed, $W_V+W_V\subseteq W_U$.
Thus addition is jointly continuous at $(0,0)$, and the translation
homeomorphisms give joint continuity everywhere.  
Since every $W_U$ is symmetric, inversion is continuous
at zero and hence everywhere.

The sets $W_U$ are solid: if $|x|\leq|y|$ and $y\in W_U$, then
$|x|\wedge u\leq|y|\wedge u\in U$, so $x\in W_U$.
Thus the group topology is locally solid.
By \eqref{eq-positive-part-contraction}, $|x^+-y^+|\leq|x-y|$.
Since every $W_U$ is solid, this inequality implies that
$x-y\in W_U$ entails $x^+-y^+\in W_U$.
Hence the positive-part map is continuous. The formulas in
\eqref{eq-lattice-operations-positive-part} now give joint continuity
of $\vee$ and $\wedge$.  Therefore $\kappa_A$ is an $\ell$-group topology.

\textbf{Step 3.} We prove that $\kappa_A$ is unbounded.
Since $F|_A=1_A$, \eqref{eq-unbounded-neighbourhood} gives
$W_U\cap A=U$.  By \eqref{eq-strong-unit-unbounded-base}
and $FF=F$, we also have
\[
 (W_U)_u=\{g:F(g)\in W_U\}=F^{-1}(W_U)=(FF)^{-1}(U)=F^{-1}(U)=W_U.
\]
Each $W_U$ is solid, so \eqref{eq-unbounded-topology-base} makes it
also a zero-neighbourhood for $\mathfrak u\kappa_A$.
Consequently $\kappa_A\subseteq\mathfrak u\kappa_A$.
Conversely, for any locally solid topology $\rho$, every basic
$\mathfrak u\rho$-neighbourhood $V_a$ is a $\rho$-neighbourhood,
since $x\mapsto |x|\wedge a$ is $\rho$-continuous.
Thus $\mathfrak u\kappa_A\subseteq\kappa_A$, proving unboundedness.
The equality $\kappa_A|_A=\tau$ follows from
Proposition~\ref{prop-kappa-lambda}.
\end{proof}

\begin{cor}
\label{cor-extremal-lgroup-extensions}
For each locally convex topological MV-algebra $(A,\tau)$,
$\kappa_A$ is the least $\ell$-group extension of $\tau$.
\end{cor}

\begin{proof}
Theorem~\ref{lem-unbounded-lift} makes $\kappa_A$ an $\ell$-group
extension of $\tau$.  Its minimality follows from
Proposition~\ref{prop-kappa-lambda}.
\end{proof}

Write $\LCMV$ for the full subcategory of locally convex objects of
$\TopMV$.  Let $\ULS$ be the category of unital abelian $\ell$-groups
with an unbounded locally solid topology and continuous unit-preserving
$\ell$-homomorphisms.  By \eqref{eq-Gammat-definition} and
Proposition~\ref{prop-interval-locally-convex}, the interval construction
restricts to a functor
\[
 \Gamma_u\colon\ULS\longrightarrow\LCMV.
\]

\begin{prop}
\label{prop-Xiu-functor}
The construction defines a functor
\[
 \Xi_u\colon\LCMV\longrightarrow\ULS,
 \qquad
 \Xi_u(A,\tau)=(G_A,u_A,\kappa_A).
\]
\end{prop}

\begin{proof}
Theorem~\ref{thm-locally-convex-reconstruction} gives the object
assignment in $\ULS$.
Let $f\colon(A,\tau_A)\longrightarrow(B,\tau_B)$ be a continuous
MV-homomorphism and put $h=\Xi f$.
Since $h$ preserves the strong unit and the group and lattice
operations, its restriction to the unit interval gives
\[
 \pi_k^B h=f\pi_k^A\qquad(k\in\mathbb Z),
\]
where the superscripts indicate the respective Mundici groups.
Every right-hand side is continuous, so the coordinate-initial
characterization in Proposition~\ref{prop-kappa-lambda} makes
$h\colon(G_A,\kappa_A)\longrightarrow(G_B,\kappa_B)$ continuous.
Identities and compositions are inherited from $\Xi$.
\end{proof}

\begin{lem}
\label{lem-recover-unbounded-topology}
Let $(G,u,\tau)$ be a locally solid topological unital abelian
$\ell$-group, and give $A=[0,u]$ its subspace topology.
Under the algebraic Mundici counit
\eqref{eq-transported-counit}, the topology $\kappa_A$
on $\Xi A$ transports to $\mathfrak u\tau$ on $G$, 
as in \eqref{eq-transported-maximal-topology}.
In particular, the counit is a homeomorphism onto
$(G,u,\tau)$ whenever $\tau$ is unbounded.
\end{lem}

\begin{proof}
Let $\varepsilon_{(G,u)}\colon G_A\longrightarrow G$
be the algebraic Mundici counit.  We identify
$(G_A,u_A)$ with $(G,u)$ via this isomorphism,
and hence regard the transported $\kappa_A$ as
a topology on $G$.

By Proposition~\ref{prop-interval-locally-convex},
$A$ is locally convex.
Choose a solid zero-neighbourhood base $\mathcal V$
for $\tau$.
The sets $V\cap A$, for $V\in\mathcal V$, form a base of
downward-closed zero-neighbourhoods in $A$: they are
downward closed by solidity and form a zero-neighbourhood
base by the subspace topology.
Theorem~\ref{thm-locally-convex-reconstruction}
therefore gives the following zero-neighbourhood base
for the transported $\kappa_A$:
\[
 W_{V\cap A}
 =\{g\in G:|g|\wedge u\in V\cap A\}
 =\{g\in G:|g|\wedge u\in V\}
 =V_u,\qquad V\in\mathcal V.
\]
By Lemma~\ref{lem-strong-unit-unbounded}, these are also
a zero-neighbourhood base for $\mathfrak u\tau$.
Thus the transported topology is $\mathfrak u\tau$.
\end{proof}

\begin{thm}
\label{thm-locally-convex-equivalence}
The functors $\Xi_u$ and $\Gamma_u$ give an equivalence
\[
 \LCMV\simeq\ULS,
\]
realizing locally convex topological MV-algebras by their least
$\ell$-group extensions.  The classical Mundici unit and counit are
the natural isomorphisms.
\end{thm}

\begin{proof}
To establish the equivalence, it suffices to show that
the classical Mundici unit and counit are homeomorphisms
for the reconstructed topologies.

For $(A,\tau)\in\LCMV$,
Theorem~\ref{thm-locally-convex-reconstruction}
gives $\kappa_A|_{[0,u_A]}=\tau$ under the
identification induced by the classical unit.
Hence the latter is a homeomorphism
\[
 \eta_A\colon(A,\tau)\longrightarrow
 \Gamma_u\Xi_u(A,\tau).
\]
Conversely, let $(G,u,\tau)\in\ULS$.
By Lemma~\ref{lem-recover-unbounded-topology},
the algebraic Mundici counit
\[
 \varepsilon_{(G,u)}\colon
 \Xi_u\Gamma_u(G,u,\tau)\longrightarrow
 (G,u,\mathfrak u\tau)
\]
is a homeomorphism.
Since $\tau=\mathfrak u\tau$ by the definition
of $\ULS$, it is an isomorphism onto $(G,u,\tau)$.

These homeomorphisms are natural, since
$\Xi_u$ and $\Gamma_u$ act on underlying
homomorphisms as the classical Mundici functors.
They therefore give the asserted equivalence.
\end{proof}

This is an alternative realization compatible with the topological
Mundici correspondence of Theorem~\ref{thm-topological-correspondence}.
Although $\Gamma_u$ is a restriction of $\Gamma_t$, the functor
$\Xi_u$ is generally not a restriction of $\Xi_t$:
the former selects the coordinate topology $\kappa_A$,
whereas the latter selects the Mundici topology $\sigma_A$
(see Example~\ref{exmp-product-cubes}).

\begin{cor}
\label{cor-Hausdorff-locally-convex-equivalence}
The equivalence of categories in
Theorem~\ref{thm-locally-convex-equivalence} restricts
to Hausdorff objects.
\end{cor}

\begin{proof}
Only preservation from the MV side needs comment.
If $(A,\tau)$ is Hausdorff, the coordinate embedding established
in the proof of Corollary~\ref{thm-retract-Tychonoff} realizes
$(G_A,\kappa_A)$ as a subspace of the Hausdorff product
$(A,\tau)^{\mathbb Z}$.  Hence $(G_A,\kappa_A)$ is Hausdorff,
so $\Xi_u$ preserves Hausdorffness.
\end{proof}

The following example gives a compact Hausdorff locally convex
topological MV-algebra $(A,\tau)$ for which
$(G_A,u_A,\kappa_A)$ and $(G_A,u_A,\sigma_A)$ are not isomorphic
as topological unital abelian $\ell$-groups.

\begin{exmp}
\label{exmp-product-cubes}
We return to the topological MV-algebra of
Example~\ref{exmp-kappa-lambda-strict}:
\[
 A=[0,1]^{\mathbb N},\qquad
 G=\ell^\infty(\mathbb N),\qquad u=\mathbf1.
\]
The pointwise topology $\rho_{\mathrm{pt}}=\kappa_A$ is an
$\ell$-group extension of the product topology $\tau$ on $A$.
Since $A$ is compact Hausdorff, it is $k_\omega$ (take the constant
decomposition $K_n=A$).  The product topology has a base of
order-convex sets obtained by restricting finitely many coordinates
to relatively open intervals in $[0,1]$, so $A$ is also locally
convex.  It therefore belongs to both classes considered in this section.
Corollary~\ref{cor-extremal-lgroup-extensions} identifies $\kappa_A$ as
the least $\ell$-group extension of $\tau$, while
Corollary~\ref{cor-compact-reconstruction} identifies the largest as
$\lambda_A=\sigma_A$, the final topology on $G$ induced by the
inclusions of the compact cubes $[-n,n]^{\mathbb N}$ for $n\geq1$,
each carrying its product topology.  Indeed, by
Proposition~\ref{prop-kappa-lambda}, the canonical stage topology
on each cube is the restriction of $\kappa_A=\rho_{\mathrm{pt}}$.
The strict inclusion proved in Example~\ref{exmp-kappa-lambda-strict}
therefore gives
\[
 \kappa_A=\rho_{\mathrm{pt}}\subsetneq\lambda_A=\sigma_A.
\]
In particular, $\rho_{\mathrm{pt}}$ and $\sigma_A$ are distinct
$\ell$-group topologies inducing the same topology on every
$I_n=[-n,n]^{\mathbb N}$, and hence on every order-bounded subset.
The interval topology therefore does not determine a unique
global $\ell$-group topology.

More is true: $(G,u,\rho_{\mathrm{pt}})$ and $(G,u,\sigma_A)$ are not
isomorphic in $\Topulg$, although $\Gamma_t$ sends both to $(A,\tau)$.
To see this, let
\[
 K=\{0\}\cup\{m\mathbf e_m:m\geq1\},
\]
where $\mathbf e_m$ is the $m$th standard basis vector.
This set is compact for the pointwise topology because
$m\mathbf e_m\to0$ pointwise, but it is not order bounded: an upper
bound in $\ell^\infty(\mathbb N)$ would have to dominate the value
$m$ in its $m$th coordinate for every $m$.  Suppose that
\[
 h\colon(G,u,\rho_{\mathrm{pt}})\longrightarrow(G,u,\sigma_A)
\]
were an isomorphism.  Then $h(K)$ would be $\sigma_A$-compact and hence,
by Theorem~\ref{thm-komega-reconstruction}%
\ref{thm-komega-reconstruction-compact}, contained in $[-nu,nu]$ for some
$n$.  Since both $h$ and $h^{-1}$ preserve the order and the
unit $u$, we have $h^{-1}([-nu,nu])=[-nu,nu]$.
Hence $K\subseteq[-nu,nu]$, contradicting the fact
that $K$ is not order bounded.

The same example also shows that $\Gamma_t$ is not full on $\Topulg$.
The identity MV-homomorphism
\[
 1_A\colon\Gamma_t(G,u,\rho_{\mathrm{pt}})
 \longrightarrow\Gamma_t(G,u,\sigma_A)
\]
is continuous because both unit intervals carry the product topology
of $A$.  By the algebraic Mundici correspondence
\eqref{eq-algebraic-Mundici-bijection}, the only
unit-preserving $\ell$-homomorphism extending $1_A$
is $1_G$.  But
\[
 1_G\colon(G,u,\rho_{\mathrm{pt}})
 \longrightarrow(G,u,\sigma_A)
\]
is not continuous, since $\rho_{\mathrm{pt}}\subsetneq\sigma_A$.
Thus $1_A$ is not the image under $\Gamma_t$ of any morphism in
$\Topulg$.  Here $(G,u,\rho_{\mathrm{pt}})$ belongs to $\ULS$ by
Theorem~\ref{lem-unbounded-lift}, but is not Mundici-fine.
The obstruction therefore concerns the unrestricted functor:
its restriction to $\TopulgM$ is full by
Theorem~\ref{thm-topological-correspondence}.
\end{exmp}

\section{Quasitopological reconstruction and the joint-continuity gap}
\label{sec-quasi-reflection}

The preceding sections reconstruct a topological unital abelian
$\ell$-group from a topological MV-algebra while preserving its
unit-interval topology.  We now consider the corresponding problem
when joint continuity of the binary operations is replaced by
separate continuity.  We first construct the least and largest
quasitopological $\ell$-group extensions of every quasitopological
MV-algebra, using initial and final topologies, respectively.
We then ask whether replacing such an extension by its topological
reflection preserves the given topology on the unit interval.
For topological inputs, we also compare these reconstructions
with the bounded-stage final and box-subspace topologies.
Throughout, all underlying topologies are ordinary topologies;
\emph{quasitopological} refers to the continuity of the operations.

\subsection{The quasitopological Mundici correspondence}
\label{sec-quasi-foundation}

The terminology for separate-continuity conditions is not uniform
in the literature.  The semitopological MV-algebras of Zahiri and
Borzooei already have continuous negation \cite{Zahiri2014}, whereas
Najafi, Rezaei, and Kouhestani distinguish the semi-, quasi-, and
paratopological conditions \cite{Najafi2017}.  To fix the objects
under discussion, we adopt the following definitions.
Here separate continuity means continuity after fixing either input.

\begin{defn}
A \emph{quasitopological MV-algebra} is an MV-algebra equipped with
a topology for which $\oplus$ is separately continuous and $\neg$
is continuous.
\end{defn}

\begin{defn}
A \emph{quasitopological unital abelian $\ell$-group} is a unital
abelian $\ell$-group equipped with a topology for which addition,
join, and meet are separately continuous and inversion is continuous.
\end{defn}

We denote their categories by
\[
 \QTopMV\quad\text{and}\quad\QTopulg,
\]
with continuous MV-homomorphisms and continuous unit-preserving
$\ell$-homomorphisms, respectively.  No separation axiom is imposed.
On the MV-algebra side, the derived operation
$a\odot b=\neg(\neg a\oplus\neg b)$ is separately continuous.
The lattice operations are separately continuous as well:
for fixed $b$, the identities
\begin{equation}
\label{eq-quasi-fixed-lattice}
 a\vee b=(a\odot\neg b)\oplus b,
 \qquad a\wedge b=b\odot(\neg b\oplus a)
\end{equation}
express $a\mapsto a\vee b$ and $a\mapsto a\wedge b$ as compositions
of continuous fixed-variable sections.
Here the varying argument occurs only once in each expression;
separate continuity alone does not ensure continuity of arbitrary
terms with repeated variables.

The filter topologies of Luan, Weber, and Yang give quasitopological
MV-algebras that need not be topological \cite{Luan2021}.
The following example exhibits this distinction for both MV-algebras
and unital abelian $\ell$-groups.

\begin{exmp}
\label{exmp-cross-quasitopological}
Consider the case $n=2$ of Example~\ref{exmp-function-algebras},
where $X$ is a two-point discrete space.  The corresponding unital
abelian $\ell$-group is $G=\mathbb R^2$, with coordinatewise addition,
the coordinatewise order, and strong unit $u=(1,1)$; its unit interval is
$A=[0,u]=[0,1]^2$.

Equip $G$ with the \emph{cross topology} $\rho_\times$: a set
$U\subseteq\mathbb R^2$ is open if and only if for every
$x,y\in\mathbb R$, both sections
\[
 \{s\in\mathbb R:(s,y)\in U\},
 \qquad
 \{t\in\mathbb R:(x,t)\in U\}
\]
are open in the usual topology on $\mathbb R$.
If $p,q\colon\mathbb R\longrightarrow\mathbb R$ are continuous
in the usual topology, then
\[
 F\colon(G,\rho_\times)\longrightarrow(G,\rho_\times),
 \qquad F(x,y)=(p(x),q(y)),
\]
is continuous: each horizontal section of $F^{-1}(U)$ is the
inverse image under $p$ of a horizontal section of $U$;
the analogous statement holds for vertical sections, using $q$.
Translations, inversion, and fixed-input lattice operations all
have this form, so $(G,u,\rho_\times)\in\QTopulg$.

Give $A$ the subspace topology $\tau_\times=\rho_\times|_A$ and
the coordinatewise MV-operations.  For fixed $a=(a_1,a_2)\in A$,
the map $x\mapsto x\oplus a$ is the restriction of $F$ with
$p(t)=\min\{t+a_1,1\}$ and $q(t)=\min\{t+a_2,1\}$.
Negation is the restriction of $F$ with $p(t)=q(t)=1-t$.
These restrictions take values in $A$ and are continuous for
$\tau_\times$.  Hence $(A,\tau_\times)\in\QTopMV$.

Neither object is topological.  The set
\[
 U=\mathbb R^2\setminus\{(t,t):t\neq0\}
\]
is a $\rho_\times$-open neighbourhood of zero, since each of its
horizontal and vertical sections omits at most one point.
Put $a_n=(1/n,0)$ and $b_n=(0,1/n)$ for $n\geq1$.
Then $(a_n)_n$ and $(b_n)_n$ converge to zero in both $G$ and $A$:
every $\rho_\times$-open neighbourhood of zero
contains a usual interval about zero on each coordinate axis.
However,
\[
 a_n+b_n=a_n\oplus b_n=(1/n,1/n)\notin U
 \qquad(n\geq2).
\]
Thus neither addition on $G$ nor MV-addition on $A$ is jointly
continuous.  Finally,
$\rho_\times$ refines the Euclidean topology on $\mathbb R^2$,
so both $(G,\rho_\times)$ and its subspace $(A,\tau_\times)$
are Hausdorff.
\end{exmp}

The unit-interval construction
(see \cite[Theorem~4.9]{Zahiri2014}) defines the functor
\begin{equation}
\label{eq-Gammaq-definition}
 \Gamma_{\mathrm q}\colon\QTopulg\longrightarrow\QTopMV,
 \qquad
 \Gamma_{\mathrm q}(G,u,\rho)=([0,u],\rho|_{[0,u]}),
 \qquad
 \Gamma_{\mathrm q}h=h|_{[0,u]}.
\end{equation}
Indeed, fixed-input truncated addition is the restriction of a
translation followed by meet with $u$.  Negation is the restriction
of inversion followed by a translation.
Morphisms restrict continuously.

For $(A,\tau)\in\QTopMV$, write $\Xi A=(G,u)$ and identify
$A=[0,u]$.  As in Definition~\ref{def-admissible-topology}, call a
topology $\rho$ on $G$ \emph{quasi-admissible} for $(A,\tau)$ if
$(G,u,\rho)\in\QTopulg$ and $\rho|_A\subseteq\tau$.
We call such a topology \emph{exact} if $\rho|_A=\tau$.

The translated truncations of
\eqref{eq-translated-truncation} remain defined by
$t_g(x)=\pi_0(x+g)$, and we again let $\kappa_A$ be the initial
topology induced by
\begin{equation}
\label{eq-quasi-initial-family}
 \bigl(t_g\colon G\longrightarrow(A,\tau)\bigr)_{g\in G}.
\end{equation}
For topological inputs this is the topology of
Proposition~\ref{prop-kappa-lambda}.  Its other initial
characterizations there use joint continuity; we do not assume
them for arbitrary quasi inputs.

As we shall see, $\kappa_A$ makes $(G,u)$ into a quasitopological
unital abelian $\ell$-group and induces $\tau$ on $A$.
We first establish the following identity.

\begin{lem}[Carry identity]
\label{lem-quasi-carry}
For $g\in G$, $a\in A$, and $k\in\mathbb Z$,
\begin{equation}
\label{eq-quasi-carry}
 \pi_k(g+a)=\pi_k(g)\oplus\bigl(\pi_{k-1}(g)\odot a\bigr).
\end{equation}
In particular,
\begin{equation}
\label{eq-translated-truncation-interval}
 t_g(a)=\pi_0(g)\oplus\bigl(\pi_{-1}(g)\odot a\bigr).
\end{equation}
\end{lem}

\begin{proof}
By Lemma~\ref{lem-coordinate-identities}\ref{coordinate-identities-interval},
the coordinates of $a\in[0,u]$ are $u$ for $j<0$, $a$ for
$j=0$, and $0$ for $j>0$.  Specializing
Lemma~\ref{lem-coordinate-identities}\ref{coordinate-identities-addition} therefore gives
\[
 \pi_k(g+a)
 =\bigoplus_{i+j=k-1}\bigl(\pi_i(g)\odot\pi_j(a)\bigr)
 =\bigl(\pi_{k-1}(g)\odot a\bigr)
   \oplus\bigoplus_{i\geq k}\pi_i(g).
\]
Here the iterated MV-sum has only finitely many nonzero terms by
Lemma~\ref{lem-coordinate-identities}\ref{eq-bounded-coordinate-tails}.  Applying
$\pi_i(g)\oplus\pi_{i+1}(g)=\pi_i(g)$ from
Lemma~\ref{lem-coordinate-identities}\ref{coordinate-identities-absorption} successively from the right
reduces the last sum to $\pi_k(g)$.  This proves
\eqref{eq-quasi-carry}; take $k=0$ for the final assertion.
\end{proof}

\begin{prop}[Minimal quasitopological reconstruction]
\label{thm-quasi-minimal-reconstruction}
For every $(A,\tau)\in\QTopMV$, the topology $\kappa_A$ on
$\Xi A=(G,u)$ satisfies
\begin{equation}
\label{eq-quasi-minimal-exactness}
 \kappa_A|_A=\tau,
 \qquad (G,u,\kappa_A)\in\QTopulg.
\end{equation}
It is the coarsest topology on $G$ with these two properties.
\end{prop}

\begin{proof}
We first show that, for every $g\in G$, the restriction
\[
 t_g|_A\colon(A,\tau)\longrightarrow(A,\tau)
\]
is continuous.  This follows from
\eqref{eq-translated-truncation-interval}: the elements $\pi_0(g)$
and $\pi_{-1}(g)$ are fixed, so separate continuity of $\oplus$
and $\odot$ suffices.
By initiality, the inclusion
$i\colon(A,\tau)\longhookrightarrow(G,\kappa_A)$ is continuous,
giving $\kappa_A|_A\subseteq\tau$.  Conversely, the continuous
map $t_0\colon(G,\kappa_A)\longrightarrow(A,\tau)$ restricts to
$1_A$, giving $\tau\subseteq\kappa_A|_A$.
Thus $\kappa_A|_A=\tau$.

We next show that $(G,u,\kappa_A)\in\QTopulg$.
By initiality, continuity of a map into $(G,\kappa_A)$ can be
checked by composing it with every $t_g$.
For $g,h,x\in G$, translations and inversion satisfy
\begin{equation}
\label{eq-translated-truncation-symmetries}
 t_g(x+h)=t_{g+h}(x),
 \qquad t_g(-x)=\neg t_{u-g}(x).
\end{equation}
The second identity follows from
$\pi_0(-y)=\neg\pi_{-1}(y)$ in
Lemma~\ref{lem-coordinate-identities}\ref{coordinate-identities-inversion}, with $y=x-g$, and
$\pi_{-1}(z)=\pi_0(z+u)$.
For fixed $g,h$, both right-hand sides are continuous in $x$
by the definition of $\kappa_A$ and continuity of negation on
$(A,\tau)$.  Initiality therefore gives continuity of inversion
and every translation.  Translation by $-h$ is the inverse of
translation by $h$, so all translations are homeomorphisms.
In particular, addition is separately continuous.

Distributivity makes $\pi_0(x)=(x\wedge u)\vee0$ preserve binary
joins and meets, and translations preserve both operations.
Consequently, for fixed $g,h\in G$,
\begin{equation}
\label{eq-translated-truncation-lattice}
 t_g(x\vee h)=t_g(x)\vee t_g(h),
 \qquad t_g(x\wedge h)=t_g(x)\wedge t_g(h).
\end{equation}
Here $t_g(h)\in A$ is fixed, so the right-hand sides are continuous
in $x$ by the fixed-variable continuity in
\eqref{eq-quasi-fixed-lattice}.  Initiality gives continuity of
$x\mapsto x\vee h$ and $x\mapsto x\wedge h$ on $(G,\kappa_A)$.
Commutativity then gives separate continuity of both lattice
operations.

Finally, let $\rho$ be a topology on $G$ such that
$(G,u,\rho)\in\QTopulg$ and $\rho|_A=\tau$.
For each $g\in G$, the map
\[
 x\mapsto t_g(x)=((x+g)\wedge u)\vee0
\]
is continuous as a map $(G,\rho)\longrightarrow(G,\rho)$,
being a composition of a translation and fixed-input lattice
operations.  Its values lie in $A$, and $\tau=\rho|_A$,
so it is also continuous as a map
$(G,\rho)\longrightarrow(A,\tau)$.
Since $\kappa_A$ is the initial topology induced by all these
maps, $\kappa_A\subseteq\rho$, proving minimality.
\end{proof}

We now construct the largest quasi-admissible topology.
For $(A,\tau)\in\QTopMV$ and $\Xi A=(G,u)$, define
$\sigma_A^{\mathrm q}$ to be the final topology on $G$ induced
by the translated unit-interval maps
\[
 \bigl(i_g\colon(A,\tau)\longrightarrow G\bigr)_{g\in G},
 \qquad i_g(a)=g+a.
\]
Explicitly, for $U\subseteq G$,
\begin{equation}
\label{eq-quasi-translated-final-test}
 U\in\sigma_A^{\mathrm q}
 \iff
 \{a\in A:g+a\in U\}\in\tau
 \quad\text{for every }g\in G.
\end{equation}

\begin{prop}[Maximal quasitopological reconstruction]
\label{prop-quasi-translated-final}
Let $(A,\tau)\in\QTopMV$ and $\Xi A=(G,u)$. Then
$\sigma_A^{\mathrm q}$ is the finest quasi-admissible topology on $G$.
\end{prop}

\begin{proof}
Let $\rho$ be any quasi-admissible topology on $G$.
By quasi-admissibility, the interval inclusion
$(A,\tau)\longhookrightarrow(G,\rho)$ is continuous.
Separate continuity of addition on $(G,\rho)$ also makes
every translation continuous.  Composing the inclusion with
translation by $g$ therefore makes every
$i_g\colon(A,\tau)\longrightarrow(G,\rho)$ continuous.
By the defining property of the final topology,
$\rho\subseteq\sigma_A^{\mathrm q}$.
It remains to show that $\sigma_A^{\mathrm q}$ is itself
quasi-admissible.  The inclusion $i_0$ is continuous by
definition, so we need only check the required operations.

For $h\in G$, the translation $T_h(x)=h+x$ satisfies
$T_hi_g=i_{g+h}$, so finality makes $T_h$ continuous.
Its inverse is $T_{-h}$, hence all translations are homeomorphisms.
Inversion is continuous by the same finality test, since
\[
 -i_g(a)=i_{-g-u}(\neg a)
\]
and negation is continuous on $(A,\tau)$.
Thus addition is separately continuous and inversion is continuous
for $\sigma_A^{\mathrm q}$.

For the lattice operations, fix $g,h\in G$.  Set
$d=h-g$ and $c=\pi_0(d)=(d\vee0)\wedge u\in A$.
We claim that, for $a\in A$,
\begin{equation}
\label{eq-quasi-translated-lattice-sections}
 (g+a)\vee h=(g\vee h)+(a\odot\neg c),
 \qquad
 (g+a)\wedge h=(g\wedge h)+(a\wedge c).
\end{equation}
To verify these identities, apply $z=(z\vee0)+(z\wedge0)$
to $z=a\wedge d$.  Distributivity and $0\leq a\leq u$ give
\[
\begin{aligned}
 a\wedge d
 &=\bigl((a\wedge d)\vee0\bigr)
   +\bigl((a\wedge d)\wedge0\bigr)
 =\bigl(a\wedge(d\vee0)\bigr)+(d\wedge0)\\
 &=\bigl(a\wedge u\wedge(d\vee0)\bigr)+(d\wedge0)
 =(a\wedge c)+(d\wedge0),
\end{aligned}
\]
where the second equality uses $a\geq0$, and the third uses $a\leq u$.
The identity $(a\vee d)+(a\wedge d)=a+d$ then gives
\[
\begin{aligned}
 a\vee d
 &=a+d-(a\wedge d)
 =a+d-\bigl((a\wedge c)+(d\wedge0)\bigr)\\
 &=(d\vee0)+\bigl(a-(a\wedge c)\bigr)
 =(d\vee0)+\bigl((a-c)\vee0\bigr),
\end{aligned}
\]
where we use $d-(d\wedge0)=d\vee0$ and
$a-(a\wedge c)=(a-a)\vee(a-c)=(a-c)\vee0$.
The definition of $\odot$ and \eqref{eq-Gamma-def} give
$a\odot\neg c=(a-c)\vee0$.
Adding $g$ to the formulas for $a\vee d$ and $a\wedge d$ now proves
\eqref{eq-quasi-translated-lattice-sections}.

Since $c$ is fixed, the maps $a\mapsto a\odot\neg c$ and
$a\mapsto a\wedge c$ are continuous on $(A,\tau)$.
The two right-hand sides of
\eqref{eq-quasi-translated-lattice-sections} therefore factor
through the continuous translated interval maps $i_{g\vee h}$
and $i_{g\wedge h}$, respectively.
Finality now gives continuity of $x\mapsto x\vee h$ and
$x\mapsto x\wedge h$ on $(G,\sigma_A^{\mathrm q})$.
By commutativity, both lattice operations are separately continuous.
Thus $\sigma_A^{\mathrm q}$ is quasi-admissible.
\end{proof}

\begin{prop}[The quasitopological Mundici adjunction]
\label{prop-quasi-adjunction}
\label{thm-quasi-exactness}
For each $(A,\tau)\in\QTopMV$, write $\Xi A=(G_A,u_A)$.
The assignments
\begin{equation}
\label{eq-Xiq-definition}
 \Xi_{\mathrm q}(A,\tau)=(G_A,u_A,\sigma_A^{\mathrm q}),
 \qquad \Xi_{\mathrm q}f=\Xi f
\end{equation}
define a left adjoint $\Xi_{\mathrm q}$ to $\Gamma_{\mathrm q}$.
Moreover,
\begin{equation}
\label{eq-quasi-exactness}
 \kappa_A\subseteq\sigma_A^{\mathrm q},
 \qquad \kappa_A|_A=(\sigma_A^{\mathrm q})|_A=\tau.
\end{equation}
The unit is a natural isomorphism, and $\Xi_{\mathrm q}$ is fully faithful.
\end{prop}

\begin{proof}
By Proposition~\ref{prop-quasi-translated-final},
$(G_A,u_A,\sigma_A^{\mathrm q})\in\QTopulg$.
To prove the adjunction, let $Y=(H,v,\rho)\in\QTopulg$ and let
\[
 f\colon(A,\tau)\longrightarrow\Gamma_{\mathrm q}Y
\]
be a continuous MV-homomorphism.
By the algebraic Mundici correspondence
\eqref{eq-algebraic-Mundici-bijection}, $f$ extends uniquely to a
unit-preserving $\ell$-homomorphism $h\colon G_A\longrightarrow H$.
For every $g\in G_A$,
\[
 h i_g(a)=h(g)+f(a)\qquad(a\in A).
\]
The right-hand side is continuous as a map from $(A,\tau)$ to
$(H,\rho)$: it is $f$, followed by the interval inclusion and
translation by $h(g)$.
Since $\sigma_A^{\mathrm q}$ is final for the maps $i_g$,
$h\colon(G_A,\sigma_A^{\mathrm q})\longrightarrow(H,\rho)$
is continuous.
Conversely, every continuous unit-preserving $\ell$-homomorphism
between these objects restricts to a continuous MV-homomorphism
from $(A,\tau)$, because the interval inclusion $i_0$ is continuous.

The algebraic Mundici hom-set bijection therefore restricts to
continuous morphisms.
In particular, applying this to a continuous MV-homomorphism
$f\colon(A,\tau)\longrightarrow(B,\tau_B)$ followed by the
interval inclusion into $\Xi_{\mathrm q}(B,\tau_B)$ makes
$\Xi f$ continuous.
Identities, compositions, and naturality are inherited from the
classical Mundici functors.  Thus the assignments define the stated
adjunction, with the classical underlying unit and counit.

Finally, Proposition~\ref{thm-quasi-minimal-reconstruction} makes
$\kappa_A$ quasi-admissible, so maximality in
Proposition~\ref{prop-quasi-translated-final} gives
$\kappa_A\subseteq\sigma_A^{\mathrm q}$.  Restricting to $A$ gives
\[
 \tau=\kappa_A|_A
 \subseteq(\sigma_A^{\mathrm q})|_A\subseteq\tau.
\]
The unit is therefore a homeomorphism and hence a natural
isomorphism.  The fully faithful left-adjoint criterion used in
Theorem~\ref{thm-topological-correspondence} completes the proof.
\end{proof}

For $X=(G,u,\rho)\in\QTopulg$, transport
$\sigma_{\Gamma_{\mathrm q}X}^{\mathrm q}$ to $G$ along the algebraic
counit, and denote it by $\widetilde\sigma_G^{\mathrm q}$.
The transport is defined in the same way as in
\eqref{eq-transported-maximal-topology}.
Call $X$ \emph{quasi-Mundici-fine} if
$\rho=\widetilde\sigma_G^{\mathrm q}$, and write $\QTopulgM$ for
the full subcategory of these objects.
Equivalently, $\rho$ is the finest topology on $G$ inducing the given
unit-interval topology and making the group an object of $\QTopulg$.

\begin{cor}[Quasitopological Mundici correspondence]
\label{thm-quasi-correspondence}
The adjunction restricts to an equivalence
\[
 \QTopMV\simeq\QTopulgM.
\]
The category $\QTopulgM$ is a replete coreflective full
subcategory of $\QTopulg$, with coreflector the corestriction of
$\Xi_{\mathrm q}\Gamma_{\mathrm q}$ and coreflection morphism its counit.
\end{cor}

\begin{proof}
By Proposition~\ref{prop-quasi-adjunction}, the unit is a natural
isomorphism.  For $X=(G,u,\rho)\in\QTopulg$, the analogue of
\eqref{eq-transported-counit-factorization} expresses the counit
as a homeomorphism followed by the continuous identity
\[
 1_G\colon(G,u,\widetilde\sigma_G^{\mathrm q})
 \longrightarrow(G,u,\rho).
\]
Hence the counit is an isomorphism if and only if
$\rho=\widetilde\sigma_G^{\mathrm q}$, that is, if and only if
$X$ is quasi-Mundici-fine.
The essential-image and coreflection arguments of
Theorem~\ref{thm-topological-correspondence} now apply, with
maximality supplied by Proposition~\ref{prop-quasi-translated-final} and interval exactness by
Proposition~\ref{prop-quasi-adjunction}.
\end{proof}

We now study the relationship between the topological and
quasitopological Mundici correspondences.  Let
\[
 J\colon\Topulg\longhookrightarrow\QTopulg
 \quad\text{and}\quad
 K\colon\TopMV\longhookrightarrow\QTopMV
\]
be the inclusions.  For a topology $\rho$ on a unital abelian
$\ell$-group $(G,u)$, define
\begin{equation}
\label{eq-lgroup-reflection}
 R_\ell(\rho)=\bigvee\{\theta\subseteq\rho:
                   \theta\text{ is an $\ell$-group topology on }G\}.
\end{equation}
For $(A,\tau)\in\QTopMV$, define
\begin{equation}
\label{eq-MV-reflection}
 R_{\mathrm{MV}}(\tau)
 =\bigvee\{\mu\subseteq\tau:(A,\mu)\in\TopMV\}.
\end{equation}
Both joins are taken in the lattice of ordinary topologies.
The join and pullback arguments of
Lemmas~\ref{lem-finest-admissible} and~\ref{lem-pullback-topology},
which apply equally to MV-algebras, show that the assignments
\[
 (G,u,\rho)\mapsto(G,u,R_\ell(\rho)),\qquad (A,\tau)\mapsto(A,R_{\mathrm{MV}}(\tau))
\]
define reflectors
\[
 \begin{aligned}
  R_\ell&\colon\QTopulg\longrightarrow\Topulg,
  &\qquad R_\ell&\dashv J,\\
  R_{\mathrm{MV}}&\colon\QTopMV\longrightarrow\TopMV,
  &\qquad R_{\mathrm{MV}}&\dashv K.
 \end{aligned}
\]
They leave the underlying homomorphisms unchanged, and their units
are the continuous identity maps to the reflected topologies.

The maximal reconstructions commute with these reflectors,
as we see below.

\begin{prop}[Reflection of maximal reconstruction]
\label{thm-quasi-reflection-factorization}
With the fixed classical reconstruction $\Xi$, the square
\begin{equation}
\label{eq-reflection-square}
\begin{tikzcd}[column sep=4.5em,row sep=3em]
 \QTopMV \arrow[r,"\Xi_{\mathrm q}"] \arrow[d,"R_{\mathrm{MV}}"']
 & \QTopulg \arrow[d,"R_\ell"] \\
 \TopMV \arrow[r,"\Xi_t"'] & \Topulg
\end{tikzcd}
\end{equation}
commutes. For topological inputs, this and
Theorem~\ref{thm-lambda-envelope} give
\begin{equation}
\label{eq-quasi-reflection-topologies}
 R_\ell(\lambda_A)=R_\ell(\sigma_A^{\mathrm q})=\sigma_A.
\end{equation}
Restricting the commuting square to topological inputs gives
$\Xi_t=R_\ell\Xi_{\mathrm q}K$.
\end{prop}

\begin{proof}
Let $(A,\tau)\in\QTopMV$ and write $(G_A,u_A)=\Xi A$.
In this proof and the following corollary, we write
$\sigma_A\bigl(R_{\mathrm{MV}}(\tau)\bigr)$ for the Mundici topology
on $G_A$ associated with the topological MV-algebra
$(A,R_{\mathrm{MV}}(\tau))$.
We need to prove that
$R_\ell(\sigma_A^{\mathrm q})=\sigma_A\bigl(R_{\mathrm{MV}}(\tau)\bigr)$.
We compare the families defining the two topologies.
For any $\ell$-group topology $\rho$ on $G_A$,
\[
 \rho\subseteq\sigma_A^{\mathrm q}
 \iff \rho|_A\subseteq\tau
 \iff \rho|_A\subseteq R_{\mathrm{MV}}(\tau).
\]
The first equivalence uses maximality from
Proposition~\ref{prop-quasi-translated-final}:
$\rho|_A\subseteq\tau$ makes $\rho$ quasi-admissible, while the
other implication follows by restriction, using
$(\sigma_A^{\mathrm q})|_A=\tau$ from
Proposition~\ref{prop-quasi-adjunction}.
For the second, $\rho|_A$ is a topological MV-algebra topology,
and $R_{\mathrm{MV}}(\tau)$ is the finest such topology contained
in $\tau$, by \eqref{eq-MV-reflection}.
Thus the joins defining $R_\ell\bigl(\sigma_A^{\mathrm q}\bigr)$ and
$\sigma_A\bigl(R_{\mathrm{MV}}(\tau)\bigr)$ are taken over the same family,
proving the required equality.
Both paths send a morphism $f$ to $\Xi f$, since the reflectors
leave underlying homomorphisms unchanged.
Hence the square commutes on morphisms as well.

For topological $\tau$, one has $R_{\mathrm{MV}}(\tau)=\tau$, so
$R_\ell(\sigma_A^{\mathrm q})=\sigma_A$.
Together with $R_\ell(\lambda_A)=\sigma_A$ from
Theorem~\ref{thm-lambda-envelope}\ref{thm-lambda-envelope-maximality}, this proves
\eqref{eq-quasi-reflection-topologies}.
Finally, restrict the commuting square along $K$ and use
$R_{\mathrm{MV}}K=1_{\TopMV}$ to obtain
$\Xi_t=R_\ell\Xi_{\mathrm q}K$.
\end{proof}

\begin{cor}
\label{cor-quasi-reflection-exactness}
For every $(A,\tau)\in\QTopMV$,
\[
 \left.R_\ell\bigl(\sigma_A^{\mathrm q}\bigr)\right|_A
 =R_{\mathrm{MV}}(\tau).
\]
In particular, this reflection induces the original topology $\tau$
on the interval if and only if $(A,\tau)\in\TopMV$.
\end{cor}

\begin{proof}
Proposition~\ref{thm-quasi-reflection-factorization} identifies the
reflected topology with $\sigma_A\bigl(R_{\mathrm{MV}}(\tau)\bigr)$.
Applying Theorem~\ref{thm-topological-representation} to
$(A,R_{\mathrm{MV}}(\tau))$ now gives
\[
 \left.\sigma_A\bigl(R_{\mathrm{MV}}(\tau)\bigr)\right|_A
 =R_{\mathrm{MV}}(\tau).
\]
Thus interval equality with $\tau$ holds
precisely when $\tau=R_{\mathrm{MV}}(\tau)$, equivalently when
$(A,\tau)$ is a topological MV-algebra.
\end{proof}

The minimal reconstruction behaves differently.

\begin{rem}
\label{rem-minimal-quasi-reflection}
For every $(A,\tau)\in\QTopMV$,
\begin{equation}
\label{eq-minimal-quasi-reflection}
 R_\ell(\kappa_A)|_A=\tau
 \iff R_\ell(\kappa_A)=\kappa_A
 \iff \kappa_A\text{ is an $\ell$-group topology}.
\end{equation}
Indeed, if the reflection retains $\tau$, minimality in
Proposition~\ref{thm-quasi-minimal-reconstruction} gives
$\kappa_A\subseteq R_\ell(\kappa_A)$, while reflection gives the
reverse inclusion.  The other implications follow from exactness
of $\kappa_A$ and the definition of the reflector.

The Hausdorff Boolean MV-algebra in
Example~\ref{exmp-kappa-not-group} is topological, but its
coordinate topology $\kappa_A$ does not make group addition
jointly continuous.  Hence
$R_\ell(\kappa_A)|_A\subsetneq\tau$.
By contrast, for the product cube $A=[0,1]^{\mathbb N}$ in
Example~\ref{exmp-product-cubes},
\[
 R_\ell(\kappa_A)=\kappa_A=\rho_{\mathrm{pt}}\subsetneq\sigma_A.
\]
Here $\kappa_A$ is already an $\ell$-group topology, but its
reflection does not recover $\sigma_A$.
Thus reflection of the minimal exact quasi reconstruction need not
retain the interval topology, and even when it does, it need not
produce $\Xi_t$.
\end{rem}

\subsection{A comparison of the topologies}
\label{sec-quasi-topological-reflection}

Throughout this subsection, the input $(A,\tau)$ is topological.
We compare the bounded-stage final topology $\lambda_A$ with
the maximal quasi topology $\sigma_A^{\mathrm q}$, introduce
the box-subspace topology $\nu_A$, and establish the full chain
of inclusions and its possible strictness.
The first comparison follows by viewing the bounded-stage final
topology as a quasi reconstruction.

\begin{prop}[Quasitopological reconstruction by bounded stages]
\label{prop-lambda-quasi-reconstruction}
Let $(A,\tau)\in\TopMV$ and $\Xi A=(G,u)$.
The bounded-stage final topology $\lambda_A$ of
Definition~\ref{def-Mundici-final-topology} satisfies
\[
 (G,u,\lambda_A)\in\QTopulg,
 \qquad \lambda_A|_A=\tau.
\]
In particular,
\begin{equation}
\label{eq-lambda-quasi-comparison}
 \lambda_A\subseteq\sigma_A^{\mathrm q}.
\end{equation}
\end{prop}

\begin{proof}
Proposition~\ref{prop-lambda-quasitopological} gives continuous
translations and inversion for $\lambda_A$.
For fixed $h\in G$, choose an integer $C\geq1$ such that
$h\in I_C$, using the strong-unit property.
For every $n\geq1$, the maps $x\mapsto x\vee h$ and
$x\mapsto x\wedge h$ restrict to continuous maps
\[
 (I_n,\tau_n)\longrightarrow
 (I_{\max\{n,C\}},\tau_{\max\{n,C\}})
\]
by the stage inclusions and lattice continuity in
Lemma~\ref{lem-stage-topology}.
The inclusion of the target stage into $(G,\lambda_A)$ is
continuous by Definition~\ref{def-Mundici-final-topology}.
Thus both maps are continuous on every stage as maps into
$(G,\lambda_A)$, and finality gives their continuity on $G$.
Commutativity gives separate continuity of the lattice operations,
so $(G,u,\lambda_A)\in\QTopulg$.
Proposition~\ref{prop-lambda-exact-space} gives
$\lambda_A|_A=\tau$.  Hence $\lambda_A$ is quasi-admissible,
and maximality in Proposition~\ref{prop-quasi-translated-final} gives
\eqref{eq-lambda-quasi-comparison}.
\end{proof}

Thus $\lambda_A$ always gives a quasi reconstruction, although
Example~\ref{exmp-lambda-not-group} shows that its addition need
not be jointly continuous.
Both $\lambda_A$ and $\sigma_A^{\mathrm q}$ induce $\tau$ on $A$,
but Proposition~\ref{prop-lgroup-extension-bounded-stages} does
not force their restrictions to agree on every bounded stage:
its proof uses joint continuity of addition, which is not
assumed for quasi reconstructions.
Using the translated-interval criterion
\eqref{eq-quasi-translated-final-test}, the following example shows
that their bounded-stage restrictions can indeed differ.
Corollary~\ref{cor-quasi-full-chain} will collect the full topology
comparison and the possible strictness of each inclusion.

\begin{exmp}[The maximal quasi topology may change bounded stages]
\label{exmp-quasi-stage-strict-left}
Let $E=\mathbb R$ or $\mathbb Q$, and take
$A=(E\cap[0,1])^{\mathbb N}$ with its product topology $\tau$.
In both cases $(A,\tau)\in\TopMV$.
These are the real and rational product cubes of
Examples~\ref{exmp-kappa-lambda-strict} and~\ref{exmp-lambda-not-group}.
In each case $G$ is the group of bounded $E$-valued sequences,
$u=\mathbf1$, and $\lambda_A|_{I_2}=\tau_2$ is the product topology on
$I_2=(E\cap[-2,2])^{\mathbb N}$.
Let $\mathbf e_m$ be the $m$th coordinate vector and put
\[
 F=\left\{\frac32\mathbf e_m:m\geq1\right\}\subseteq I_2.
\]
The sequence $((3/2)\mathbf e_m)_{m\geq1}$ converges to $0$ in $\tau_2$, but
$0\notin F$; hence $F$ is not $\tau_2$-closed.

We claim that $F$ is $\sigma_A^{\mathrm q}$-closed.
Every translated interval is
$g+A=\prod_{k\geq1}(E\cap[g_k,g_k+1])$.
If $0\in g+A$, each upper endpoint is at most $1$, so
$(g+A)\cap F=\varnothing$.  If $0\notin g+A$, some coordinate
$k$ has $0\notin[g_k,g_k+1]$; hence $(g+A)\cap F$ contains
at most $(3/2)\mathbf e_k$.
In either case $(g+A)\cap F$ is finite, and therefore closed
in the translated copy of the Hausdorff space $(A,\tau)$.
By \eqref{eq-quasi-translated-final-test}, $F$ is closed in
$\sigma_A^{\mathrm q}$.  Consequently,
\[
 \lambda_A|_{I_2}\subsetneq\sigma_A^{\mathrm q}|_{I_2}.
\]
In particular, quasi maximality need not preserve the forced
bounded-stage topologies of $\ell$-group extensions.
\end{exmp}

Thus the bounded-stage uniqueness of
Proposition~\ref{prop-lgroup-extension-bounded-stages} depends
essentially on joint continuity of addition: it fails even for
exact quasi reconstructions of topological inputs.

We next introduce the box-subspace topology $\nu_A$, which lies
between $\sigma_A$ and $\lambda_A$.
Let $(A,\tau)\in\TopMV$ and $\Xi A=(G,u)$.
Give $A^{\mathbb Z}$ the box topology $\tau_{\mathrm{box}}$,
whose basic open sets are $\prod_{k\in\mathbb Z}U_k$ with every
$U_k$ open in $A$.
The homeomorphism $\Theta$ in
\eqref{eq-coordinate-negation-homeomorphism} shows that the
coordinate embeddings $\pi$ and $\delta$ from
Subsection~\ref{subsec-bounded-control-machinery} induce the same
initial topology on $G$.  We denote this topology by $\nu_A$
and call it the \emph{box-subspace topology}.
The coordinate boxes form a zero-neighbourhood base for $\nu_A$.
By Lemma~\ref{lem-coordinate-identities}\ref{coordinate-identities-inversion}
and \eqref{eq-Mundici-deviation}, we have
$\delta_k(-g)=\delta_{-k-1}(g)$.
Thus inversion is continuous for $\nu_A$, and the symmetric
coordinate boxes also form a zero-neighbourhood base.

The following proposition shows that $\nu_A$ makes $(G,u)$ into a 
quasitopological $\ell$-group and that its group reflection is $\sigma_A$.
This description is retrospective; the reconstruction in
Section~\ref{sec-group-envelope} does not use the following proposition.

\begin{prop}[The coordinate-box quasitopological reconstruction]
\label{prop-box-quasi-reconstruction}
Let $(A,\tau)\in\TopMV$ and $\Xi A=(G,u)$.
Then $(G,u,\nu_A)\in\QTopulg$.
Moreover, its lattice operations are jointly continuous, and
\[
 \kappa_A\subseteq\sigma_A=\gamma_A
 \subseteq\nu_A\subseteq\lambda_A.
\]
Therefore, the group and $\ell$-group reflections of $\nu_A$
coincide with the Mundici topology:
\begin{equation}
\label{eq-box-quasi-reflections}
 R_{\mathrm{grp}}(\nu_A)
 =R_\ell(\nu_A)=\sigma_A.
\end{equation}
\end{prop}

\begin{proof}
\textbf{Step 1.} We prove continuity of every fixed translation.
For $a\in A$, the carry identity
\eqref{eq-quasi-carry} and the equality $h-a=h+\neg a-u$ give
\[
 \begin{aligned}
 \pi_k(h+a)&=\pi_k(h)\oplus(\pi_{k-1}(h)\odot a),\\
 \pi_k(h-a)&=\pi_{k+1}(h)\oplus(\pi_k(h)\odot\neg a).
 \end{aligned}
\]
Write a fixed $x\in G$ as a finite signed sum of elements of $A$
using \eqref{eq-deviation-decomposition}.
If this sum has $r$ terms, each iteration changes the coordinate
index by at most one.  After replacing negative-index coordinates
by their complements, we therefore obtain continuous maps
$T_{x,k}\colon A^{2r+1}\longrightarrow A$ such that
\[
 \delta_k(x+h)=T_{x,k}
   (\delta_{k-r}(h),\ldots,\delta_{k+r}(h))
 \qquad(k\in\mathbb Z).
\]
The same $r$ works for every $k$.
Here continuity of the finite terms uses the topological MV hypothesis.
Fix $h_0\in G$ and a basic $\nu_A$-open neighbourhood
\[
 D(\mathcal U)=\bigcap_{k\in\mathbb Z}\delta_k^{-1}(U_k)
\]
of $x+h_0$, where each $U_k$ is open in $A$.
For each $k$, continuity of $T_{x,k}$ gives open neighbourhoods
$V_{k,j}$ of $\delta_j(h_0)$, for $|j-k|\leq r$, such that
\[
 T_{x,k}\left(\prod_{j=k-r}^{k+r}V_{k,j}\right)\subseteq U_k.
\]
For each input index $j$, put
\[
 V_j=\bigcap_{|k-j|\leq r}V_{k,j}.
\]
This is an open neighbourhood of $\delta_j(h_0)$, since the
intersection is finite.  Thus $D((V_j)_j)$ is a
$\nu_A$-open neighbourhood of $h_0$, and the coordinate
formula above gives $x+D((V_j)_j)\subseteq D(\mathcal U)$.
Thus all translations are continuous, and applying the argument to
$-x$ makes them homeomorphisms.

\textbf{Step 2.} We prove joint continuity of the lattice operations.
Lemma~\ref{lem-coordinate-identities}\ref{coordinate-identities-lattice}
and \eqref{eq-Mundici-deviation} give
\[
 \delta_k(x\vee y)=
 \begin{cases}
  \delta_k(x)\vee\delta_k(y),&k\geq0,\\
  \delta_k(x)\wedge\delta_k(y),&k<0,
 \end{cases}
\]
with $\vee$ and $\wedge$ interchanged for $\delta_k(x\wedge y)$.
The induced operations on the box product are jointly continuous:
for each output-coordinate neighbourhood, choose one neighbourhood
for each of the two inputs by joint lattice continuity on $A$, then
assemble the choices into two box neighbourhoods.
Initiality of $\delta$ proves joint lattice continuity on $G$.
Together with translations and inversion this gives membership in
$\QTopulg$.

\textbf{Step 3.} We establish the topology comparisons and the
reflection identities.
The box topology refines the product topology, so the coordinate
characterization of $\kappa_A$ gives
$\kappa_A\subseteq\nu_A$.
Every symmetric coordinate box is $\lambda_A$-open by
Lemma~\ref{lem-deviation-decomposition}, and its translates are
$\lambda_A$-open by Proposition~\ref{prop-lambda-quasitopological}.
As these translates form a base for
$\nu_A$, we have $\nu_A\subseteq\lambda_A$.
Since $\kappa_A|_{I_n}=\lambda_A|_{I_n}=\tau_n$ by
\eqref{eq-lambda-stage-restriction}, we also have
$\nu_A|_{I_n}=\tau_n$ for every $n\geq1$.

By Lemma~\ref{lem-coordinate-box-refinement}, every
$\gamma_A$-zero-neighbourhood contains a symmetric coordinate box;
hence $\gamma_A\subseteq\nu_A$.
The inclusions $\gamma_A\subseteq\nu_A\subseteq\lambda_A$
and maximality of $\gamma_A$ therefore give
$R_{\mathrm{grp}}(\nu_A)=\gamma_A$.
By Corollary~\ref{cor-reconstruction-topology-chain},
$\gamma_A=\sigma_A$ is an $\ell$-group topology.  Hence it is also the greatest
$\ell$-group topology below $\nu_A$, proving
$R_\ell(\nu_A)=\sigma_A$.
\end{proof}

The box-subspace construction thus gives an explicit upper bound
for all $\ell$-group extensions that is contained in $\lambda_A$.
The following example shows that it still need not supply joint
continuity of addition.

\begin{exmp}[The coordinate-box topology need not be a group topology]
\label{exmp-box-not-group}
Let $(B,\tau_B)$ be the Hausdorff topological Boolean MV-algebra
constructed in Example~\ref{exmp-kappa-not-group}: its elements are
the finite and cofinite subsets of
$S=\mathbb N_{\geq1}\times\{0,1\}$.
Its open zero-neighbourhood base is
\[
 \mathcal U_m=
 \left\{\bigcup_{n\in F}P_n:
 F\subseteq\{m,m+1,\ldots\}\text{ is finite}\right\},
 \qquad P_n=\{(n,0),(n,1)\}.
\]
Recall that $b_n=\mathbf1_{P_n}\to0$ in $B$, whereas
$c_n=\mathbf1_{\{(n,0)\}}\leq b_n$ does not belong to
$\mathcal U_1$ for any $n$.  This topology is metrizable, since
its Hausdorff Boolean group has a countable open-subgroup base.
Put $A=B^{\mathbb N_{\geq1}}$ with its product topology $\tau$.
Then $A$ is a Hausdorff metrizable topological Boolean MV-algebra.

Let $H=\mathbb Z\mathbf1+\mathbb Z^{(S)}$ be the Mundici group
of $B$ from Example~\ref{exmp-kappa-not-group}.  The Mundici group
of $A$ may be realized as
\[
 G=\left\{(x_j)_{j\geq1}\in H^{\mathbb N_{\geq1}}:
       \text{for some }M\geq1,\ 
       |x_j|\leq M\mathbf1\text{ for every }j\right\},
 \qquad u=(\mathbf1)_{j\geq1}.
\]
Indeed, this is a unital abelian $\ell$-group with unit interval
$A$; the bounded-coordinate decomposition generates every element
from that interval.  For $b\in B$, write $\iota_j(b)\in A\subseteq G$
for the element equal to $b$ in factor $j$ and to $0$ elsewhere.

For $k\geq1$ set
$W_k=\operatorname{pr}_k^{-1}(\mathcal U_1)\subseteq A$,
and set $W_k=A$ for $k\leq0$.  Then
\[
 Q(\mathcal W)=\bigcap_{k\in\mathbb Z}\delta_k^{-1}(W_k)
\]
is a $\nu_A$-open zero-neighbourhood.
We show that no two $\nu_A$-zero-neighbourhoods have their
sum contained in $Q(\mathcal W)$.

Let $N,N'$ be any such neighbourhoods.
Choose coordinate boxes $Q(\mathcal V)\subseteq N$ and
$Q(\mathcal V')\subseteq N'$.
The product-open zero-neighbourhood $V_0\subseteq A$ contains
a basic cylinder restricting only finitely many factors.
Choose $j\geq1$ outside those factors.  Then
$\iota_j(c_n)\in V_0$ for every $n$.
For this fixed $j$, we have $\iota_j(b_n)\to0$ in $A$, so choose
$n$ sufficiently large that
\[
 \iota_j(b_n)\in\bigcap_{i=0}^{j-1}V'_i.
\]
Set $g=\iota_j(c_n)$ and $h=j\iota_j(b_n)$, where the latter is
an integer multiple in $G$.  Pointwise truncation gives
\[
 \delta_i(g)=
 \begin{cases}
  \iota_j(c_n),&i=0,\\
  0,&i\ne0,
 \end{cases}
 \qquad
 \delta_i(h)=
 \begin{cases}
  \iota_j(b_n),&0\leq i<j,\\
  0,&\text{otherwise}.
 \end{cases}
\]
Thus $g\in Q(\mathcal V)\subseteq N$ and
$h\in Q(\mathcal V')\subseteq N'$.
In factor $j$, the sum $g+h$ takes value $j+1$ at $(n,0)$,
value $j$ at $(n,1)$, and $0$ elsewhere.  Hence
\[
 \delta_j(g+h)=\iota_j(c_n)\notin W_j,
 \qquad g+h\notin Q(\mathcal W).
\]
Addition is therefore not jointly continuous at $(0,0)$ for
$\nu_A$.  By Proposition~\ref{prop-box-quasi-reconstruction},
\[
 \sigma_A=\gamma_A\subsetneq\nu_A.
\]
In this construction, $j$ can be chosen arbitrarily large.
The sum then satisfies $\delta_j(g+h)=\iota_j(c_n)\notin W_j$,
although $g$ and $h$ each have only finitely many nonzero deviations.
\end{exmp}

\begin{cor}[Comparison and strictness of the reconstruction topologies]
\label{cor-quasi-full-chain}
For every $(A,\tau)\in\TopMV$, the reconstruction topologies satisfy
\begin{equation}
\label{eq-quasi-full-chain}
 \kappa_A\subseteq\sigma_A=\gamma_A
 \subseteq\nu_A\subseteq\lambda_A
 \subseteq\sigma_A^{\mathrm q}.
\end{equation}
All induce $\tau$ on $A$, and all except possibly
$\sigma_A^{\mathrm q}$ induce $\tau_n$ on every bounded stage $I_n$.
All four inclusions can be strict simultaneously, even for a
Hausdorff metrizable input.
For the rational product cube
$A=(\mathbb Q\cap[0,1])^{\mathbb N}$ with its product topology,
\begin{equation}
\label{eq-quasi-full-chain-strict}
 \kappa_A\subsetneq\sigma_A=\gamma_A=\nu_A
 \subsetneq\lambda_A\subsetneq\sigma_A^{\mathrm q},
\end{equation}
while Example~\ref{exmp-box-not-group} gives
$\sigma_A\subsetneq\nu_A$.
For the usual interval $A=[0,1]$, all the topologies in
\eqref{eq-quasi-full-chain} coincide.
\end{cor}

\begin{proof}
Corollary~\ref{cor-reconstruction-topology-chain} and
Proposition~\ref{prop-box-quasi-reconstruction} give the comparisons
through $\lambda_A$, together with the bounded-stage restrictions.
Proposition~\ref{prop-lambda-quasi-reconstruction} gives the
remaining comparison \eqref{eq-lambda-quasi-comparison}, and
Proposition~\ref{prop-quasi-adjunction} gives
$\sigma_A^{\mathrm q}|_A=\tau$.
The rational-cube comparison and the simultaneous strictness of
all four inclusions are proved in
Example~\ref{exmp-all-reconstruction-inclusions-strict} below.
Example~\ref{exmp-box-not-group} gives $\sigma_A\subsetneq\nu_A$.

For $A=[0,1]$, Remark~\ref{rem-reconstruction-strictness} gives
equality through $\lambda_A$, hence also for $\nu_A$.
The translated-interval criterion
\eqref{eq-quasi-translated-final-test} makes
$\sigma_A^{\mathrm q}$ usual as well: every point of $\mathbb R$
lies in the usual interior of a translated unit interval, so a set
whose intersection with every such interval is relatively open is
usual-open.
\end{proof}

\begin{exmp}[Simultaneous strictness in the full comparison]
\label{exmp-all-reconstruction-inclusions-strict}
Let $A_1=(\mathbb Q\cap[0,1])^{\mathbb N}$ be the rational product
cube, and let $A_2=B^{\mathbb N_{\geq1}}$ be the Boolean product
from Example~\ref{exmp-box-not-group}.
With their product topologies, $C=A_1\times A_2$ is a Hausdorff
metrizable topological MV-algebra for which
\begin{equation}
\label{eq-all-reconstruction-inclusions-strict}
 \kappa_C\subsetneq\sigma_C=\gamma_C
 \subsetneq\nu_C\subsetneq\lambda_C
 \subsetneq\sigma_C^{\mathrm q}.
\end{equation}

We first prove $\nu_{A_1}=\sigma_{A_1}$.
For the real product cube $A_{\mathbb R}=[0,1]^{\mathbb N}$,
Corollary~\ref{cor-compact-reconstruction} gives
$\sigma_{A_{\mathbb R}}=\lambda_{A_{\mathbb R}}$.
Proposition~\ref{prop-box-quasi-reconstruction} therefore yields
$\nu_{A_{\mathbb R}}=\sigma_{A_{\mathbb R}}$.
The rational interval product $A_1$ has the subspace topology in
$A_{\mathbb R}$.  This remains true for their box powers:
a box of relatively open sets is the intersection of an ambient
open box with the smaller box power.
The rational Mundici group is an $\ell$-subgroup of the real one,
and its deviation embedding is the restriction of the real
deviation embedding.  Consequently, $\nu_{A_1}$ is the subspace
topology inherited from the real Mundici group with topology
$\nu_{A_{\mathbb R}}$.
Since the latter is an $\ell$-group topology, so is $\nu_{A_1}$.
It therefore equals its group reflection, which is $\sigma_{A_1}$
by Proposition~\ref{prop-box-quasi-reconstruction}.
Remark~\ref{rem-reconstruction-strictness} gives
$\kappa_{A_1}\subsetneq\sigma_{A_1}\subsetneq\lambda_{A_1}$, while
Example~\ref{exmp-quasi-stage-strict-left} gives
$\lambda_{A_1}\subsetneq\sigma_{A_1}^{\mathrm q}$.

We now pass to the product $C$.
Write $\Xi A_j=(G_j,u_j)$ for $j=1,2$.
The Mundici group of $C$ can be identified algebraically with
$(G_1\times G_2,(u_1,u_2))$, whose unit interval is $C$.
For each factor, let $p_j$ be the group projection and let $i_j$
insert an element in that factor and $0$ in the other.
We show that $p_j$ and $i_j$ are continuous for each of the
corresponding topologies
\[
 \kappa,\quad\lambda,\quad\sigma=\gamma,
 \quad\nu,\quad\sigma^{\mathrm q}.
\]
Although $i_j$ does not preserve the strong unit, its continuity
follows directly from these constructions.

For $\kappa$ and $\nu$, the bounded coordinates are computed
componentwise on the group product.
Projection and insertion of the constant coordinates of $0$ are
continuous for both the product and box coordinate topologies.
For $\lambda$, the bounded stages and their topologies satisfy
\[
 I_n^C=I_n^{A_1}\times I_n^{A_2},\qquad
 \tau_n^C=\tau_n^{A_1}\times\tau_n^{A_2},
\]
by the initial definition of the stage topologies.
Thus $p_j$ and $i_j$ are continuous on every bounded stage;
finality makes them continuous for the bounded-stage final
topologies.  For either map, composition with the continuous identity
from its target $\lambda$ topology to its target $\sigma$ topology
gives a continuous homomorphism into a topological group.
Since $\sigma=R_{\mathrm{grp}}(\lambda)$ by
Corollary~\ref{cor-reconstruction-topology-chain}, the universal
property of group reflection makes it continuous from its source
$\sigma$ topology as well.

For $\sigma^{\mathrm q}$, use the translated-interval criterion
\eqref{eq-quasi-translated-final-test}.
On a translated unit interval of $G_j$, the insertion has the form
$a\mapsto i_j(g)+i_j(a)$.
The interval map $a\mapsto i_j(a)$ is continuous into $C$ with
its product topology; inclusion into the reconstructed group
and translation by $i_j(g)$ are continuous.
Likewise, on a translated interval of the product group,
$p_j(g+a)=p_j(g)+p_j(a)$ is continuous because the interval
projection $C\longrightarrow A_j$ is continuous.
Finality therefore gives the two required continuous maps.

Since $p_ji_j=1_{G_j}$, each $i_j$ is a topological embedding
for all five constructions.  Equality of two consecutive
topologies on $C$ would therefore force equality of their
restrictions to both factors.
For $A_1$, the first, third, and fourth inclusions are strict
by the rational-cube comparison proved above; for $A_2$,
Example~\ref{exmp-box-not-group} gives strictness of the second.
Hence all four inclusions in
\eqref{eq-all-reconstruction-inclusions-strict} are strict.
\end{exmp}

Since $\nu_A$ already has jointly continuous lattice operations,
one might try to prove the topological representation theorem by
applying group reflection to it.  Such a proof would still have to
show that the coarser topology induces $\tau$ on $A$ and that the
lattice operations are continuous for that topology.
Joint lattice continuity for $\nu_A$ does not immediately imply
continuity after reflection, since the topology on $G_A\times G_A$
also becomes coarser.
The proof of Proposition~\ref{prop-box-quasi-reconstruction} uses
Section~\ref{sec-group-envelope} to identify the reflected topology
with $\sigma_A$, for which both properties have already been proved.
Thus the box description does not provide an independent proof
of the representation theorem.

\section{Summary}
\label{sec-summary}

For every topological MV-algebra $(A,\tau)$, we have constructed
the finest $\ell$-group topology $\sigma_A$ on its Mundici group
$(G_A,u_A)$ that induces $\tau$ on $A=[0,u_A]$
(Corollary~\ref{cor-reconstruction-topology-chain}).
Taking unit intervals and reconstructing with $\sigma_A$ give
the equivalence (Theorem~\ref{thm-topological-correspondence})
\[ 
 \TopMV\simeq\TopulgM.
\]
Here the group-side objects are Mundici-fine: their topology is
the finest $\ell$-group extension of their unit-interval topology.

An $\ell$-group extension need not be unique.  Any two extensions
of $\tau$ induce the same topology $\tau_n$ on each bounded stage
\[ 
 I_n=[-nu_A,nu_A],
\]
but they can differ on the whole group
(Proposition~\ref{prop-lgroup-extension-bounded-stages} and
Example~\ref{exmp-product-cubes}).

Every extension lies between the coordinate topology $\kappa_A$
and the bounded-stage final topology $\lambda_A$.
The coordinate topology is the coarsest topology making all bounded
coordinate maps into $A$ continuous.  The bounded-stage final
topology declares a set open exactly when its intersection with
each $I_n$ is open in $\tau_n$.
Neither need make addition jointly continuous
(Examples~\ref{exmp-kappa-not-group} and~\ref{exmp-lambda-not-group}).

The finest group topology contained in $\lambda_A$ is its group
reflection $\gamma_A$.
Proposition~\ref{prop-kappa-rho-gamma} and
Theorem~\ref{thm-rho-lgroup} show that this coarsening preserves
every bounded-stage topology and has jointly continuous lattice
operations.  It coincides with $\sigma_A$ by
Corollary~\ref{cor-reconstruction-topology-chain}.
Table~\ref{tab-reconstruction-topologies} records the construction
and role of each topology.

\begin{table}[!htb]
\centering
\small
\renewcommand{\arraystretch}{1.25}
\begin{tabular}{@{}>{\raggedright\arraybackslash}p{0.10\textwidth}
                  >{\raggedright\arraybackslash}p{0.43\textwidth}
                  >{\raggedright\arraybackslash}p{0.41\textwidth}@{}}
\hline
Topology & Construction & Role in reconstruction\\
\hline
$\kappa_A$ &
Product-subspace topology via $\pi$ or $\delta$;
equivalently, the initial topology of translated truncations
\eqref{eq-translated-initial-family}. &
Common lower bound for all $\ell$-group extensions
(Proposition~\ref{prop-kappa-lambda}); least exact quasi topology
(Proposition~\ref{thm-quasi-minimal-reconstruction}).\\
$\lambda_A$ &
Bounded-stage finalization of $\kappa_A$
(Definition~\ref{def-Mundici-final-topology}). &
Finest topology inducing $\tau_n$ on every $I_n$
(Proposition~\ref{prop-lambda-exact-space}); upper bound for all
admissible $\ell$-group topologies
(Theorem~\ref{thm-lambda-envelope}).\\
$\gamma_A$ &
Topological-group reflection of $\lambda_A$
\eqref{eq-gamma-A-def}, with a finite-sum coordinate-box base
(Proposition~\ref{prop-coordinate-box-topology}). &
Supplies joint group continuity; bounded-output and reverse-tail
control prove stage preservation and lattice continuity
(Proposition~\ref{prop-kappa-rho-gamma} and
Theorem~\ref{thm-rho-lgroup}).\\
$\sigma_A$ &
Largest admissible $\ell$-group topology
(Lemma~\ref{lem-finest-admissible}), called the Mundici topology. &
Equals $\gamma_A$ and is the largest $\ell$-group extension
(Corollary~\ref{cor-reconstruction-topology-chain}); yields
$\TopMV\simeq\TopulgM$
(Theorem~\ref{thm-topological-correspondence}).\\
$\nu_A$ &
Initial topology induced by either $\pi$ or $\delta$ into
$(A^{\mathbb Z},\tau_{\mathrm{box}})$. &
Explicit upper bound contained in $\lambda_A$, with joint lattice
continuity and group reflection $\sigma_A$
(Proposition~\ref{prop-box-quasi-reconstruction}).\\
$\sigma_A^{\mathrm q}$ &
Final topology of translated unit intervals
\eqref{eq-quasi-translated-final-test}. &
Largest quasi-admissible topology
(Proposition~\ref{prop-quasi-translated-final}); gives the quasi correspondence
(Corollary~\ref{thm-quasi-correspondence}) and recovers $\sigma_A$
by $\ell$-group reflection
(Proposition~\ref{thm-quasi-reflection-factorization}).\\
\hline
\end{tabular}
\caption{The reconstruction topologies for a topological MV-algebra.}
\label{tab-reconstruction-topologies}
\end{table}

Two further topologies describe the separate-continuity setting.
The coordinate-box topology $\nu_A$ has continuous inversion and
jointly continuous lattice operations, while its addition is separately
continuous but need not be jointly continuous
(Example~\ref{exmp-box-not-group}).
The topology $\sigma_A^{\mathrm q}$, which gives the quasitopological
Mundici correspondence (Corollary~\ref{thm-quasi-correspondence}), is the
finest topology inducing $\tau$ on $A$ for which inversion is continuous
and addition and the lattice operations are separately continuous.
For any $\ell$-group extension $\rho$ of $\tau$, all these topologies
fit into the comparison
(Corollary~\ref{cor-quasi-full-chain})
\[
 \kappa_A\subseteq\rho\subseteq\sigma_A=\gamma_A
 \subseteq\nu_A\subseteq\lambda_A
 \subseteq\sigma_A^{\mathrm q}.
\]
Every topology in this chain recovers $\tau$ on $A$.
All except possibly $\sigma_A^{\mathrm q}$ also agree on every
bounded stage $I_n$.  The topology $\sigma_A^{\mathrm q}$ can be
strictly finer on a bounded stage
(Example~\ref{exmp-quasi-stage-strict-left}).

None of the four inclusions between $\kappa_A$, $\sigma_A$,
$\nu_A$, $\lambda_A$, and $\sigma_A^{\mathrm q}$ can in general
be replaced by equality.
Table~\ref{tab-reconstruction-strictness} lists examples for each
strict inclusion.  Taking a finite product of these examples gives
a single Hausdorff metrizable topological MV-algebra for which
all four inclusions are strict simultaneously
(Example~\ref{exmp-all-reconstruction-inclusions-strict}).

\begin{table}[!htb]
\centering
\small
\renewcommand{\arraystretch}{1.25}
\begin{tabular}{@{}>{\raggedright\arraybackslash}p{0.3\textwidth}
                  >{\raggedright\arraybackslash}p{0.31\textwidth}
                  >{\raggedright\arraybackslash}p{0.13\textwidth}@{}}
\hline
Strict inclusion & Input with product topology & Reference\\
\hline
$\kappa_A\subsetneq\sigma_A$ &
$(\mathbb Q\cap[0,1])^{\mathbb N}$ &
Remark~\ref{rem-reconstruction-strictness}\\
$\sigma_A\subsetneq\nu_A$ &
$B^{\mathbb N_{\geq1}}$ &
Example~\ref{exmp-box-not-group}\\
$\nu_A\subsetneq\lambda_A$ &
$(\mathbb Q\cap[0,1])^{\mathbb N}$ &
Example~\ref{exmp-all-reconstruction-inclusions-strict}\\
$\lambda_A\subsetneq\sigma_A^{\mathrm q}$ &
$(\mathbb Q\cap[0,1])^{\mathbb N}$ &
Example~\ref{exmp-quasi-stage-strict-left}\\
All four inclusions simultaneously &
$(\mathbb Q\cap[0,1])^{\mathbb N}\times B^{\mathbb N_{\geq1}}$ &
Example~\ref{exmp-all-reconstruction-inclusions-strict}\\
\hline
\end{tabular}
\caption{Witnesses to strictness in
Corollary~\ref{cor-quasi-full-chain}.}
\label{tab-reconstruction-strictness}
\end{table}

Under additional hypotheses, the bounds $\kappa_A$ and $\lambda_A$
themselves become $\ell$-group topologies, with the extremal properties
listed in Table~\ref{tab-extremal-extensions}.
If $(A,\tau)$ is $k_\omega$, then $\sigma_A=\lambda_A$,
so the Mundici topology is the finest topology inducing
$\tau_n$ on every $I_n$
(Theorem~\ref{thm-komega-reconstruction} and
Proposition~\ref{prop-lambda-exact-space}).
If $(A,\tau)$ is locally convex, then $\kappa_A$ is the least
$\ell$-group extension
(Corollary~\ref{cor-extremal-lgroup-extensions}).
Even when both hypotheses hold, the least and largest extensions
can differ (Example~\ref{exmp-product-cubes}).

\begin{table}[!htb]
\centering
\small
\renewcommand{\arraystretch}{1.25}
\begin{tabular}{@{}>{\raggedright\arraybackslash}p{0.2\textwidth}
                  >{\raggedright\arraybackslash}p{0.15\textwidth}
                  >{\raggedright\arraybackslash}p{0.37\textwidth}@{}}
\hline
Hypothesis on $A$ & Topology & Extremal property\\
\hline
$k_\omega$ & $\lambda_A=\sigma_A$ &
Largest topology inducing $\tau_n$ on every $I_n$\\
Locally convex & $\kappa_A$ &
Least $\ell$-group extension\\
\hline
\end{tabular}
\caption{Extremal $\ell$-group extensions under additional hypotheses.}
\label{tab-extremal-extensions}
\end{table}

\section*{Acknowledgements}

The authors acknowledge support from the National Natural Science
Foundation of China (Grant Nos. 12371463 and 12671561).

\section*{Declaration of generative AI and AI-assisted technologies
in the manuscript preparation process}

During the preparation of this work, the authors used ChatGPT to assist with 
the review of mathematical arguments and citations, manuscript organization, 
language editing, consistency checking, and submission-readiness review. 
After using this tool, the authors reviewed and edited the content as needed 
and take full responsibility for the content of the article.

\bibliographystyle{abbrv}
\bibliography{lili}

\end{document}